\documentclass[preprint,3p]{elsarticle}
\journal{J. Comput. Phys.}

\usepackage{stmaryrd}
\usepackage{amsmath,amssymb}
\usepackage{tabularx,epsfig,color,tabularx}
\usepackage{epstopdf}
\usepackage{algorithm,algorithmic,multirow}
\usepackage{graphicx}
\usepackage{epsfig}
\usepackage{epstopdf}
\usepackage{caption}
\usepackage{subcaption}
\usepackage{mathtools}

\input psfig.sty

\renewcommand{\theequation}{\arabic{section}.\arabic{equation}}

\newcommand{\bx}{{\bf x} }

\newcommand{\bB}{{\bf B} }
\newcommand{\bk}{{\bf k} }
\newcommand{\bJ}{{\bf J} }
\newcommand{\bee}{{\bf e} }
\newcommand{\beU}{{ U} }

\newcommand{\eps}{\varepsilon}

\newtheorem{thm}{Theorem}[section]
\newtheorem{lemma}{Lemma}
\newtheorem{rmk}{Remark}[section]

\newcommand{\be}{\begin{equation}}
\newcommand{\ee}{\end{equation}}
\newcommand{\ba}{\begin{array}}
\newcommand{\ea}{\end{array}}
\newcommand{\bea}{\begin{eqnarray}}
\newcommand{\eea}{\end{eqnarray}}
\newcommand{\beas}{\begin{eqnarray*}}
\newcommand{\eeas}{\end{eqnarray*}}

\newcommand{\vep}{\varepsilon}

\begin{document}

\begin{frontmatter}

\title{Numerical methods and comparison for the Dirac equation in the nonrelativistic limit regime}
\author[1]{Weizhu Bao}
 \ead{matbaowz@nus.edu.sg}

\address[1]{Department of Mathematics,  National University of
Singapore, Singapore 119076, Singapore}
\ead[url]{http://www.math.nus.edu.sg/\~{}bao/}

\author[2,3]{Yongyong Cai\corref{5}} 
\ead{yongyong.cai@gmail.com}
\address[2]{Beijing Computational Science Research Center, Beijing 100094,
P. R. China}
\address[3]{Department of Mathematics, Purdue University, West Lafayette,
IN 47907, USA }
\cortext[5]{Corresponding author.}

\author[1]{Xiaowei Jia}
\ead{A0068124@nus.edu.sg}

\author[4]{Qinglin Tang}
\ead{tqltql2010@gmail.com}
\address[4]{Institut Elie Cartan de Lorraine, Universit\'e de Lorraine, Inria Nancy-Grand Est,
F-54506 Vandoeuvre-l\`es-Nancy Cedex, France}




\begin{abstract}
We analyze rigorously error estimates and compare numerically spatial/temporal
resolution of various numerical methods for the discretization of the Dirac equation
in the nonrelativistic limit regime, involving a small dimensionless parameter $0<\varepsilon\ll 1$
which is inversely proportional to the speed of light.  In this limit
regime, the solution is highly oscillatory in time, i.e. there are propagating waves
with wavelength $O(\varepsilon^2)$ and $O(1)$ in time and space, respectively.
We begin with several frequently used finite difference time domain (FDTD)
methods and obtain rigorously their error estimates in the
nonrelativistic limit regime by paying particular attention to
how error bounds depend explicitly on mesh size $h$ and time step $\tau$ as well
as the small parameter $\varepsilon$. Based on the error bounds, in order
to obtain `correct' numerical solutions in the nonrelativistic limit regime,
i.e. $0<\varepsilon\ll 1$,
the FDTD methods share the same $\varepsilon$-scalability on time step and mesh size as:
$\tau=O(\varepsilon^3)$ and $h=O(\sqrt{\varepsilon})$. Then we propose and analyze two numerical methods
for the discretization of the Dirac equation by using the Fourier
spectral discretization for spatial derivatives combined
with the exponential wave
integrator and time-splitting technique
for temporal derivatives, respectively.
Rigorous error bounds for the two numerical methods
show that their $\varepsilon$-scalability is improved
to $\tau=O(\varepsilon^2)$ and $h=O(1)$ when $0<\varepsilon\ll 1$.
Extensive numerical results
are reported to support our error estimates.
\end{abstract}


\begin{keyword}
Dirac equation, nonrelativistic limit regime, finite difference time domain method,
exponential wave integrator spectral method, time splitting spectral method,
$\varepsilon$-scalability
\end{keyword}

\end{frontmatter}

\section{Introduction}\setcounter{equation}{0}
The Dirac equation, which plays an important role in particle physics, is a relativistic wave equation derived
by the British physicist Paul Dirac in 1928 \cite{Dirac1,Dirac2,Dirac3,Thaller}. It provided a description of
elementary spin-$1/2$ massive particles, such as electrons and positrons,
consistent with both the principle of quantum
mechanics and the theory of special relativity.
It was the first theory to fully account for relativity in the context
of quantum mechanics.  It addressed the fine details of the hydrogen spectrum in a
completely rigorous way and predicted the existence of a new form of matter, antimatter \cite{Ad}.
Since the graphene was first produced in the lab in 2003
\cite{AMPGMKWTNLG, NGMJKGDF, NGMJZDGF, NJZMSZMBKG, SGMHBKN},
the Dirac equation has been extensively adopted to study theoretically the structures and/or  dynamical
properties of graphene and graphite as well as two dimensional (2D) materials \cite{NGPNG}.
This experimental advance renewed extensively
the research interests on the mathematical analysis and numerical simulations
of the Dirac equation and/or the (nonlinear) Schr\"{o}dinger equation without/with
external potentials, especially the honeycomb lattice potential \cite{AZ, FW}.

We consider the three dimensional (3D) Dirac equation for describing the time evolution of spin-$1/2$ massive
particles, such as electrons and positrons, within external time-dependent electromagnetic potentials \cite{Dirac1,Dirac2}
\begin{equation}
\label{3DME}
i\hbar\partial_t\Psi(t,\bx) = \Bigl[-ic\hbar\sum_{j=1}^3\alpha_j\partial_j
+mc^2\beta \Bigr] \Psi(t,\bx)+e\Bigl[V(t,\bold{x})I_4 -\sum_{j=1}^3A_j(t,\mathbf{x})\alpha_j\Bigr]\Psi(t,\bx).
\end{equation}
Here, $i=\sqrt{-1}$,  $t$ is time, $\bx=(x_1,x_2,x_3)^T\in {\mathbb R}^3$ (equivalently written as $\bx=(x,y,z)^T$)
is the spatial coordinate vector, $\partial_k=\frac{\partial}{\partial x_k}$ ($k=1,2,3$),
 $\Psi :=\Psi(t,\bx)=(\psi_1(t,\bx),\psi_2(t,\bx), \psi_3(t,\bx), \psi_4(t,\bx))^T\in\mathbb{C}^4$
 is the complex-valued vector
wave function of the ``spinorfield''.
$I_n$ is the $n\times n$ identity matrix for $n\in {\mathbb N}$,
$V:=V(t,\bx)$ is the real-valued electrical potential and
${\bf A}:={\bf A}(t,\bx)=(A_1(t,\bx), A_2(t,\bx), A_3(t,\bx))^T$ is the real-valued magnetic potential vector,
and hence the electric field is given by ${\bf E}(t,\bx)=-\nabla V-\partial_t {\bf A}$ and
the magnetic field is given by ${\bf B}(t,\bx)={\rm curl}\, {\bf A}=\nabla \times {\bf A}$.
The physical constants are: $c$ for the speed of light, $m$ for the particle's rest mass,
$\hbar$ for the Planck constant and $e$ for the unit charge. In addition, the $4\times 4$
matrices $\alpha_1$, $\alpha_2$, $\alpha_3$ and $\beta$ are defined as
\be \label{alpha}
\alpha_1=\left(\begin{array}{cc}
\mathbf{0} & \sigma_1  \\
\sigma_1 & \mathbf{0}  \\
\end{array}
\right),\qquad
\alpha_2=\left(\begin{array}{cc}
\mathbf{0} & \sigma_2 \\
\sigma_2 & \mathbf{0} \\
\end{array}
\right), \qquad
\alpha_3=\left(\begin{array}{cc}
\mathbf{0} & \sigma_3 \\
\sigma_3 & \mathbf{0} \\
\end{array}
\right),\qquad
\beta=\left(\begin{array}{cc}
I_{2}& \mathbf{0} \\
\mathbf{0} & -I_{2} \\
\end{array}
\right),
\ee
with  $\sigma_1$, $\sigma_2$,  $\sigma_3$
 (equivalently written $\sigma_x$, $\sigma_y$, $\sigma_z$)
being the Pauli matrices  defined as
\be\label{Paulim}
\sigma_{1}=\left(
\begin{array}{cc}
0 & 1  \\
1 & 0  \\
\end{array}
\right), \qquad
\sigma_{2}=\left(
\begin{array}{cc}
0 & -i \\
i & 0 \\
\end{array}
\right),\qquad
\sigma_{3}=\left(
\begin{array}{cc}
1 & 0 \\
0 & -1 \\
\end{array}
\right).
\ee

In order to scale the Dirac equation (\ref{3DME}), we introduce
\be\label{scale}
\tilde{t}=\frac{t}{t_s}, \quad \tilde\bx=\frac{\bx}{x_s},
\quad \tilde\Psi(\tilde t, \tilde \bx)=x_s^{3/2}\,\Psi(t,\bx), \quad
\tilde V(\tilde t, \tilde \bx)=\frac{V(t,\bx)}{A_s},\quad  \tilde A_j(\tilde t, \tilde \bx)
=\frac{A_{j}(t, \bx)}{A_s},\ \ j = 1, 2, 3,
\ee
where $x_s$, $t_s$ and $A_s$ are the dimensionless length unit, time unit and potential unit, respectively,
satisfying $t_s=\frac{mx_s^2}{\hbar}$ and $A_s=\frac{m v^2}{e}$ with $v=\frac{x_s}{t_s}$ being the wave speed.
Plugging (\ref{scale}) into (\ref{3DME}), multiplying by $\frac{t_sx_s^{3/2}}{\hbar}$,
and then removing all $\tilde{ }$, we obtain
the following dimensionless Dirac equation in 3D
\be
\label{SDE}
i\partial_t\Psi(t,\bx)=\Bigl[-\frac{i}{\varepsilon}\sum_{j=1}^{3}\alpha_j
\partial_j+\frac{1}{\varepsilon^2}\beta\Bigr]\Psi(t,\bx)+
\Bigl[V(t,\bx)I_4-\sum_{j=1}^{3}A_j(t,\bx)\alpha_j\Bigr]\Psi(t,\bx), \qquad \bx\in{\mathbb R}^3,
\ee
where $\varepsilon$ is a dimensionless parameter  inversely proportional to the speed of light given by
\begin{equation}\label{eps}
0<\varepsilon := \frac{x_s}{t_s\,c}=\frac{v}{c}\le 1.
\end{equation}

We remark here that if one chooses the dimensionless length unit
$x_s=\frac{\hbar}{mc}$, $ t_s=\frac{x_s}{c}$ and $A_s=\frac{mc^2}{e}$
in (\ref{scale}),  then $\varepsilon =1$ in (\ref{eps}) and Eq. (\ref{SDE}) with  $\varepsilon =1$
takes the form often appearing in the 
literature \cite{Abe, BL, bolj, BHM, Esteban, FS, Hamm, HJMSZ}. This choice of $x_s$ is appropriate
when the wave speed is at the same order of the speed of light. However, when the wave speed is
much smaller than the speed of light, a different choice of $x_s$ is more appropriate. Note that
the choice of $x_s$ determines the observation scale of the time evolution of the particles and decides:
(i) which phenomena are `visible' by asymptotic analysis, and (ii) which phenomena
can be resolved by discretization by specified spatial/temporal grids.
In fact, there are two important parameter regimes: One is  $\varepsilon =1$
($\Longleftrightarrow x_s=\frac{\hbar}{mc}$, $t_s=\frac{x_s}{c}$ and $A_s=\frac{mc^2}{e}$), then Eq. (\ref{SDE}) describes the case
that wave speed is at the same order of the speed of light; the other one is $0<\varepsilon \ll 1$,
then Eq. (\ref{SDE}) is in the nonrelativistic limit regime.

Similarly to the dimension reduction of the nonlinear Schr\"{o}dinger equation
and/or the Schr\"{o}dinger-Poisson equations with/without anisotropic external potentials
\cite{BC1},  when  the initial data $\Psi(0,\bx)$ and the electromagnetic potentials
$V(t,\bx)$ and ${\bf A}(t,\bx)$ are independent of $z$ and thus
the wave function $\Psi$ is formally assumed  to be independent of $z$, or when the electromagnetic potentials
$V(t,\bx)$ and ${\bf A}(t,\bx)$ are strongly confined in the $z$-direction and thus
 $\Psi$ is formally assumed to be concentrated on the $xy$-plane,  then the 3D Dirac equation
(\ref{SDE}) can be reduced to the Dirac equation in 2D with $\bx=(x,y)^T\in{\mathbb R}^2$ as
\be
\label{SDE2d}
i\partial_t\Psi(t,\bx)=\Bigl[-\frac{i}{\varepsilon}\sum_{j=1}^{2}\alpha_j
\partial_j+\frac{1}{\varepsilon^2}\beta\Bigr]\Psi(t,\bx)+
\Bigl[V(t,\bx)I_4-\sum_{j=1}^{2}A_j(t,\bx)\alpha_j\Bigr]\Psi(t,\bx), \quad \bx\in{\mathbb R}^2.
\ee
This 2D Dirac equation has been widely used to model
the electron structure and/or dynamical properties of graphene since
they share the same dispersion relation on the Dirac points
\cite{AMPGMKWTNLG, NGMJKGDF, NGMJZDGF, NJZMSZMBKG, SGMHBKN,FW,FW2,FGNMPN}.
Similarly, under the proper assumptions on the initial data and the external electromagnetic
potential, the 3D Dirac equation
(\ref{SDE}) can be reduced to the Dirac equation in 1D with $\Psi=\Psi(t,x)$ as
\be
\label{SDE1d}
i\partial_t\Psi(t,x)=\Bigl[-\frac{i}{\varepsilon}\alpha_1
\partial_x+\frac{1}{\varepsilon^2}\beta\Bigr]\Psi(t,x)+\Bigl[V(t,x)I_4-A_1(t,x)\alpha_1\Bigr]\Psi(t,x), \qquad x\in{\mathbb R}.
\ee
In fact, the Dirac equation in 3D (\ref{SDE}), in 2D (\ref{SDE2d}) and in 1D (\ref{SDE1d}) can
be written in a unified way in $d$-dimensions ($d=1,2,3$)
\be
\label{SDEd}
i\partial_t\Psi(t,\bx)=\Bigl[-\frac{i}{\varepsilon}\sum_{j=1}^{d}\alpha_j
\partial_j+\frac{1}{\varepsilon^2}\beta\Bigr]\Psi(t,\bx)
+\Bigl[V(t,\bx)I_4-\sum_{j=1}^{d}A_j(t,\bx)\alpha_j\Bigr]\Psi(t,\bx), \qquad \bx\in{\mathbb R}^d,
\ee
and the initial condition for dynamics is given as
\be
\Psi(t=0,\bx)=\Psi_0(\bx), \qquad \bx\in{\mathbb R}^d.
\ee
The Dirac equation (\ref{SDEd}) is dispersive and time symmetric.
Introducing the position density $\rho_j$ for the $j$-component ($j=1,2,3,4$)
and the total density $\rho$ as well as the current density $\bJ(t,\bx)=(J_1(t,\bx),J_2(t,\bx)$, $J_3(t,\bx))^T$
\begin{equation} \label{obser11}
\rho(t,\bx)=\sum_{j=1}^4\rho_j(t,\bx)=\Psi^*\Psi,
\quad \rho_j(t,\bx)=|\psi_j(t,\bx)|^2, \quad 1\leq j\leq 4;\quad
J_l(t,\bx)=\frac{1}{\varepsilon}\Psi^*\alpha_l\Psi,\quad l=1,2,3,
\end{equation}
where $\overline{f}$ denotes the complex conjugate of  $f$ and
$\Psi^*=\overline{\Psi}^T$, then the following conservation law can be obtained from the Dirac equation
(\ref{SDEd})
\begin{equation} \label{cons}
\partial_t\rho(t,\bx)+\nabla\cdot \bJ(t,\bx) = 0, \qquad \bx\in{\mathbb R}^d, \quad t\ge0.
\end{equation}
Thus the Dirac equation (\ref{SDEd}) conserves the total mass as
\be\label{norm}
\|\Psi(t,\cdot)\|^2:=\int_{{\mathbb R}^d}|\Psi(t,\bx)|^2\,d\bx=
\int_{{\mathbb R}^d}\sum_{j=1}^4|\psi_j(t,\bx)|^2\,d\bx\equiv \|\Psi(0,\cdot)\|^2
=\|\Psi_0\|^2, \qquad t\ge0.
\ee
If the electric potential $V$ is perturbed by a real constant $V^0$, e.g. $V(t,\bx)\to V(t,\bx)+V^0$,
then the solution $\Psi(t,\bx)\to e^{-iV^0t}\Psi(t,\bx)$  which
implies the density of each component $\rho_j$ ($j=1,2,3,4$) and the total
density $\rho$ unchanged. When $d=1$, if the magnetic potential $A_1$
is perturbed by a real constant $A_1^0$, e.g. $A_1(t,\bx)\to A_1(t,\bx)+A_1^0$,
then the solution $\Psi(t,\bx)\to e^{iA_1^0t\alpha_1}\Psi(t,\bx)$  which
implies the total density $\rho$ unchanged; but this property is not valid when $d=2,3$.
In addition, when the electromagnetic potentials are time-independent,
i.e. $V(t,\bx)=V(\bx)$ and $A_j(t,\bx)=A_j(\bx)$ for $j=1,2,3$,
 the following energy functional is also conserved
\be\label{engery60}
E(t):=\int_{\mathbb{R}^d}\left(-\frac{i}{\varepsilon}
\sum_{j=1}^d\Psi^*\alpha_j\partial_j\Psi+\frac{1}{\varepsilon^2}\Psi^*\beta\Psi+V(\bx)|\Psi|^2-
\sum_{j=1}^dA_j(\bx)\Psi^*\alpha_j\Psi\right)d\bx\equiv E(0),\quad t\ge0.
\ee
Furthermore,  if the external electromagnetic potentials are constants, i.e.
$V(t,\bx)\equiv V^0$ and $A_j(t,\bx)\equiv A_j^0$ for $j=1,2,3$ with
${\bf A}^0=(A_1^0,\ldots,A_d^0)^T$,
the Dirac equation (\ref{SDEd}) admits the plane wave solution as
$\Psi(t,\bx)={\bf B}\,e^{i(\bk\cdot\bx-\omega t)}$,
where the time frequency $\omega$, amplitude vector $\bB\in {\mathbb R}^4$
and spatial wave number $\bk=(k_1,\ldots,k_d)^T\in {\mathbb R}^d$ satisfy the following
{\sl dispersion relation}
\be\label{disp}
\omega \bB=\Bigl[\sum_{j=1}^{d}\left(\frac{k_j}{\varepsilon}-A_j^0\right)\alpha_j
+\frac{1}{\varepsilon^2}\beta
+V^0I_{4}\Bigr]\bB,
\ee
which immediately implies the dispersion relation of the Dirac equation (\ref{SDEd}) as
\be\label{disper}
\omega:=\omega(\bk)=V^0\pm \frac{1}{\vep^2}\sqrt{1+\vep^2\left|\bk-\vep {\bf A}^0\right|^2},\qquad
\bk\in {\mathbb R}^d.
\ee

 Plugging (\ref{alpha}) and (\ref{Paulim}) into (\ref{SDE2d}), the 2D Dirac equation (\ref{SDE2d})
can be decoupled as
\be
\label{comp1d2}
\begin{split}
&i\partial_t\psi_1=-\frac{i}{\varepsilon}\left(\partial_x-i\partial_y\right)\psi_4
+\frac{1}{\varepsilon^2}\psi_1+V(t,\bx)\psi_1-\left[A_1(t,\bx)-iA_2(t,\bx)\right]\psi_4,\\
&i\partial_t\psi_4 =-\frac{i}{\varepsilon}\left(\partial_x+i\partial_y\right)\psi_1
-\frac{1}{\varepsilon^2}\psi_4
+V(t,\bx)\psi_4-\left[A_1(t,\bx)+iA_2(t,\bx)\right]\psi_1, \qquad \bx\in{\mathbb R}^2,
\end{split}
\ee
\be \label{comp1d3}
\begin{split}
&i\partial_t\psi_2=-\frac{i}{\varepsilon}\left(\partial_x+i\partial_y\right)\psi_3
+\frac{1}{\varepsilon^2}\psi_2+V(t,\bx)\psi_2-\left[A_1(t,\bx)+iA_2(t,\bx)\right]\psi_3,\\
&i\partial_t\psi_3 =-\frac{i}{\varepsilon}\left(\partial_x-i\partial_y\right)\psi_2
-\frac{1}{\varepsilon^2}\psi_3
+V(t,\bx)\psi_3-\left[A_1(t,\bx)-iA_2(t,\bx)\right]\psi_2,\qquad \bx\in{\mathbb R}^2.
\end{split}
\ee
Eq. (\ref{comp1d3}) will collapse to (\ref{comp1d2}) under the transformation
$y\to -y$ and $A_2\to -A_2$. Thus, in 2D, the Dirac equation (\ref{SDE2d})
can be reduced to the following simplified PDEs with $\Phi=\Phi(t,\bx)=(\phi_1(t,\bx),\phi_2(t,\bx))^T
\in{\mathbb C}^2$
\begin{equation}
\label{SDED2}
i\partial_t\Phi(t,\bx)=\Bigl[-\frac{i}{\varepsilon}\left(\sigma_1\partial_x+\sigma_2\partial_y\right)
+\frac{1}{\varepsilon^2}\sigma_3\Bigr]\Phi(t,\bx) +
\Bigl[V(t,\bx)I_2-A_1(t,\bx)\sigma_1-A_2(t,\bx)\sigma_2\Bigr]\Phi(t,\bx),\quad \bx\in{\mathbb R}^2,
\end{equation}
where $\Phi=(\psi_1,\psi_4)^T$ (or $\Phi=(\psi_2,\psi_3)^T$ under the transformation
$y\to -y$ and $A_2\to -A_2$). Similarly, in 1D, the Dirac equation (\ref{SDE1d})
can be reduced to the following simplified PDEs with $\Phi=\Phi(t,x)=(\phi_1(t,x),\phi_2(t,x))^T$
\begin{equation}
\label{SDED21d}
i\partial_t\Phi(t,x)=\Bigl[-\frac{i}{\varepsilon}\sigma_1\partial_x
+\frac{1}{\varepsilon^2}\sigma_3\Bigr]\Phi(t,x) +\Bigl[V(t,x)I_2-A_1(t,x)\sigma_1\Bigr]\Phi(t,x),
\qquad x\in{\mathbb R},
\end{equation}
where $\Phi=(\psi_1,\psi_4)^T$ (or $\Phi=(\psi_2,\psi_3)^T$).
Again,  the Dirac equation  in 2D (\ref{SDED2}) and in 1D (\ref{SDED21d}) can
be written in a unified way in $d$-dimensions ($d=1,2$)
\be
\label{SDEdd}
i\partial_t\Phi(t,\bx)=\Bigl[-\frac{i}{\varepsilon}\sum_{j=1}^{d}\sigma_j
\partial_j+\frac{1}{\varepsilon^2}\sigma_3\Bigr]\Phi(t,\bx)
+\Bigl[V(t,\bx)I_2-\sum_{j=1}^{d}A_j(t,\bx)\sigma_j\Bigr]\Phi(t,\bx), \quad \bx\in{\mathbb R}^d,
\ee
and the initial condition for dynamics is given as
\be\label{SDEini}
\Phi(t=0,\bx)=\Phi_0(\bx), \qquad \bx\in{\mathbb R}^d.
\ee
The Dirac equation (\ref{SDEdd}) is dispersive and time symmetric.
By introducing the position density $\rho_j$ for the $j$-th component ($j=1,2$)
and the total density $\rho$ as well as the current density $\bJ(t,\bx)=(J_1(t,\bx),J_2(t,\bx))^T$
\be\label{obser12}
\rho(t,\bx)=\sum_{j=1}^2\rho_j(t,\bx)=\Phi^\ast \Phi,
\quad \rho_j(t,\bx)=|\phi_j(t,\bx)|^2, \quad
J_j(t,\bx)=\frac{1}{\varepsilon}\Phi^*\sigma_j\Phi, \quad j=1,2,
\ee
 the conservation law (\ref{cons})
is also satisfied \cite{BHM}.  In addition, the Dirac equation (\ref{SDEdd})
conserves the total mass as
\be\label{normd}
\|\Phi(t,\cdot)\|^2:=\int_{{\mathbb R}^d}|\Phi(t,\bx)|^2\,d\bx=\int_{{\mathbb R}^d}\sum_{j=1}^2|\phi_j(t,\bx)|^2\,d\bx\equiv \|\Phi(0,\cdot)\|^2=\|\Phi_0\|^2, \qquad t\ge0.
\ee
Again, if the electric potential $V$ is perturbed by a real constant $V^0$, e.g. $V(t,\bx)\to V(t,\bx)+V^0$,
the solution $\Phi(t,\bx)\to e^{-iV^0t}\Phi(t,\bx)$  which
implies the density of each component $\rho_j$ ($j=1,2$) and the total
density $\rho$ unchanged. When $d=1$, if the magnetic potential $A_1$
is perturbed by a real constant $A_1^0$,  e.g. $A_1(t,\bx)\to A_1(t,\bx)+A_1^0$,
the solution $\Phi(t,\bx)\to e^{iA_1^0t\sigma_1}\Phi(t,\bx)$
implying the total density $\rho$ unchanged; but this property is not valid when $d=2$.
When the electromagnetic potentials are time-independent,
i.e. $V(t,\bx)=V(\bx)$ and $A_j(t,\bx)=A_j(\bx)$ for $j=1,2$,
 the following energy functional is also conserved
\be \label{engery65}
E(t):=\int_{\mathbb{R}^d}\left(-\frac{i}{\varepsilon}
\sum_{j=1}^d\Phi^*\sigma_j\partial_j\Phi+\frac{1}{\varepsilon^2}\Phi^*\sigma_3\Phi+V(\bx)|\Phi|^2-
\sum_{j=1}^dA_j(\bx)\Phi^*\sigma_j\Phi\right)d\bx\equiv E(0),\quad t\ge0.
\ee
Furthermore,  if the external electromagnetic potentials are constants, i.e.
$V(t,\bx)\equiv V^0$ and $A_j(t,\bx)\equiv A_j^0$ for $j=1,2$,
the Dirac equation (\ref{SDEdd}) admits the plane wave solution as
$\Phi(t,\bx)={\bf B}\,e^{i(\bk\cdot\bx-\omega t)}$,
where the time frequency $\omega$, amplitude vector $\bB\in {\mathbb R}^2$
and spatial wave number $\bk=(k_1,\ldots,k_d)^T\in {\mathbb R}^d$ satisfy the following
{\sl dispersion relation}
\be\label{dispp}
\omega \bB=\Bigl[\ \sum_{j=1}^{d}\left(\frac{k_j}{\varepsilon}-A_j^0\right)\sigma_j
+\frac{1}{\varepsilon^2}\sigma_3
+V^0I_{2}\Bigr]\bB,
\ee
which again implies the dispersion relation (\ref{disper})
of the Dirac equation (\ref{SDEdd}) for $d=2,1$.

\begin{figure}[htb]
\centerline{\psfig{figure=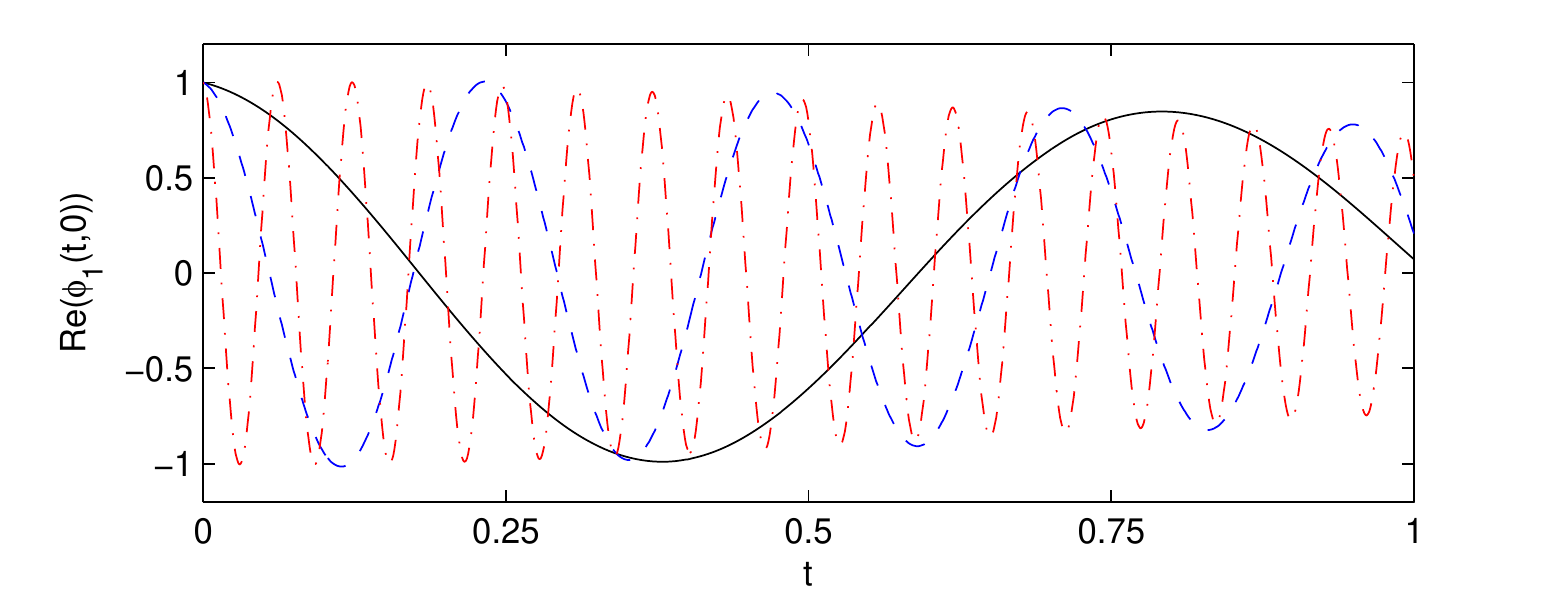,height=5cm,width=13cm,angle=0}}
\centerline{\psfig{figure=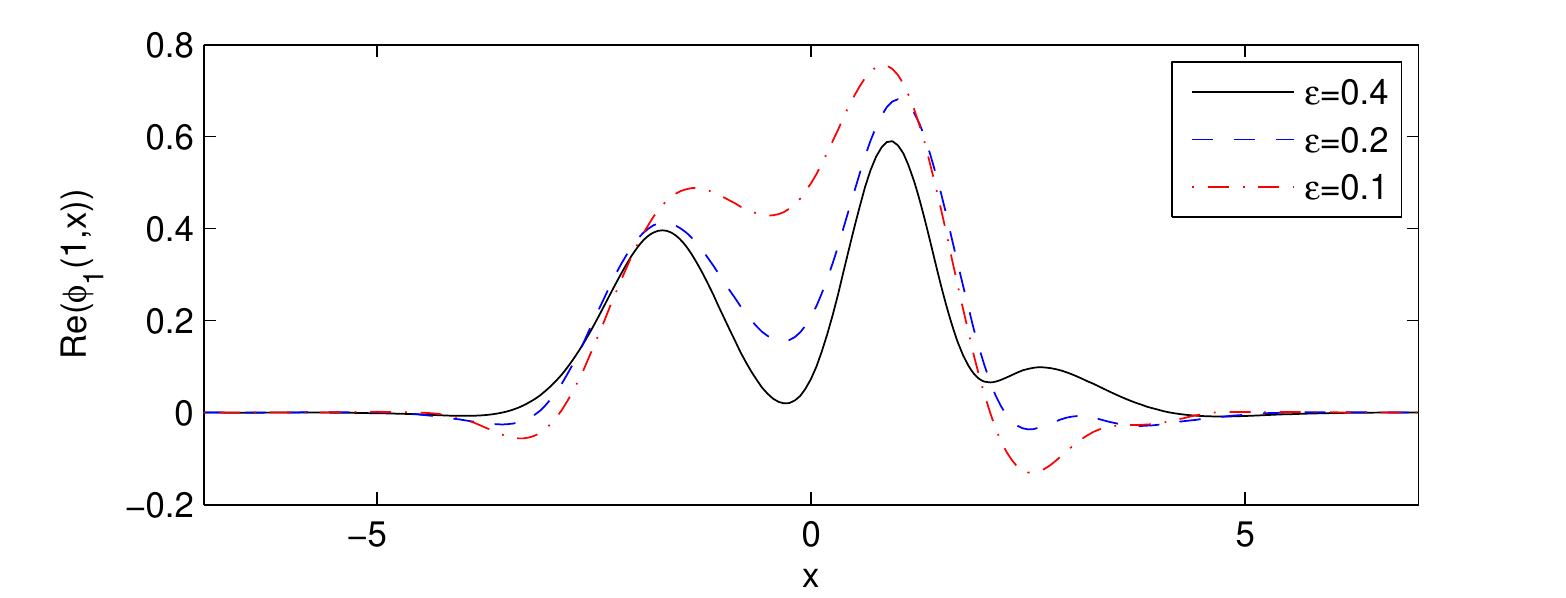,height=5cm,width=13cm,angle=0}}
\caption{The solution $\phi_1(t=1,x)$ and $\phi_1(t,x=0)$ of the Dirac equation
(\ref{SDEdd}) with $d=1$ for different $\varepsilon$. ${\rm Re}(f)$ denotes the real part
of $f$.}
\label{Oscillation}
\end{figure}

For the Dirac equation (\ref{SDEd}) with $\varepsilon=1$, i.e. $O(1)$-speed of light regime,
there are extensive analytical and numerical results
in the literatures. For the existence and multiplicity of bound states and/or standing wave solutions,
we refer to \cite{Esteban0,Das,Das2,Dol,Ges,Ves} and references therein.
For the analysis of the classical/semiclassical limits via
the Wigner transform techniques, we refer to \cite{GMMP,AS,boke,Bo,Kamm,spohn,ST4} and references therein.
For the numerical methods and comparison such as the finite difference time domain (FDTD) methods and
the Gaussian beam methods, we refer to \cite{Ant,AS,XST, WHJY, WT,Deuflhard,S,WT,Gosse} and references therein.
However, for the Dirac equation (\ref{SDEd}) with $0<\varepsilon\ll 1$,
i.e. nonrelativistic limit regime (or the scaled speed of light goes to infinity),
the analysis and efficient computation of the Dirac equation (\ref{SDEd}) (or (\ref{SDEdd}))
are mathematically rather complicated.
The main difficulty is due to that the solution is highly oscillatory in time
and the corresponding energy functionals (\ref{engery60}) and (\ref{engery65})
are indefinite \cite{BMP, Esteban} and become unbounded when $\varepsilon\to0$.
There are extensive mathematical analysis of the (semi)-nonrelativistic
limit of the Dirac equation (\ref{SDEd}) to the Pauli equation \cite{Hun,BMS2, BMP,CC,Foldy,Gri,NM,Mau,N,Sch,White}
and/or the Schr\"{o}dinger equation when  $\varepsilon\to 0$ \cite{BMP}.
These rigorous analytical results show that the solution
propagates waves with wavelength $O(\varepsilon^2)$ and
$O(1)$ in time and space, respectively, when $0<\varepsilon\ll 1$.
In fact, the oscillatory structure of the solution of the Dirac equation (\ref{SDEd})
when $0<\varepsilon\ll 1$ can be formally observed from its dispersion relation
(\ref{disp}) (or (\ref{dispp})). To illustrate this further, Figure \ref{Oscillation}
shows the solution of the Dirac equation (\ref{SDEdd}) with $d=1$,
$V(t,x)=\frac{1-x}{1+x^2}$, $A_1(t,x)=\frac{(1+x)^2}{1+x^2}$ and $\Phi_0(x)=\left(\exp(-x^2/2),\exp(-(x-1)^2/2)\right)^T$ for different $\varepsilon$.
This highly oscillatory nature of the solution of (\ref{SDEd})  (or
(\ref{SDEdd})) causes severe numerical burdens in practical computation,
making the numerical approximation of (\ref{SDEd})  (or
(\ref{SDEdd})) extremely challenging and costly in the
nonrelativistic regime $0<\varepsilon\ll 1$.
In \cite{HJMSZ}, the resolution of the time-splitting
Fourier pseudospectral (TSFP) method was studied
for the Maxwell-Dirac equation in the nonrelativistic limit regime.

Recently,
different numerical methods were proposed and analyzed for the efficient
computation of the Klein-Gordon equation in the nonrelativistic limit regime \cite{BCZ,BD,FS2}
and/or highly oscillatory dispersive partial differential equations (PDEs) \cite{BC2,BC3,BDZ,BDZ2}.
To our knowledge, so far there are few results on the
numerics of the Dirac equation in the nonrelativistic limit regime.
The aim of this paper is to study the efficiency of the frequently used  FDTD
and TSFP methods applied
to the Dirac equation in the nonrelativistic limit regime, to propose the
exponential wave integrator Fourier pseudospectral (EWI-FP) method and
to compare their resolution capacities in this regime. We
start with the detailed analysis on the stability and convergence of
several standard implicit/semi-implicit/explicit FDTD methods \cite{D}. Here we pay particular attention to
how the error bounds depend explicitly on the small parameter $\eps$
in addition to the mesh size $h$ and time step $\tau$. Based on the
estimates, in order to obtain `correct' numerical approximations
when $0<\eps\ll1$, the meshing strategy requirement
($\eps$-scalability) for those frequently used FDTD methods is:
$h = O(\sqrt{\eps})$ and $\tau = O(\eps^3)$,
which suggests that the standard FDTD methods are
computationally expensive for the Dirac equation (\ref{SDEd}) as $0<\eps\ll 1$. To
relax the $\eps$-scalability, we then propose the EWI-FP method and compare
it with the TSFP method,  whose $\eps$-scalability are optimal for both time and space in view
of the inherent oscillatory nature. The key ideas of the EWI-FP
are: (i) to apply the Fourier pseudospectral  discretization for spatial derivatives;
and (ii) to adopt the exponential wave integrator (EWI) for integrating
the ordinary differential equations (ODEs)
in phase space \cite{Gautschi0,Gautschi-type-3}
which was well demonstrated in the literatures that it has favorable properties
compared to standard time integrators for oscillatory  differential
equations \cite{Gautschi0,Gautschi-type-3,Gautschi-type-2,Gautschi-type-5}.
Rigorous error
estimates show that the $\eps$-scalability of the EWI-FP method  is
$h = O(1)$, and $\tau= O(\eps^2)$
for the Dirac equation with  external
electromagnetic potentials,  meanwhile,
the $\eps$-scalability of TSFP method
is $h = O(1)$ and $\tau= O(\eps^2)$.
Thus, the EWI-FP and TSFP offer compelling advantages
over commonly used FDTD methods in temporal and spatial resolution when
$0<\eps\ll 1$.

The rest of this paper is organized as follows. In Section \ref{sec2}, several second-order FDTD methods
are reviewed and their stabilities and convergence are analyzed in the nonrelativistic limit regime.
In Section \ref{sec3}, an exponential wave integrator Fourier pseudospectral
method is proposed and analyzed rigorously. In Section
\ref{sec4}, a time-splitting Fourier pseudospectral method
is reviewed and analyzed rigorously. In Section
\ref{sec5}, numerical comparison results are reported. Finally,
some concluding remarks are drawn in Section \ref{sec6}. The mathematical proofs of the error
estimates are given in the appendices, where extensions of EWI-FP and TSFP to higher dimensions are also presented.
Throughout the paper, we adopt the standard notations of Sobolev spaces, use
the notation $p \lesssim q$ to represent that there exists a generic
constant $C$ which is independent of $h$, $\tau$ and $\eps$ such that
$|p|\le C\,q$.


\section{FDTD methods and their analysis}\label{sec2}
\setcounter{equation}{0}
\setcounter{table}{0}
\setcounter{figure}{0}

In this section, we apply the commonly used  FDTD methods to the Dirac equation (\ref{SDEd}) (or
(\ref{SDEdd}))
and analyze their stabilities and convergence in the nonrelativistic limit regime.
For simplicity of notations, we shall only present the numerical methods and their analysis
for (\ref{SDEdd}) in 1D. Generalization to (\ref{SDEd}) and/or higher dimensions is
straightforward and results remain valid without modifications.
Similarly to most works in the literatures for the analysis and computation of the
Dirac equation (cf. \cite{BL, BHM, Hamm0, Hamm, HJMSZ, NSG,  WHJY, WT} and references therein),
in practical computation, we truncate the whole space
problem onto an interval $\Omega = (a,b)$ with
periodic boundary conditions, which is large enough such that
the truncation error is negligible.
In 1D, the Dirac equation (\ref{SDEdd}) with periodic boundary conditions collapses to
\begin{align}
\label{dir}
i\partial_t\Phi(t,x)=&\Bigl[-\frac{i}{\varepsilon}\sigma_1
\partial_x+\frac{1}{\varepsilon^2}\sigma_3\Bigr]\Phi(t,x)+\Bigl[V(t,x)I_2-A_1(t,x)\sigma_1\Bigr]\Phi(t,x),
 \quad x\in\Omega, \quad t>0,\\
\Phi(t,a)=&\Phi(t,b),\quad\partial_x\Phi(t,a) =\partial_x \Phi(t,b),\quad t\geq 0, \qquad
\Phi(0,x) = \Phi_0(x),\quad x\in\overline{\Omega},
\label{dir2}
\end{align}
where $\Phi_0(a)=\Phi_0(b)$ and $\Phi_0^\prime(a)=\Phi_0^\prime(b)$.

\subsection{FDTD  methods}

Choose mesh size $h:=\Delta x =\frac{b-a}{M}$ with $M$ being an even positive integer,
time step $\tau := \Delta t >0$ and denote the grid points and time steps as:
$$
x_j:=a+jh,\quad j=0, 1, \ldots ,M;\qquad t_n:= n\tau,\quad n=0, 1, 2,\ldots .
$$
Denote $X_M=\{\beU=(U_0,U_1,...,U_M)^T\ |\ U_j\in {\mathbb C}^2, j=0,1,\ldots,M, \ U_0=U_M\}$
and  we always use $U_{-1}=U_{M-1}$ and $U_{M+1}=U_1$ if they are involved. For any $\beU\in X_M$,
we denote its Fourier representation as
\be
U_j=\sum_{l=-M/2}^{M/2-1} \widetilde{U}_l\, e^{i\mu_l(x_j-a)}
=\sum_{l=-M/2}^{M/2-1} \widetilde{U}_l\, e^{2ijl\pi/M}, \qquad j=0,1,\ldots, M,
\end{equation}
where $\mu_l$ and $\widetilde{U}_l\in{\mathbb C}^2$  are defined as
\be\mu_l=\frac{2l\pi}{b-a},\qquad \widetilde{U}_l=\frac{1}{M}\sum_{j=0}^{M-1}U_j\, e^{-2ijl\pi/M}, \qquad
 l=-\frac{M}{2}, \ldots,\frac{M}{2}-1.\label{dftco}
\ee
The
 standard $l^2$-norm in $X_M$ is given as
\begin{eqnarray}
\|\beU\|^2_{l^2}=h\sum^{M-1}_{j=0}|U_j|^2, \qquad  \beU\in X_M.
\end{eqnarray}
Let $\Phi_{j}^{n}$  be the numerical approximation of $\Phi(t_n,x_j)$ and
 $V_j^n=V(t_n,x_j)$, $V_j^{n+1/2}=V(t_n+\tau/2,x_j)$, $A_{1,j}^n=A_1(t_n,x_j)$
and $A_{1,j}^{n+1/2}=A_1(t_n+\tau/2,x_j)$  for $0\le j\le M$ and $n\ge0$.
Denote $\Phi^n=\left(\Phi_0^n, \Phi_1^n, \ldots, \Phi_M^n\right)^T\in X_M$ as
the solution vector at $t=t_n$.
Introduce the finite difference discretization operators for $j=0,1,\ldots, M$ and $n\ge0$ as:
\begin{equation*}
\delta_t^+\Phi_{j}^{n}=\frac{\Phi_{j}^{n+1}-\Phi_{j}^{n}}{\tau},\qquad
\delta_t\Phi_{j}^{n}=\frac{\Phi_{j}^{n+1}-\Phi_{j}^{n-1}}{2\tau}, \qquad\delta_x\Phi_{j}^{n}=\frac{\Phi_{j+1}^{n}-\Phi_{j-1}^{n}}{2h},\qquad
\Phi_{j}^{n+\frac{1}{2}}=\frac{\Phi_{j}^{n+1}+\Phi_{j}^{n}}{2}.
\end{equation*}

Here we consider several frequently used FDTD methods to discretize the Dirac equation
(\ref{dir}) for $j=0,1,\ldots, M-1$.

\medskip

{\sl I. Leap-frog finite difference (LFFD) method} 
\begin{eqnarray}
\label{efd1}
i\delta_t \Phi_{j}^{n}=\Bigl[-\frac{i}{\varepsilon}\sigma_1 \delta_x +\frac{1}{\varepsilon^2}\sigma_3\Bigr]\Phi_{j}^{n}+\Bigl[V_j^nI_2-A_{1,j}^n
\sigma_1\Bigr]\Phi_{j}^{n},\quad n\ge1.
\end{eqnarray}

\medskip

{\sl II. Semi-implicit finite difference (SIFD1) method}
\begin{eqnarray}
\label{sifd1}
i\delta_t \Phi_{j}^{n}=-\frac{i}{\varepsilon}\sigma_1\delta_x \Phi_{j}^{n}+\frac{1}{\varepsilon^2}\sigma_3\frac{\Phi_{j}^{n+1}+\Phi_{j}^{n-1}}{2}
+\Bigl[V_j^{n}I_2-A_{1,j}^{n}\sigma_1\Bigr]\frac{\Phi_{j}^{n+1}+\Phi_{j}^{n-1}}{2},\quad n\ge1.
\end{eqnarray}

\medskip

{\sl III. Another semi-implicit finite difference (SIFD2) method}
\begin{eqnarray}
\label{sifd2}
i\delta_t \Phi_{j}^{n}=\Bigl[-\frac{i}{\varepsilon}\sigma_1\delta_x +\frac{1}{\varepsilon^2}\sigma_3\Bigr]\frac{\Phi_{j}^{n+1}+\Phi_{j}^{n-1}}{2}
+\Bigl[V_j^nI_2- A_{1,j}^n\sigma_1\Bigr]\Phi_{j}^{n},\quad n\ge1.
\end{eqnarray}

\medskip

{\sl IV. Crank-Nicolson finite difference (CNFD) method}
\begin{eqnarray}
\label{cnfd1}
i\delta_t^+\Phi_{j}^{n}=\Bigl[-\frac{i}{\varepsilon}\sigma_1\delta_x +\frac{1}{\varepsilon^2}\sigma_3\Bigr]\Phi_{j}^{n+1/2}
+\Bigl[V_j^{n+1/2}I_2-A_{1,j}^{n+1/2}\sigma_1\Bigr]\Phi_{j}^{n+1/2},\quad n\ge0.
\end{eqnarray}

\bigskip

\noindent The initial and boundary conditions in (\ref{dir2}) are discretized as:
\begin{eqnarray}
\label{bdc}
\Phi_{M}^{n+1}=\Phi_{0}^{n+1},\quad \Phi_{-1}^{n+1}=\Phi_{M-1}^{n+1},\quad n\geq 0,
\qquad \Phi_{j}^{0}=\Phi_0(x_j),\quad j=0,1,...,M.
\end{eqnarray}
Using Taylor expansion and noticing (\ref{dir}), the first step
for the LFFD (\ref{efd1}), SIFD1  (\ref{sifd1}) and SIFD2 (\ref{sifd2})
can be computed as
\be\label{bdct1}
\Phi_j^1=\Phi_j^0 +\tau\left[-\frac{1}{\tau}\sin\left(\frac{\tau}{\eps}\right)\sigma_1\Phi_0^\prime(x_j)-i
\left(\frac{1}{\tau}\sin\left(\frac{\tau}{\eps^2}\right)\sigma_3
+V_j^0I_2-A_{1,j}^0\sigma_1\right)\Phi_j^0\right], \quad j=0,1,\ldots, M.
\ee
In the above, we adapt $\frac{1}{\tau}\sin\left(\frac{\tau}{\eps}\right)$ and
$\frac{1}{\tau}\sin\left(\frac{\tau}{\eps^2}\right)$
instead of $\frac{1}{\eps}$ and $\frac{1}{\eps^2}$
such that (\ref{bdct1}) is second order in term of $\tau$ for any fixed $0<\eps\le 1$ and
$\|\Phi^1\|_{\infty}:=\max\limits_{0\le j\le M}|\Phi_j^1|\lesssim 1$ for $0<\eps\le 1$. We remark here
that they can be simply replaced by $1$ when  $\eps=1$.

The above four methods are all time symmetric, i.e. they are unchanged under $\tau\leftrightarrow -\tau$
and $n+1\leftrightarrow n-1$ in the LFFD, SIFD1 and SIFD2 methods or $n+1\leftrightarrow n$
in the CNFD method, and the memory cost is the same at $O(M)$. The LFFD method (\ref{efd1}) is explicit and its computational cost per step
is $O(M)$. In fact, it might be the
simplest and most efficient discretization for the Dirac equation when $\eps=1$ and thus
it has been widely used 
when $\eps=1$. 
The SIFD1 method (\ref{sifd1}) is implicit, however at each time step for $n\ge1$,
the corresponding linear system
is decoupled and  can be solved explicitly for $j=0,1,\ldots,M-1$ as
\be
\Phi_j^{n+1}=\left[(i-\tau V_j^n)I_2-\frac{\tau}{\eps^2}\sigma_3+\tau A_{1,j}^n\sigma_1\right]^{-1}
\left[\left((i+\tau V_j^n)I_2+\frac{\tau}{\eps^2}\sigma_3-\tau A_{1,j}^n\sigma_1\right)\Phi_j^{n-1}
-\frac{2i\tau}{\eps}\sigma_1\delta_x\Phi_j^n\right],
\ee
and thus its computational cost per step is $O(M)$.

The SIFD2 method (\ref{sifd2}) is implicit, however at each time step for $n\ge1$,
the corresponding linear system
is decoupled in phase (Fourier) space and  can be solved explicitly in phase space
for $l=-M/2,\ldots,M/2-1$ as
\be
\widetilde{(\Phi^{n+1})}_l=\left(iI_2-\frac{\tau \sin(\mu_l h)}{\eps h}\sigma_1-\frac{\tau}{\eps^2}\sigma_3\right)^{-1}
\left[\left(iI_2+\frac{\tau \sin(\mu_l h)}{\eps h}\sigma_1+\frac{\tau}{\eps^2}\sigma_3\right)\widetilde{(\Phi^{n-1})}_l+2\tau\widetilde{(G^n\Phi^n)}_l\right],
\ee
where $G^n=(G_0^n,G_1^n,\ldots,G_M^n)^T\in X_M$ with
$G_j^n=- A_{1,j}^n\sigma_1+V_j^nI_2$ for $j=0,1,\ldots,M$,
and thus its computational cost per step is $O(M\ln M)$. The CNFD method (\ref{cnfd1})
is implicit and at each time step for $n\ge0$, the corresponding linear system is coupled and
needs to be solved via either a direct solver or an iterative
solver, and thus its computational cost per step  depends on the linear system solver,
which is usually much larger than $O(M)$, especially in 2D and 3D.
Based on the computational cost per time step, the LFFD method is the most
efficient one and the CNFD method is the most expensive one.

\subsection{Linear stability analysis}

In order to carry out the linear stability analysis for the FDTD methods
via the von Neumann method \cite{D}, we assume that $A_1(t,x)\equiv A_1^0$ and
$V(t,x)\equiv V^0$ with $A_1^0$ and $V^0$ being two real constants
in the Dirac equation (\ref{dir}). Then we have the
following results for the FDTD methods:

\begin{lemma}
\label{THM_STA_FDTD}

(i) The LFFD method (\ref{efd1}) is stable under the stability condition
\begin{equation} \label{sbclffd}
0<\tau\leq \frac{\eps^2h}{|V^0|\eps^2h+\sqrt{h^2+\eps^2(1+\eps h|A_1^0|)^2}}, \qquad
h>0, \quad 0<\eps\le 1.
\end{equation}

(ii) The  SIFD1 method (\ref{sifd1}) is stable  under the stability condition
\begin{equation} \label{scsifd1}
0<\tau\le\varepsilon h, \qquad h>0, \qquad 0<\eps\le 1.
\end{equation}

(iii) The  SIFD2 method (\ref{sifd2}) is stable  under the stability condition
\begin{equation} \label{scsifd2}
0<\tau\le \frac{1}{|V^0|+|A_1^0|}, \qquad h>0, \quad 0<\eps\le 1.
\end{equation}

(iv) The  CNFD method (\ref{cnfd1}) is unconditionally stable,
i.e. it is stable for any $\tau,h>0$ and $0<\varepsilon \le 1$.
\end{lemma}

{\sl Proof}: (i) Plugging
\begin{equation}\label{Phistc}
\Phi^n_{j}=\sum_{l=-M/2}^{M/2-1} \xi_l^n\, \widetilde{(\Phi^0)}_l\,e^{i\mu_l(x_j-a)}
=\sum_{l=-M/2}^{M/2-1} \xi_l^n \, \widetilde{(\Phi^0)}_l\, e^{2ijl\pi/M}, \qquad j=0,1,\ldots, M, \quad n\ge0,
\end{equation}
with $\xi_l\in{\mathbb C}$ and $\widetilde{(\Phi^0)}_l$ being
the amplification factor and the Fourier coefficient at $n=0$, respectively, of the $l$-th mode in the phase space
 into (\ref{efd1}),
using the orthogonality of the Fourier series, we obtain
\begin{equation}\label{xilffd}
\left|(\xi_l^2-1)I_2-2i\tau
\xi_l\left(A_1^0 \sigma_1-V^0I_2-\frac{1}{\eps^2}\sigma_3-\frac{\sin(\mu_l h)}{\eps h}\sigma_1\right)
\right|=0,\qquad l = -\frac{M}{2},..., \frac{M}{2}-1.
\end{equation}
Substituting (\ref{Paulim}) into (\ref{xilffd}), we get that the amplification factor $\xi_l$ satisfies
\begin{equation}
\label{cha}
\xi_l^2-2i\tau\theta_l\xi_l-1=0, \qquad l = -\frac{M}{2},..., \frac{M}{2}-1,
\end{equation}
where
\[\theta_l=-V^0\pm\frac{1}{\varepsilon^2h}\sqrt{h^2+\eps^2
\left(A_1^0\eps h-\sin(\mu_l h)\right)^2},
\qquad l = -\frac{M}{2},..., \frac{M}{2}-1.
\]
Then the  stability condition for the LFFD method (\ref{efd1}) becomes
\begin{equation}
|\xi_l|\leq1 \iff |\tau\theta_l|\leq1,\qquad l = -\frac{M}{2},..., \frac{M}{2}-1,
\end{equation}
which immediately implies the condition (\ref{sbclffd}).

(ii) Similarly to (i), plugging (\ref{Phistc}) into the SIFD1 method (\ref{sifd1}), we have
\begin{equation}
\left|(\xi_l^2-1)I_2-i\tau (\xi_l^2+1)
\left(A_1^0 \sigma_1-V^0I_2-\frac{1}{\eps^2}\sigma_3\right)
+\frac{2i\tau\xi_l \sin(\mu_l h)}{\eps h}\sigma_1\right|=0,\qquad l = -\frac{M}{2},..., \frac{M}{2}-1.
\end{equation}
Noticing (\ref{Paulim}), under the  condition (\ref{scsifd1}), we can get
$|\xi_l|\le 1$ for $l = -\frac{M}{2},..., \frac{M}{2}-1$,
and thus it is stable.

(iii) Similarly to (i), plugging (\ref{Phistc}) into the SIFD2 method (\ref{sifd2}), we have
\begin{equation}
\left|(\xi_l^2-1)I_2+i\tau (\xi_l^2+1)
\left(\frac{1}{\eps^2}\sigma_3+\frac{\sin(\mu_l h)}{\eps h}\sigma_1\right)
-2i\tau \xi_l(A_1^0 \sigma_1-V^0I_2)\right|=0,\quad l = -\frac{M}{2},..., \frac{M}{2}-1.
\end{equation}
Noticing (\ref{Paulim}), under the condition (\ref{scsifd2}), we obtain
\[|\xi_l|\le 1, \qquad l = -\frac{M}{2},..., \frac{M}{2}-1,\]
and thus it is stable.

(iv) Similarly to (i), plugging (\ref{Phistc}) into the CNFD method (\ref{cnfd1}), we obtain
\begin{equation}
\left|(\xi_l-1)I_2+\frac{i\tau}{2} (\xi_l+1)
\left(\frac{1}{\eps^2}\sigma_3-A_1^0 \sigma_1-V^0I_2+\frac{\sin(\mu_l h)}{\eps h}\sigma_1\right)
\right|=0,\quad l = -\frac{M}{2},..., \frac{M}{2}-1.
\end{equation}
Noticing (\ref{Paulim}), we have for $l = -\frac{M}{2},..., \frac{M}{2}-1$,
\begin{equation}
|\xi_l|=\left|\frac{2+i\tau \theta_l}{2-i\tau\theta_l}\right|=1, \quad\theta_l=V^0\pm\frac{1}{\varepsilon^2h}\sqrt{h^2+\eps^2\left(A_1^0\eps h-
\sin(\mu_l h)\right)^2}.
\end{equation}
Thus it is unconditionally stable. \hfill $\Box$


\subsection{Mass and energy conservation}
For the CNFD method (\ref{cnfd1}),  we have the following conservative properties.

\begin{lemma}
\label{Mass_STA_FDTD}
The  CNFD (\ref{cnfd1}) conserves
the mass in the discretized level, i.e.
\be\label{cons976}
\|\Phi^n\|_{l^2}^2:=h\sum_{j=0}^{M-1}|\Phi_j^{n}|^2\equiv  h\sum_{j=0}^{M-1}|\Phi_j^{0}|^2
=\|\Phi^0\|_{l^2}^2=h\sum_{j=0}^{M-1}|\Phi_0(x_j)|^2,\qquad n\ge0.
\ee
Furthermore, if $V(t,x)=V(x)$ and $A_1(t,x)=A_1(x)$ are time independent, the
CNFD (\ref{cnfd1}) conserves the energy as well,
\be\label{cons:energy}
\begin{split}
E_h^n=&-\frac{ih}{\eps}\sum\limits_{j=1}^{M-1}(\Phi_j^{n})^*\sigma_1\delta_x\Phi_j^{n}
+\frac{h}{\eps^2}\sum\limits_{j=0}^{M-1}(\Phi_j^{n})^*\sigma_3\Phi_j^{n}
+h\sum\limits_{j=0}^{M-1}V_j(\Phi_j^{n})^*\sigma_3\Phi_j^{n}
-h\sum\limits_{j=0}^{M-1}A_{1,j}(\Phi_j^{n})^*\sigma_1\Phi_j^{n}\\
\equiv& E_h^0,\qquad n\ge0,
\end{split}
\ee
where $V_j=V(x_j)$ and $A_{1,j}=A_1(x_j)$ for $j=0,1,\ldots,M$.
\end{lemma}

{\sl Proof}: (i) Firstly, we prove the mass conservation (\ref{cons976}).
Multiplying both sides of (\ref{cnfd1}) from left by $h\tau\, (\Phi_j^{n+1/2})^*$
and taking the imaginary part, we have
\be \label{Phic78}
h|\Phi_j^{n+1}|^2=h|\Phi_j^{n}|^2-\frac{\tau h}{2\eps}\left[
(\Phi_j^{n+1/2})^*\sigma_1\delta_x\Phi_j^{n+1/2}+(\Phi_j^{n+1/2})^T\sigma_1\delta_x
\overline{\Phi}_j^{n+1/2}\right],
\qquad j=0,1,\ldots,M-1.
\ee
Summing (\ref{Phic78}) for $j=0,1,\ldots,M-1$ and noticing (\ref{Paulim}), we get
\bea
\|\Phi^{n+1}\|_{l^2}^2&=&\|\Phi^{n}\|_{l^2}^2-\frac{\tau h}{2\eps}\sum_{j=0}^{M-1}\left[ (\Phi_j^{n+1/2})^*\,
\sigma_1 \delta_x \Phi_j^{n+1/2}+(\Phi_j^{n+1/2})^T\,
\sigma_1 \delta_x \overline{\Phi}_j^{n+1/2}\right]\nonumber\\
&=&\|\Phi^{n}\|_{l^2}^2-\frac{\tau }{2\eps}\sum_{j=0}^{M-1}\Bigl[ (\Phi_j^{n+1/2})^*\,
\sigma_1 \Phi_{j+1}^{n+1/2}+(\Phi_j^{n+1/2})^T\,
\sigma_1 \overline{\Phi}_{j+1}^{n+1/2}\nonumber\\
&&\qquad \qquad \qquad \qquad -(\Phi_{j+1}^{n+1/2})^*\,
\sigma_1 \Phi_{j}^{n+1/2}-(\Phi_{j+1}^{n+1/2})^T\,
\sigma_1 \overline{\Phi}_{j}^{n+1/2}\Bigr]\nonumber\\
&=&\|\Phi^{n}\|_{l^2}^2, \qquad n\ge0,
\eea
which immediately implies (\ref{cons976}) by induction.

(ii) Secondly, we prove the energy conservation (\ref{cons:energy}). Multiplying both sides of (\ref{cnfd1}) from left by $2h\, (\Phi_j^{n+1}-\Phi_j^n)^*$
and taking the real part, we have
\begin{align}
&-h\,\text{Re}\left[\frac{i}{\eps}(\Phi_j^{n+1}-\Phi_j^n)^*\sigma_1\delta_x(\Phi_j^{n+1}+\Phi_j^n)\right]
+\frac{h}{\eps^2}\left[(\Phi_j^{n+1})^*\sigma_3\Phi_j^{n+1}-(\Phi_j^{n})^*\sigma_3\Phi_j^{n}\right]\nonumber\\
&\qquad+hV_j(|\Phi_j^{n+1}|^2-|\Phi_j^n|^2)-
hA_{1,j}\left[(\Phi_j^{n+1})^*\sigma_1\Phi_j^{n+1}-(\Phi_j^{n})^*\sigma_1\Phi_j^{n}\right]=0.\label{eq:enercons1}
\end{align}
Summing (\ref{eq:enercons1}) for $j=0,1,\ldots,M-1$ and noticing the summation by parts formula, we have
\begin{align*}
h\sum\limits_{j=0}^{M-1}\text{Re}\left(\frac{i}{\eps}(\Phi_j^{n+1}-\Phi_j^n)^*\sigma_1\delta_x(\Phi_j^{n+1}+\Phi_j^n)\right)
=\frac{i h}{\eps}\sum\limits_{j=0}^{M-1}(\Phi_j^{n+1})^*\sigma_1\delta_x\Phi_j^{n+1}-
\frac{i h}{\eps}\sum\limits_{j=0}^{M-1}(\Phi_j^n)^*\sigma_1\delta_x\Phi_j^n,
\end{align*}
and
\begin{align}
&-\frac{ih}{\eps}\sum\limits_{j=0}^{M-1}(\Phi_j^{n+1})^*\sigma_1\delta_x\Phi_j^{n+1}+
\frac{ih}{\eps}\sum\limits_{j=0}^{M-1}(\Phi_j^n)^*\sigma_1\delta_x\Phi_j^n
+\frac{h}{\eps^2}\sum\limits_{j=0}^{M-1}\left((\Phi_j^{n+1})^*\sigma_3\Phi_j^{n+1}-(\Phi_j^{n})^*\sigma_3\Phi_j^{n}\right)\nonumber\\
&\qquad+h\sum\limits_{j=0}^{M-1}V_j(|\Phi_j^{n+1}|^2-|\Phi_j^n|^2)-
h\sum\limits_{j=0}^{M-1}A_{1,j}\left((\Phi_j^{n+1})^*\sigma_1\Phi_j^{n+1}-(\Phi_j^{n})^*\sigma_1\Phi_j^{n}\right)=0,\label{eq:enercons2}
\end{align}
which immediately implies (\ref{cons:energy}).

\hfill $\Box$

\subsection{Error estimates}

Let  $0<T<T^*$ with $T^*$ being the maximal existence time of the solution, and denote $\Omega_{T}=[0,T]\times \Omega$. Motivated by the nonrelativistic limit
of the Dirac equation \cite{BMP} and the dispersion relation (\ref{dispp}),
we assume that the exact solution of (\ref{dir}) satisfies $\Phi\in C^{3}([0,T]; (L^{\infty}(\Omega))^2)\cap C^{2}([0,T]; (W_p^{1,\infty}(\Omega))^2)
\cap C^{1}([0,T]; (W_p^{2,\infty}(\Omega))^2)\cap C([0,T]; (W_p^{3,\infty}(\Omega))^2)$ and
\be
(A) \qquad \qquad \left\|\frac{\partial^{r+s}}{\partial t^{r}
\partial x^{s}}\Phi\right\|_{L^\infty([0,T];(L^{\infty}(\Omega))^2)}
\lesssim \frac{1}{\varepsilon^{2r}},
\quad0\leq r\leq 3,\ 0\leq r+s\leq 3, \qquad 0<\eps\le 1,\hskip1cm
\ee
where $W_p^{m,\infty}(\Omega)=\{u\ |\ u\in W^{m,\infty}(\Omega),\  \partial_x^l u(a)=\partial_x^l u(b),\
l=0,\ldots,m-1\}$ for $m\ge1$ and here the boundary values are understood in the trace sense. In the subsequent discussion,
we will omit $\Omega$ when referring to the space norm taken on $\Omega$.
In addition, we assume the electromagnetic potentials
$V\in C(\overline{\Omega}_T)$ and $A_1\in C(\overline{\Omega}_T)$ and denote
\begin{eqnarray}
(B) \qquad \qquad \qquad \qquad V_{\rm max}:=\max_{(t,x)\in\overline{\Omega}_T}|V(t,x)|,
\qquad A_{1,\rm max}:=\max_{(t,x)\in\overline{\Omega}_T}|A_1(t,x)|.\qquad \qquad \qquad \qquad
\end{eqnarray}
Define the grid error function $\bee^n=(\bee_0^n,\bee_1^n,\ldots,\bee_M^n)^T\in X_M$ as:
\begin{equation}
\label{err}
\bee_{j}^n = \Phi(t_n,x_j) - \Phi_{j}^n, \qquad j=0,1,\ldots,M, \quad n\ge0,
\end{equation}
with $\Phi_j^n$ being the approximations obtained from the FDTD methods.

For the CNFD (\ref{cnfd1}), we can establish the following error bound (see its proof in Appendix A).

\begin{thm}
\label{thm_cnfd}
Under the assumptions (A) and (B), there exist constants $h_{0}>0$ and $\tau_0>0$
sufficiently small and independent of $\varepsilon$, such that for any
$0<\varepsilon\leq1$, $0<h\leq h_{0}$ and $0<\tau\leq\tau_0$,
we have the following error estimate for the CNFD  (\ref{cnfd1}) with (\ref{bdc})
\begin{equation} \label{ebcnfd97}
\|\bee^{n}\|_{l^{2}}\lesssim \frac{h^{2}}{\varepsilon}+\frac{\tau^{2}}{\varepsilon^{6}},
\qquad 0\leq n\leq\frac{T}{\tau}.
\end{equation}
\end{thm}

For the LFFD  (\ref{efd1}), we assume the stability condition
\begin{equation}\label{stclffd8}
0<\tau\leq \frac{\eps^2h}{\eps^2hV_{\rm max}+\sqrt{h^2+\eps^2(1+\eps hA_{1,\rm max})^2}},
\qquad h>0, \quad 0<\eps\le 1,
\end{equation}
and establish the following error estimate (see its proof in Appendix B).

\begin{thm}
\label{thm_efd}
Under the assumptions (A) and (B), there exist constants $h_{0}>0$ and $\tau_0>0$
sufficiently small and independent of $\varepsilon$, such that for any
$0<\varepsilon\leq1$, when $0<h\leq h_{0}$ and $0<\tau\leq\tau_0$ and under
the stability condition (\ref{stclffd8}),
we have the following error estimate for the LFFD  (\ref{efd1}) with (\ref{bdc}) and
(\ref{bdct1})
\begin{equation}
\label{eblffd}
\|\bee^{n}\|_{l^{2}}\lesssim \frac{h^{2}}{\varepsilon}+\frac{\tau^{2}}{\varepsilon^{6}},
\qquad 0\leq n\leq\frac{T}{\tau}.
\end{equation}
\end{thm}

Similarly to the proofs of the LFFD and CNFD methods, error estimates for SIFD1 (\ref{sifd1}) and SIFD2 (\ref{sifd2}) can be derived
and the details are omitted here for brevity.
For the  SIFD2 (\ref{sifd2}), we assume the stability condition
\begin{equation}\label{stclffd9}
0<\tau\leq \frac{1}{V_{\rm max}+A_{1,\rm max}}, \qquad h>0, \quad 0<\eps\le 1,
\end{equation}
and establish the following error estimates.

\begin{thm}
\label{thm_sifd1}
Under the assumptions (A) and (B), there exist constants $h_{0}>0$ and $\tau_0>0$
sufficiently small and independent of $\varepsilon$, such that for any
$0<\varepsilon\leq1$, when $0<h\leq h_{0}$ and $0<\tau\leq\tau_0$ and under
the stability condition (\ref{scsifd1}),
we have the following error estimate for the SIFD1  (\ref{sifd1}) with (\ref{bdc}) and
(\ref{bdct1})
\begin{equation}
\|\bee^{n}\|_{l^{2}}\lesssim \frac{h^{2}}{\varepsilon}+\frac{\tau^{2}}{\varepsilon^{6}},
\qquad 0\leq n\leq\frac{T}{\tau}.
\end{equation}
\end{thm}

\begin{thm}
\label{thm_sifd2}
Under the assumptions (A) and (B), there exist constants $h_{0}>0$ and $\tau_0>0$
sufficiently small and independent of $\varepsilon$, such that for any
$0<\varepsilon\leq1$, when $0<h\leq h_{0}$ and $0<\tau\leq\tau_0$ and under
the stability condition (\ref{stclffd9}),
we have the following error estimate for the SIFD2  (\ref{sifd2}) with (\ref{bdc}) and
(\ref{bdct1})
\begin{equation}
\|\bee^{n}\|_{l^{2}}\lesssim \frac{h^{2}}{\varepsilon}+\frac{\tau^{2}}{\varepsilon^{6}},
\qquad 0\leq n\leq\frac{T}{\tau}.
\end{equation}
\end{thm}

Based on Theorems \ref{thm_cnfd}-\ref{thm_sifd2}, the four FDTD methods studied here share
the same temporal/spatial resolution capacity in the nonrelativistic limit regime. In fact,
given an accuracy bound $\delta>0$, the $\eps$-scalability of the four FDTD methods is:
\be
\tau=O\left(\eps^3 \sqrt{\delta}\right)=O(\eps^3), \qquad  h=O\left(\sqrt{\delta\eps}\right)=O\left(\sqrt{\eps}\right),
\qquad 0<\eps\ll1.
\ee


\section{An EWI-FP method and its analysis} \label{sec3}
\setcounter{equation}{0}
\setcounter{table}{0}
\setcounter{figure}{0}

In this section, we propose an exponential wave integrator Fourier pseudospectral (EWI-FP)
method  to solve the Dirac equation (\ref{SDEd}) (or (\ref{SDEdd}))
and establish its stability and convergence in the nonrelativistic limit regime.
Again, for simplicity of notations, we shall only present the numerical method and its analysis
for (\ref{dir}) in 1D. Generalization to (\ref{SDEd}) and/or higher dimensions is
straightforward and the results remain valid without modifications (see  generalizations in Appendix D).

\subsection{The EWI-FP method}

Denote
\[
Y_{M}=Z_M\times Z_M, \qquad \hbox{with}\quad
Z_M={\rm span}\left\{\phi_l(x)=e^{i\mu_l(x-a)},\ l=-\frac{M}{2}, -\frac{M}{2}+1,\ldots, \frac{M}{2}-1\right\}.
\]
Let $[C_p(\overline{\Omega})]^2$ be the function space consisting of all periodic vector function $U(x):\ \overline{\Omega}=[a,b]\to {\mathbb C}^2$. For any $U(x)\in [C_p(\overline{\Omega})]^2$ and $U\in X_M$,
define $P_M:\ [L^2(\Omega)]^2\rightarrow Y_M$ as the standard projection operator \cite{ST}, $I_M:\ [C_p(\overline{\Omega})]^2\rightarrow Y_M$ and $I_M:X_M\rightarrow Y_M$ as the standard interpolation operator \cite{ST}, i.e.
\be
(P_MU)(x)=\sum_{l=-M/2}^{M/2-1}\widehat{U}_l\, e^{i\mu_l(x-a)},\quad (I_MU)(x)=\sum_{l=-M/2}^{M/2-1}\widetilde{U}_l\, e^{i\mu_l(x-a)},\qquad
a\leq x\leq b,
\ee
with
\be\label{fouriercoef}
\widehat{U}_l=\frac{1}{b-a}\int_a^bU(x)\,e^{-i\mu_l(x-a)}\,dx,\quad \widetilde{U}_l=\frac{1}{M}\sum_{j=0}^{M-1}U_j\, e^{-2ijl\pi/M},\qquad
l=-\frac{M}{2}, -\frac{M}{2}+1,\ldots, \frac{M}{2}-1,
\ee
where $U_j=U(x_j)$ when $U$ is a function.

The Fourier spectral discretization for the Dirac equation (\ref{dir}) is as follows:

\noindent Find $\Phi_{M}(t,x)\in Y_M$, i.e.
\begin{eqnarray}
\label{FPR}
\Phi_{M}(t,x)=\sum^{M/2-1}_{l=-M/2}\widehat{(\Phi_M)}_l(t)\,
e^{i\mu_l(x-a)},\qquad a\le x\le b, \qquad t\ge0,
\end{eqnarray}
such that for $a<x<b$ and $t>0$,
\be
\label{FTD}
i\partial_t\Phi_{M}(t,x)=\left[-\frac{i}{\varepsilon}\sigma_1\partial_x
+\frac{1}{\varepsilon^2}\sigma_3\right]\Phi_{M}(t,x)+P_M(V\Phi_M)(t,x)-\sigma_1 P_M(A_1\Phi_M)(t,x).
\ee
Substituting (\ref{FPR}) into (\ref{FTD}), noticing the orthogonality of $\phi_l(x)$,
we get for $l=-\frac{M}{2},\ldots,\frac{M}{2}-1$,
\be
i\,\frac{d}{dt}\widehat{(\Phi_M)}_{l}(t)=\left[\frac{\mu_l}{\varepsilon}\sigma_1+\frac{1}{\varepsilon^2}\sigma_3\right]
\widehat{(\Phi_M)}_{l}(t)+\widehat{(V\Phi_M)}_l(t)-\sigma_1\widehat{(A_1\Phi_{M})}_l(t)=0, \qquad t\ge0.
\ee
For each  $l$ ($l=-\frac{M}{2}, -\frac{M}{2}+1,\ldots, \frac{M}{2}-1$),
when $t$ is near $t=t_n$ $(n\geq 0)$, we rewrite the above ODEs  as
\be\label{ODE765}
i\,\frac{d}{ds}\widehat{(\Phi_M)}_{l}(t_n+s)=\frac{1}{\eps^2}\Gamma_l\,
\widehat{(\Phi_M)}_{l}(t_n+s)+\widehat{F}_l^n(s),\qquad s\in{\mathbb R},
\ee
where
\begin{equation}\label{eq:fdef}
\widehat{F}_l^n(s)=\widehat{(G\Phi_M)}_l(t_n+s),\qquad G(t,x)=V(t,x)I_2-\sigma_1A_1(t,x), \qquad
s,t\in{\mathbb R},
\end{equation}
 and $\Gamma_l=\mu_l\eps\sigma_1+\sigma_3=Q_l\, D_l\, (Q_l)^*$ with
\begin{equation}\label{eq:Gamma}
\Gamma_l=\begin{pmatrix}
1 &\mu_l\eps\\
\mu_l\eps &-1
\end{pmatrix}, \quad Q_l=\begin{pmatrix}
\frac{1+\delta_l}{\sqrt{2\delta_l(1+\delta_l)}} &-\frac{\eps\mu_l}{\sqrt{2\delta_l(1+\delta_l)}}\\
\frac{\eps\mu_l}{\sqrt{2\delta_l(1+\delta_l)}} &\frac{1+\delta_l}{\sqrt{2\delta_l(1+\delta_l)}}
\end{pmatrix}, \quad D_l=\begin{pmatrix}
\delta_l &0\\
0 &-\delta_l\\
\end{pmatrix},\quad \delta_l=\sqrt{1+\eps^2\mu_l^2}.
\end{equation}
Solving the above ODE (\ref{ODE765}) via the integrating factor method, we obtain
\begin{eqnarray}
\label{EWIS}
\widehat{(\Phi_M)}_l(t_n+s)=e^{-is\Gamma_l/\eps^2}\widehat{(\Phi_M)}_l(t_n)
-i\int_0^{s}e^{i(w-s)\Gamma_l/\eps^2}\widehat{F}_l^n(w)\,dw, \qquad s\in{\mathbb R}.
\end{eqnarray}
Taking $s=\tau$ in (\ref{EWIS}) we have
\begin{equation}\label{EWIn}
\widehat{(\Phi_M)}_l(t_{n+1})=e^{-i\tau\Gamma_l/\eps^2}\widehat{(\Phi_M)}_l(t_{n})
-i\int_0^{\tau}e^{\frac{i(w-\tau)}{\eps^2}\Gamma_l}
\widehat{F}_l^n(w)dw.
\end{equation}
To obtain an explicit  numerical method with second  order accuracy in time, we approximate the integrals in (\ref{EWIn})
via the Gautschi-type rules \cite{Gautschi0, Gautschi-type-3, Gautschi-type-2}, which have been widely
used  for integrating highly oscillatory ODEs \cite{ASS,BD, Gautschi0, Gautschi-type-3, Gautschi-type-2,Iser,Iser2}, as
\be\label{eq:inteapp0}
\int_0^{\tau}e^{\frac{i(w-\tau)}{\eps^2}\Gamma_l}\widehat{F}_l^0(w)\,dw
\approx\int_0^{\tau}e^{\frac{i(w-\tau)}{\eps^2}\Gamma_l}\,dw\,\widehat{F}_l^0(0)
=-i\eps^2\Gamma_l^{-1}\left[I_2-e^{-\frac{i\tau}{\eps^2}\Gamma_l}\right]
\widehat{F}_l^0(0),
\ee
and
\begin{align}\label{eq:inteappn}
&\int_0^{\tau}e^{\frac{i(w-\tau)}{\eps^2}\Gamma_l}\widehat{F}_l^n(w)dw
\approx\int_0^{\tau}e^{\frac{i(w-\tau)}{\eps^2}\Gamma_l}\left(\widehat{F}_l^n(0)+
w\delta_t^-\widehat{F}_l^n(0)\right)dw\nonumber\\
&=-i\eps^2\Gamma_l^{-1}\left[I_2-e^{-\frac{i\tau}{\eps^2}\Gamma_l}\right]\widehat{F}_l^n(0)+
\left[-i\eps^2\tau\Gamma_l^{-1}+\eps^4\Gamma_l^{-2}\left( I_2-e^{-\frac{i\tau}{\eps^2}\Gamma_l}\right)\right]\delta_t^-\widehat{F}_l^n(0), \qquad n\ge1,
\end{align}
where we have approximated the time derivative $\partial_t\widehat{F}_l^n(s)$ at $s=0$ for
$n\ge1$ by the finite difference as
\be
\partial_t\widehat{F}_l^n(0)\approx\delta_t^-\widehat{F}_l^n(0)=\frac{\widehat{F}_l^{n}(0)-\widehat{F}_l^{n-1}(0)}{\tau}.
\ee
Now, we are ready to describe our scheme. Let $\Phi_M^n(x)$ be the approximation of  $\Phi_{M}(t_n,x)$ ($n\ge0$).
Choosing $\Phi_M^0(x)=(P_M\Phi_0)(x)$,  an {\sl exponential wave integrator  Fourier spectral} (EWI-FS)
discretization for the Dirac equation (\ref{dir}) is to update the numerical approximation $\Phi^{n+1}_M(x)\in Y_M$ ($n=0,1,\ldots$) as
\begin{equation}\label{AL1}
\Phi_M^{n+1}(x)=
\sum_{l=-M/2}^{M/2-1}\widehat{(\Phi_M^{n+1})}_l\, e^{i\mu_l(x-a)},\qquad a\leq x\leq b,\qquad n \ge0,
\end{equation}
where for $l = -\frac{M}{2}, .., \frac{M}{2}-1$,
\be\label{Coe2}
\widehat{(\Phi_M^{n+1})}_l=\left\{\ba{ll}
e^{-i\tau\Gamma_l/\eps^2}\widehat{(\Phi_M^0)}_l
-i\eps^2\Gamma_l^{-1}\left[I_2-e^{-\frac{i\tau}{\eps^2}\Gamma_l}\right]
\widehat{(G(t_0)\Phi_M^0)}_l, &n=0,\\
e^{-i\tau\Gamma_l/\eps^2}\widehat{(\Phi_M^n)}_l
-iQ_l^{(1)}(\tau)
\, \widehat{(G(t_n)\Phi_M^n)}_l-iQ_l^{(2)}(\tau)\delta_t^-\widehat{\left(G(t_n)\Phi_M^n\right)}_l, &n\ge1,\\
\ea\right.
\ee
with $G(t)$ denoting $G(t,x)$ and the matrices $Q_l^{(1)}(\tau)$ and $Q_l^{(2)}(\tau)$ given as
\be
Q_l^{(1)}(\tau)=-i\eps^2\Gamma_l^{-1}\left[I_2-e^{-\frac{i\tau}{\eps^2}\Gamma_l}\right],
\quad
Q_l^{(2)}(\tau)=-i\eps^2\tau\Gamma_l^{-1}+\eps^4\Gamma_l^{-2}\left( I_2-e^{-\frac{i\tau}{\eps^2}\Gamma_l}\right).
\ee

The above procedure is not suitable in practice due to the difficulty in computing the Fourier coefficients through integrals in (\ref{fouriercoef}).
 Here we present an efficient implementation by choosing $\Phi_M^0(x)$ as the interpolant of $\Phi_0(x)$ on the grids $\left\{x_j, j = 0, 1,\ldots, M\right\}$ and approximate the integrals
in (\ref{fouriercoef}) by a quadrature rule.

Let $\Phi_j^n$  be the numerical approximation of $\Phi(t_n,x_j)$ for $j = 0, 1, 2, \ldots, M$ and $n\geq 0$, and denote $\Phi^n\in X_M$ as the vector with components $\Phi_j^n$.
 Choosing $\Phi^0_j=\Phi_0(x_j)$ ($j = 0, 1,\ldots, M$), an {\sl EWI Fourier pseudospectral} (EWI-FP) method for computing $\Phi^{n+1}$ for $n\geq 0$ reads
\begin{equation}\label{AL4}
\Phi_{j}^{n+1} = \sum_{l=-M/2}^{M/2-1}\widetilde{(\Phi^{n+1})}_le^{2ijl\pi/M}, \quad j = 0, 1, ..., M,
\end{equation}
where
\begin{equation}\label{AL5}
\widetilde{(\Phi^{n+1})}_l=\left\{\ba{ll}
e^{-i\tau\Gamma_l/\eps^2}\widetilde{(\Phi_0)}_l-i\eps^2\Gamma_l^{-1}
\left[I_2-e^{-\frac{i\tau}{\eps^2}\Gamma_l}\right]\widetilde{(G(t_0)\Phi_0)}_l, &n=0,\\
e^{-i\tau\Gamma_l/\eps^2}\widetilde{(\Phi^{n})}_l-iQ_l^{(1)}(\tau)\widetilde{(G(t_n)\Phi^{n})}_l
-iQ_l^{(2)}(\tau)\delta_t^-\widetilde{\left(G(t_n)\Phi^n\right)}_l, &n\ge1.\\
\ea\right.
\end{equation}
The EWI-FP (\ref{AL4})-(\ref{AL5}) is explicit, and can be computed efficiently by the
fast Fourier transform (FFT). The memory
cost is $O(M)$ and the computational cost per time step is $O(M \log M)$.

\subsection{Linear stability analysis}

To consider the linear stability, we assume that in the Dirac equation (\ref{dir}),
the external potential fields are constants, i.e. $A_1(t,x)\equiv A_1^0$ and
$V(t,x)\equiv V^0$ with $A_1^0$ and $V^0$ being two real constants. In this case, we adopt the Von Neumann stability requirement that
the errors grow exponentially at most. Then we have

\begin{lemma}
The EWI-FP method (\ref{AL4})-(\ref{AL5}) and  EWI-FS method (\ref{AL1})-(\ref{Coe2}) are  stable under the
stability condition
\begin{equation} \label{sbcexp}
0<\tau\lesssim 1, \quad 0<\eps\le 1.
\end{equation}
\end{lemma}

{\sl Proof}: We shall only prove the EWI-FS case  (\ref{AL1})-(\ref{Coe2}),
as the EWI-FP method case (\ref{AL5}) is quite the same.
Similarly to the proof of Lemma \ref{THM_STA_FDTD}, noticing (\ref{Coe2}),
(\ref{eq:fdef}), (\ref{EWIn}) and
(\ref{eq:inteappn}), we find that
\begin{equation}\label{eq:ewifpstab}
\xi_l^2(\widetilde{\Phi^0})_l= \xi_le^{-i\tau\Gamma_l/\eps^2}(\widetilde{\Phi^0})_l
-i\int_0^{\tau}e^{i(w-\tau)\Gamma_l/\eps^2}(V^0I_2-A_1^0\sigma_1)\left(\xi_l+
\frac{w}{\tau}(\xi_l-1)\right)(\widetilde{\Phi^0})_l\,dw.
\end{equation}
Denoting $C=|V^0|+|A_1^0|$, taking the $l^2$ norms of the vectors on  both sides of (\ref{eq:ewifpstab}) and then dividing both sides by the $l^2$ norm
of $(\widetilde{\Phi^0})_l$,  in view of the properties of $e^{-is\Gamma_l/\eps^2}$, we get
\be
|\xi_l|^2\leq \left(1+C\tau+\frac{C}{2}\tau\right)|\xi_l|+\frac{C}{2}\tau,
\ee
which implies
\be
\left(|\xi_l|-\frac{1+3C\tau/2}{2}\right)^2\leq \frac{1+5C\tau+9C^2\tau^2/4}{4}\leq \frac{(1+5C\tau/2)^2}{4}.
\ee
Thus, we obtain
\be
|\xi_l|\leq 1+4C\tau,\quad l=-\frac{M}{2},\ldots,\frac{M}{2}-1,
\ee
and it follows that the EWI-FS (\ref{AL1})-(\ref{Coe2}) is stable under the stability condition (\ref{sbcexp}).
$\hfill\Box$

\subsection{Error estimates}

In order to obtain an error estimate for the EWI methods (\ref{AL1})-(\ref{Coe2}) and (\ref{AL4})-(\ref{AL5}),
motivated by the results in \cite{BMS2,CC}, we assume that there exists an integer $m_0\geq 2$ such that the exact solution $\Phi(t,x)$ of the Dirac equation (\ref{dir}) satisfies
\begin{align*}
(C)\quad \|\Phi\|_{ L^\infty([0,T]; (H_p^{m_0})^2)}\lesssim 1,\qquad
\|\partial_t\Phi\|_{L^{\infty}([0, T];(L^2)^2)}\lesssim \frac{1}{\varepsilon^2},\quad \|\partial_{tt}\Phi\|_{L^{\infty}([0,T];(L^2)^2)}\lesssim\frac{1}{\varepsilon^4},\hskip4cm
\end{align*}
where $H^k_p(\Omega)=\{u\ |\ u\in H^{k}(\Omega),\  \partial_x^l u(a)=\partial_x^l u(b),\
l=0,\ldots,k-1\}$.
In addition, we assume the electromagnetic potentials satisfy
\begin{equation*}
(D)\hskip5cm \|V\|_{W^{2,\infty}([0,T];L^\infty)}+\|A_1\|_{W^{2,\infty}([0,T];L^\infty)}\lesssim 1.\hskip4cm
\end{equation*}
The following estimate can be established (see its proof in
Appendix C).
\begin{thm}
\label{thm_EWI}
Let $\Phi_M^n(x)$ be the approximation obtained from the EWI-FS (\ref{AL1})-(\ref{Coe2}). Under the
assumptions (C) and (D), there exists $h_0>0$ and $\tau_0>0$ sufficiently small
and independent of $\varepsilon$ such that, for any $0<\varepsilon\leq1$,
when $0<h\leq h_0$ and $0<\tau\leq\tau_0$ satisfying the stability condition
(\ref{stclffd9}), we have the following error estimate
\begin{equation}
\label{thm_eq_EWI}
\|\Phi(t_n,x)-\Phi_M^n(x)\|_{L^2}\lesssim\frac{\tau^2}{\varepsilon^4}+h^{m_0},\quad 0\leq n\leq\frac{T}{\tau}.
\end{equation}
\end{thm}
\begin{rmk}
The same error estimate in Theorem \ref{thm_EWI} holds for the
EWI-FP (\ref{AL4})-(\ref{AL5}) and the proof is quite similar to that of Theorem \ref{thm_EWI}.
\end{rmk}


\section{A TSFP method and its analysis} \label{sec4}
\setcounter{equation}{0}
\setcounter{table}{0}
\setcounter{figure}{0}

In this section, we present a time-splitting Fourier pseudospectral (TSFP)
method  to solve the Dirac equation(\ref{SDEd}) (or (\ref{SDEdd})) which has
been proposed and studied for the Maxwell-Dirac equation \cite{BL,HJMSZ}.
Again, for simplicity of notations, we shall only present the numerical method and its analysis
for (\ref{dir}) in 1D. Generalization to (\ref{SDEd}) and/or higher dimensions is
straightforward and results remain valid without modifications (see generalizations in Appendix D).


From time $t = t_n$ to time $t = t_{n+1}$,
the Dirac equation (\ref{dir}) is split into two steps. One solves first
\be
\label{dir1st1}
i\partial_t\Phi(t,x)=\left[-\frac{i}{\varepsilon}\sigma_1
\partial_x+\frac{1}{\varepsilon^2}\sigma_3\right]\Phi(t,x), \quad x\in\Omega,
\ee
with the periodic boundary condition (\ref{dir2})
for the time step of length $\tau$, followed by solving
\be
\label{dir1st2}
i\partial_t\Phi(t,x)=\left[-A_1(t,x)\sigma_1
+V(t,x)I_2\right]\Phi(t,x), \quad x\in\Omega,
\ee
for the same time step. Eq. (\ref{dir1st1}) will be first discretized
in space by the Fourier spectral method and then
integrated (in phase or Fourier space) in
time {\sl exactly} \cite{BL}. For the ODEs (\ref{dir1st2}), we can integrate
{\sl analytically} in time  as
\be
\Phi(t,x)=e^{-i\int_{t_n}^t \left[ V(s,x)\, I_2-A_1(s,x)\, \sigma_1\right]ds}\,\Phi(t_n,x),
\qquad a\le x\le b, \quad t_n\le t\le t_{n+1}.
\ee
In practical computation, from time $t=t_n$ to $t=t_{n+1}$, one
often combines the splitting steps via the standard Strang splitting \cite{St}
-- which results in a second order time-splitting Fourier pseudospectral (TSFP)
method -- as
\begin{equation}\label{eq:tsfp}
\begin{split}
&\Phi_j^{(1)}=\sum_{l=-M/2}^{M/2-1} e^{-i\tau \Gamma_l/2\eps^2}\,\widetilde{(\Phi^n)}_l\, e^{i\mu_l(x_j-a)}
=\sum_{l=-M/2}^{M/2-1} Q_l\, e^{-i \tau D_l/2\eps^2}\,(Q_l)^*\,
\widetilde{(\Phi^n)}_l\, e^{\frac{2ijl\pi}{M}},\\
&\Phi_j^{(2)}=e^{-i\int_{t_n}^{t_{n+1}} G(t,x_j)\,dt}\,\Phi_j^{(1)} = P_j\, e^{-i \Lambda_j}\, P_j^*\, \Phi_j^{(1)},
\qquad \qquad j=0,1,\ldots,M,\qquad n\ge0,\\
&\Phi_j^{n+1}=\sum_{l=-M/2}^{M/2-1} e^{-i \tau \Gamma_l/2\eps^2}\,\widetilde{(\Phi^{(2)})}_l\, e^{i\mu_l(x_j-a)}
=\sum_{l=-M/2}^{M/2-1} Q_l\, e^{-i \tau D_l/2\eps^2}\,(Q_l)^*\,
\widetilde{(\Phi^{(2)})}_l\, e^{\frac{2ijl\pi}{M}},
\end{split}
\end{equation}
where
$\int_{t_n}^{t_{n+1}}G(t,x_j)dt=V_j^{(1)}\, I_2-
A_{1,j}^{(1)}\,\sigma_1=P_j\, \Lambda_j\, P_j^*$ with $V_j^{(1)}=\int_{t_n}^{t_{n+1}}V(t,x_j)dt$,
$A_{1,j}^{(1)}=\int_{t_n}^{t_{n+1}}A_1(t,x_j)dt$, $\Lambda_j={\rm diag}(\Lambda_{j,-},\Lambda_{j,+})$
with $\Lambda_{j,\pm}=V_j^{(1)}\pm A_{1,j}^{(1)}$, and  $P_j=I_2$ if $A_{1,j}^{(1)}=0$ and otherwise
\be\label{eq:Gjn}
 P_j=P^{(0)}:=\begin{pmatrix}
\frac{1}{\sqrt{2}} &\frac{1}{\sqrt{2}}\\
-\frac{1}{\sqrt{2}} &\frac{1}{\sqrt{2}}\end{pmatrix}.
\ee

\begin{rmk}
Again, if the definite integrals in $\int_{t_n}^{t_{n+1}}\Lambda(t,x_j)\,dt$ cannot be evaluated analytically, we can
evaluate them numerically via the Simpson's quadrature rule as
\begin{align*}&\int_{t_n}^{t_{n+1}}A_1(t,x_j)\,dt\approx
\frac{\tau}{6}\left[A_1(t_n,x_j)+4A_1\left(t_n+\frac{\tau}{2},x_j\right)+A_1(t_{n+1},x_j)\right],\\
&\int_{t_n}^{t_{n+1}}V(t,x_j)\,dt\approx
\frac{\tau}{6}\left[V(t_n,x_j)+4V\left(t_n+\frac{\tau}{2},x_j\right)+V(t_{n+1},x_j)\right].
\end{align*}
\end{rmk}

\begin{lemma}
\label{Mass_STA_TSFP}
The  TSFP (\ref{eq:tsfp}) conserves
the mass in the discretized level, i.e.
\be
\|\Phi^n\|_{l^2}^2:=h\sum_{j=0}^{M-1}|\Phi_j^{n}|^2\equiv  h\sum_{j=0}^{M-1}|\Phi_j^{0}|^2
=\|\Phi^0\|_{l^2}^2=h\sum_{j=0}^{M-1}|\Phi_0(x_j)|^2,\qquad n\ge0.
\ee
\end{lemma}
{\sl Proof:} The proof is quite standard and similar to that of Lemma \ref{Mass_STA_FDTD}. We omit it here.
$\hfill\Box$

From Lemma \ref{Mass_STA_TSFP}, we conclude that the TSFP (\ref{eq:tsfp}) is unconditionally stable. In addition,
under proper assumptions of the exact solution $\Phi(t,x)$ and electromagnetic potentials, it is easy to show the following error estimate
via the formal Lie calculus introduced in \cite{Lubich},
\be\label{tssp987}
\|\Phi(t_n,x)-I_M(\Phi^n)\|_{L^2}\lesssim h^{m_0}+\frac{\tau^2}{\eps^4}, \qquad
0\le n\le \frac{T}{\tau},
\ee
where $m_0$ depends on the regularity of $\Phi(t,x)$. We omit the details here for brevity.

\section{Numerical comparison and applications} \label{sec5}

In this section, we compare the accuracy of
different numerical methods including the FDTD, EWI-FP and TSFP methods
for the Dirac equation (\ref{SDEdd}) in 1D in terms of the mesh size $h$, time step $\tau$
and the parameter $0<\varepsilon\le 1$. We will pay particular attention
to the $\eps$-scalability of different methods
in the nonrelativistic limit regime, i.e. $0<\varepsilon\ll 1$.
Then we simulate the dynamics of the Dirac equation (\ref{SDEdd}) in 2D with a honeycomb lattice potential
by the TSFP method.

\subsection{Comparison of spatial/temporal resolution}
To test the accuracy, we choose the electromagnetic potentials in the Dirac equation (\ref{SDEdd}) with
$d=1$ as
\begin{eqnarray}
A_1(t,x) = \frac{(x+1)^2}{1+x^2},\qquad V(t,x) = \frac{1-x}{1+x^2}, \qquad x\in{\mathbb R}, \quad t\ge0,
\end{eqnarray}
and the initial data as
\begin{equation}
\phi_1(0,x) = e^{-x^2/2}, \quad\phi_2(0,x) = e^{-(x-1)^2/2}, \qquad x\in{\mathbb  R}.
\end{equation}
The problem is solved numerically on an interval $\Omega=(-16, 16)$
with periodic boundary conditions on $\partial\Omega$.
The `reference' solution  $\Phi(t,x)=(\phi_1(t, x),\phi_2(t, x))^T$
is obtained numerically by using the TSFP method with a small time step and a very fine mesh size,
e.g. $\tau_e = 10^{-7}$ and $h_e = 1/16$ or $h_e = 1/4096$ for the comparison of
the EWI-FP/TSFP methods or the FDTD methods, respectively.
Denote $\Phi^n_{h,\tau}$ as the numerical solution
obtained by a numerical method with  mesh size $h$ and time step $\tau$.
In order to quantify the convergence, we introduce
\[e_{h,\tau}(t_n)=\|\Phi^n-\Phi(t_n,\cdot)\|_{l^2}=\sqrt{h\sum_{j=0}^{M-1}|\Phi^n_j-\Phi(t_n,x_j)|^2}.\]

\begin{table}[t!]
\def\temptablewidth{1\textwidth}
\vspace{-12pt}
\caption{Spatial and temporal error analysis of
the LFFD method for the Dirac equation (\ref{SDEdd}) in 1D.}
{\rule{\temptablewidth}{1pt}}

\begin{tabular*}{\temptablewidth}{@{\extracolsep{\fill}}cccccc}
Spatial Errors  & $h_0=1/8$   & $h_0/2$   &$h_0/2^2$ &  $h_0/2^3$ &  $h_0/2^4$ \\
\hline
$\varepsilon_0=1$ &  1.06E-1  &  2.65E-2  & 6.58E-3  &  1.64E-3 & 4.10E-4 \\
{\rm order}             &	--	&	2.00	&	2.01	&	2.00	& 2.00\\
\hline
$\varepsilon_0/2$ &  9.06E-2  &  2.26E-2  &  5.64E-3  &  1.41E-3 & 3.51E-4 \\
{\rm order}             &	--	&	2.00	&	2.00	&	2.00	& 2.00\\
\hline
$\varepsilon_0/2^2$ &  8.03E-2  &  2.02E-2  &  5.04E-3  &  1.25E-3 &  3.05E-4 \\
{\rm order}             &	--	&	1.99	&	2.00	&	2.01	& 2.02\\
\hline
$\varepsilon_0/2^3$ &  9.89E-2  &  2.47E-2  & 6.17E-3  &  1.54E-3 & 3.85E-4 \\
{\rm order}             &	--	&	2.00	&	2.00	&	2.00	& 2.00\\
\hline
$\varepsilon_0/2^4$ &  9.87E-2  &  2.48E-2  & 6.18E-3  &  1.54E-3 & 3.83E-4 \\
{\rm order}             &	--	&	2.00	&	2.00	&	2.00	& 2.01\\
\end{tabular*}
{\rule{\temptablewidth}{1pt}}
\begin{tabular*}{\temptablewidth}{@{\extracolsep{\fill}}cccccc}
Temporal Errors & $\begin{array}{c} \tau_0=0.1\\ h_0=1/8\\ \end{array}$  & $\begin{array}{c} \tau_0/8\\ h_0/8\delta_1(\vep)\\ \end{array}$  &  $\begin{array}{c} \tau_0/8^2\\ h_0/8^2\delta_2(\vep)\\ \end{array}$   & $\begin{array}{c} \tau_0/8^3\\  h_0/8^3\delta_3(\vep)\\ \end{array}$  & $\begin{array}{c} \tau_0/8^4\\ h_0/
8^4\delta_4(\vep)\\ \end{array}$  \\
\hline
  $\varepsilon_0=1$  &  \underline{1.38E-1}  &  1.99E-3  &  3.11E-5 & 4.86E-7 & 7.59E-9\\
  {\rm order}             &	--	&	2.04	&	2.00	&	2.00	& 2.00\\
\hline
  $\varepsilon_0/2$  &  unstable  &  \underline{1.14E-2}  &  1.77E-4 & 2.77E-6 & 4.32E-8\\
  {\rm order}             &	--	&	--	&	2.00	&	2.00	&2.00\\
\hline
  $\varepsilon_0/2^2$  &  unstable  &  4.59E-1  & \underline{7.01E-3} & 1.05E-4 & 1.64E-6\\
  {\rm order}             &	--	&	--	&	2.01	&	2.02	& 2.00\\
\hline
  $\varepsilon_0/2^3$  &  unstable  &  unstable  &  4.14E-1 & \underline{6.42E-3} & 1.00E-4\\
  {\rm order}             &	--	&		--&		--&	2.00	& 2.00	\\
\hline
  $\varepsilon_0/2^4$  &  unstable  &  unstable  &  unstable & 4.04E-1 & \underline{6.00E-3}\\
   {\rm order}             &	--	&	--	&	--	&		--&	2.02\\
\end{tabular*}
{\rule{\temptablewidth}{1pt}}
\label{table_linear_LFFD}
\end{table}

\begin{table}[t!]
\def\temptablewidth{1\textwidth}
\vspace{-12pt}
\caption{Spatial and temporal error analysis of
the SIFD1 method for the Dirac equation (\ref{SDEdd}) in 1D. }
{\rule{\temptablewidth}{1pt}}
\begin{tabular*}{\temptablewidth}{@{\extracolsep{\fill}}cccccc}
Spatial Errors & $h_0=1/8$   & $h_0$/2   &$h_0/2^2$     &  $h_0/2^3$ &  $h_0/2^4$   \\
\hline
  $\varepsilon_0=1$ &  1.06E-1  &  2.65E-2  & 6.58E-3  &  1.64E-3 & 4.10E-4 \\
{\rm order} & --	& 2.00	&	2.01 &	2.00 &	2.00 \\
\hline
$\varepsilon_0/2$ &  9.06E-2  &  2.26E-2  &  5.64E-3  &  1.41E-3 & 3.51E-4 \\
{\rm order} & --	&2.00	&2.00	&	2.00&2.00	\\
\hline
$\varepsilon_0/2^2$ &  8.03E-2  &  2.02E-2  &  5.04E-3  &  1.25E-3 &  3.05E-4 \\
{\rm order} & --	& 1.99	&2.00	&2.01	&2.02	\\
\hline
$\varepsilon_0/2^3$ &  9.89E-2  &  2.47E-2  & 6.17E-3  &  1.54E-3 & 3.85E-4 \\
{\rm order} & --	&2.00	&2.00	&	2.00 &2.00	\\
\hline
$\varepsilon_0/2^4$ &  9.87E-2  &  2.48E-2  & 6.18E-3  &  1.54E-3 & 3.83E-4 \\
{\rm order} & --	&	1.99&	2.00	& 2.00	&	2.01\\
\end{tabular*}
{\rule{\temptablewidth}{1pt}}
\begin{tabular*}{\temptablewidth}{@{\extracolsep{\fill}}cccccc}
Temporal Errors & $\begin{array}{c} \tau_0=0.1\\ h_0=1/8\\ \end{array}$  & $\begin{array}{c} \tau_0/8\\ h_0/8\delta_1(\vep)\\ \end{array}$  &  $\begin{array}{c} \tau_0/8^2\\ h_0/8^2\delta_2(\vep)\\ \end{array}$   & $\begin{array}{c} \tau_0/8^3\\ h_0/8^3\delta_3(\vep)\\ \end{array}$  & $\begin{array}{c} \tau_0/8^4\\ h_0/
8^4\delta_4(\vep)\\ \end{array}$ \\
\hline
  $\varepsilon_0=1$  &  \underline{1.44E-1}  &  2.09E-3  &  3.27E-5 & 5.11E-7 & 7.98E-9\\
  {\rm order} & --	&	2.03	&	2.00	&	2.00	 &	2.00	\\
\hline
  $\varepsilon_0/2$  &  unstable  &  \underline{2.99E-2}  &  4.67E-4 & 7.30E-6 & 1.14E-7\\
    {\rm order} & --	&	--	&	2.00	&	2.00	 &	2.00	\\
\hline
  $\varepsilon_0/2^2$  &  unstable  &  8.18E-1  & \underline{1.54E-2} & 2.41E-4 & 3.77E-6\\
    {\rm order} & --	&	--	&	1.91	&	2.00	 &	2.00	\\
\hline
  $\varepsilon_0/2^3$  &  unstable  &  unstable  &  7.99E-1 & \underline{1.31E-2} & 2.05E-4\\
    {\rm order} & --	&	--	&	--	&	1.98	 &	2.00	\\
\hline
  $\varepsilon_0/2^4$  &  unstable  &  unstable  &  4.19E-1 & 7.97E-1 & \underline{1.26E-2}\\
    {\rm order} & --	&	--	&	--	&	-0.31	 &	1.99	\\
\end{tabular*}
{\rule{\temptablewidth}{1pt}}
\label{table_linear_SIFD1}
\end{table}

Table \ref{table_linear_LFFD} lists spatial errors $e_{h,\tau_e}(t=2)$ with different $h$
(upper part) and temporal errors $e_{h_e,\tau}(t=2)$ with different $\tau$
(lower part) for the LFFD method (\ref{efd1}).  Tables \ref{table_linear_SIFD1}-\ref{table_linear_TSSP}
show similar results for the SIFD1 method (\ref{sifd1}), SIFD2 method (\ref{sifd2}), CNFD method (\ref{cnfd1}),
EWI-FP method (\ref{AL4})-(\ref{AL5}) and TSFP method (\ref{eq:tsfp}), respectively.
For the LFFD and SIFD1 methods,
due to the stability condition and accuracy requirement, we take
\[\delta_j(\vep) = \left\{\begin{array}{ll} \vep^2 &\vep_0/2^j\le \vep \le 1,\\
\vep_0^2/4^j &0<\vep<\vep_0/2^j,\\
 \end{array} \right.  \qquad j=0,1,\ldots\]
in Tables \ref{table_linear_LFFD} and \ref{table_linear_SIFD1}.
For comparison, Table  \ref{table_linear_comp}
depicts temporal errors
of different numerical methods when $\vep=1$ for different $\tau$,
Table  \ref{table_linear_scal} depicts temporal errors of different numerical methods
under different $\eps$-scalability.

\begin{table}[t!]
\def\temptablewidth{1\textwidth}
\vspace{-12pt}
\caption{Spatial and temporal error analysis of
the SIFD2 method for the Dirac equation (\ref{SDEdd}) in 1D. }
{\rule{\temptablewidth}{1pt}}
\begin{tabular*}{\temptablewidth}{@{\extracolsep{\fill}}cccccc}
Spatial Errors & $h_0=1/8$ & $h_0/2$ & $h_0/2^2$ & $h_0/2^3$ & $h_0/2^4$ \\
\hline
 $\varepsilon_0=1$ &  1.06E-1  &  2.65E-2  & 6.58E-3  &  1.64E-3 & 4.10E-4 \\
 {\rm order} &	 --	&	2.00	&	2.01	&	2.00	&	2.00			\\
 \hline
$\varepsilon_0/2$ &  9.06E-2  &  2.26E-2  &  5.64E-3  &  1.41E-3 & 3.51E-4 \\
 {\rm order} &	--	&	2.00	&	2.00	&	2.00	&		2.00		\\
 \hline
$\varepsilon_0/2^2$ &  8.03E-2  &  2.02E-2  &  5.04E-3  &  1.25E-3 &  3.05E-4 \\
 {\rm order} &	--	&	1.99	&	2.00	&	2.01	&	2.02			\\
 \hline
$\varepsilon_0/2^3$ &  9.89E-2  &  2.47E-2  & 6.17E-3  &  1.54E-3 & 3.85E-4 \\
 {\rm order} &	--	&	2.00	&	2.00	&	2.00	&		2.00		\\
 \hline
$\varepsilon_0/2^4$ &  9.87E-2  &  2.48E-2  & 6.18E-3  &  1.54E-3 & 3.83E-4 \\
 {\rm order} &	--	&	1.99	&	2.00	&	2.00	&	2.01			\\
 \end{tabular*}
{\rule{\temptablewidth}{1pt}}
\begin{tabular*}{\temptablewidth}{@{\extracolsep{\fill}}cccccc}
Temporal Errors & $\tau_0$=0.1  & $\tau_0$/8 & $\tau_0/8^2$ & $\tau_0/8^3$ & $\tau_0/8^4$ \\
\hline
  $\varepsilon_0=1$  &  \underline{1.72E-1}  &  2.59E-3  &  4.05E-5 & 6.33E-7 & 9.89E-9\\
   {\rm order} &	--	&	2.01	&    2.00		&	2.00		&		2.00			\\
 \hline
  $\varepsilon_0/2$  &  1.69  &  \underline{3.57E-2}  &  5.58E-4 & 8.72E-6 & 1.36E-7\\
   {\rm order} &	--	&	1.86	&	 2.00		&	2.00		&		2.00			\\
 \hline
  $\varepsilon_0/2^2$  &  2.59  &  8.66E-1  & \underline{1.63E-2} & 2.55E-4 & 3.98E-6\\
   {\rm order} &	--	&	0.52	&	1.91	&	2.00	&	2.00			\\
 \hline
  $\varepsilon_0/2^3$  &  2.67  &  2.89  &  8.43E-1 & \underline{1.37E-2} & 2.14E-4\\
   {\rm order} &	--	&	-0.04	&	0.59	&	1.98	&	2.00			\\
 \hline
  $\varepsilon_0/2^4$  &  3.07  &  3.56  &  5.19E-1 & 8.37E-1 & \underline{1.28E-2}\\
   {\rm order} &	--	&	-0.07	&	0.93	&	-0.23	&	2.01			\\
\end{tabular*}
{\rule{\temptablewidth}{1pt}}
\label{table_linear_SIFD2}
\end{table}

\begin{table}[t!]
\def\temptablewidth{1\textwidth}
\vspace{-12pt}
\caption{Spatial and temporal error analysis of
the CNFD method for the Dirac equation (\ref{SDEdd}) in 1D.}
{\rule{\temptablewidth}{1pt}}
\begin{tabular*}{\temptablewidth}{@{\extracolsep{\fill}}cccccc}
Spatial Errors & $h_0$=1/8   & $h_0/2$   &$h_0/2^2$     &  $h_0/2^3$  &  $h_0/2^3$  \\
\hline
$\varepsilon_0=1$ &  1.06E-1  &  2.65E-2  & 6.58E-3  &  1.64E-3 & 4.10E-4 \\
 {\rm order} &	--	&	2.00	&	2.01	&	2.00	&	2.00			\\
 \hline
$\varepsilon_0/2$ &  9.06E-2  &  2.26E-2  &  5.64E-3  &  1.41E-3 & 3.51E-4 \\
 {\rm order} &	--	&	2.00	&	2.00	&	2.00	&	2.00		\\
 \hline
$\varepsilon_0/2^2$ &  8.03E-2  &  2.02E-2  &  5.04E-3  &  1.25E-3 &  3.05E-4 \\
 {\rm order} &	--	&	1.99	&	2.00	&	2.01	&	2.02			\\
 \hline
$\varepsilon_0/2^3$ &  9.89E-2  &  2.47E-2  & 6.17E-3  &  1.54E-3 & 3.85E-4 \\
 {\rm order} &	--	&	2.00	&	2.00	&	2.00	&	2.00		\\
 \hline
$\varepsilon_0/2^4$ &  9.87E-2  &  2.48E-2  & 6.18E-3  &  1.54E-3 & 3.83E-4 \\
 {\rm order} &	--	&	1.99	&	2.00	&	2.00	&	2.01			\\
 \end{tabular*}
{\rule{\temptablewidth}{1pt}}
\begin{tabular*}{\temptablewidth}{@{\extracolsep{\fill}}cccccc}
Temporal Errors & $\tau_0$=0.1  & $\tau_0$/8 & $\tau_0/8^2$ & $\tau_0/8^3$  & $\tau_0/8^4$ \\
\hline
  $\varepsilon_0=1$  &  \underline{5.48E-2}  &  8.56E-4  &  1.34E-5 & 2.09E-7 & 3.27E-9\\
   {\rm order} &	--	&	2.00	&	2.00	&	2.00	&	2.00			\\
 \hline
  $\varepsilon_0/2$  &  3.90E-1  &  \underline{6.63E-3}  &  1.77E-4 & 2.77E-6 & 4.32E-8\\
   {\rm order} &	--	&	1.96	&	1.74	&	2.00	&	2.00			\\
 \hline
  $\varepsilon_0/2^2$  &  1.79  &  2.27E-1  & \underline{3.55E-3} & 1.56E-5 & 2.44E-7\\
   {\rm order} &	--	&	0.99	&	2.00	&	2.61	&	2.00			\\
 \hline
  $\varepsilon_0/2^3$  &  3.10  &  4.69E-1  &  2.06E-1 & \underline{3.22E-3} & 5.03E-5\\
   {\rm order} &	--	&	0.91	&	0.40	&	2.00	&	2.00			\\
 \hline
  $\varepsilon_0/2^4$  &  2.34  &  1.83  &  8.05E-1 & 2.04E-1 & \underline{3.19E-3}\\
   {\rm order} &	--	&	0.12	&	0.39	&	0.66	&	2.00			\\
\end{tabular*}
{\rule{\temptablewidth}{1pt}}
\label{table_linear_FDTDg}
\end{table}

\begin{table}[t!]
\def\temptablewidth{1\textwidth}
\vspace{-12pt}
\caption{Spatial and temporal error analysis of
the EWI-FP method for the Dirac equation (\ref{SDEdd}) in 1D.}
{\rule{\temptablewidth}{1pt}}
\begin{tabular*}{\temptablewidth}{@{\extracolsep{\fill}}cccccc}
Spatial Errors & $h_0$=2   & $h_0$/2   &$h_0/2^2$ &  $h_0/2^3$  &  $h_0/2^4$ \\
\hline
  $\varepsilon_0=1$  &  1.10 & 2.43E-1 & 2.99E-3 & 2.79E-6 & 1.00E-8 \\
   {\rm order} &	--	&	2.13	&	9.02	&	32.74	&	16.70	\\
 \hline
  $\varepsilon_0/2$  &  1.06 & 1.46E-1 & 1.34E-3 & 9.61E-7 & 5.90E-9 \\
   {\rm order} &	--	&	2.69	&	10.44	&	37.34	&	12.76	\\
 \hline
  $\varepsilon_0/2^2$  &  1.11 & 1.43E-1 & 9.40E-4 & 5.10E-7 & 7.02E-9 \\
   {\rm order} &	--	&	2.79	&	12.33	&	42.93	&	8.52			\\
 \hline
  $\varepsilon_0/2^3$  &  1.15 & 1.44E-1 & 7.89E-4 & 3.62E-7 & 6.86E-9 \\
   {\rm order} &	--	&	2.83	&	13.51	&	46.69	&	7.26			\\
 \hline
  $\varepsilon_0/2^4$  &  1.18 & 1.45E-1 & 7.63E-4 & 2.91E-7 & 8.46E-9 \\
   {\rm order} &	--	&	2.85	&	13.79	&	51.21	&	5.86			\\
   \end{tabular*}
{\rule{\temptablewidth}{1pt}}
\begin{tabular*}{\temptablewidth}{@{\extracolsep{\fill}}cccccc}
Temporal Errors & $\tau_0$=0.1  & $\tau_0$/4 & $\tau_0/4^2$ & $\tau_0/4^3$ & $\tau_0/4^4$\\
\hline
  $\varepsilon_0 = 1$  &  \underline{1.40E-1}  &  8.51E-3  &  5.33E-4  &  3.34E-5  &  2.09E-6\\
   {\rm order} &	--	&	2.02	&	2.00	&	2.00	&	2.00			\\
 \hline
  $\varepsilon_0/2    $  &  4.11E-1  &  \underline{2.37E-2}  &  1.49E-3  &  9.29E-5  &  5.81E-6\\
   {\rm order} &	--	&	2.06	&	2.00	& 2.00		&	2.00			\\
 \hline
  $\varepsilon_0/2^2  $  &  6.03  &  1.88E-1  &  \underline{1.18E-2}  &  7.38E-4  &  4.62E-5\\
   {\rm order} &	--	&	2.50	&	2.00	& 2.00		&	2.00				\\
 \hline
  $\varepsilon_0/2^3  $  &  2.21  &  3.98  &  1.60E-1  &  \underline{1.01E-2}  &  6.31E-4\\
   {\rm order} &	--	&	-0.42	&	2.32	&	2.00	&		2.00		\\
 \hline
  $\varepsilon_0/2^4  $  &  2.16  &  2.09  &  3.58  &  1.53E-1  &  \underline{9.69E-3}\\
   {\rm order} &	--	&	0.02	&	-0.39	&	2.27	&	1.99			\\
\end{tabular*}
{\rule{\temptablewidth}{1pt}}
\label{table_linear_EWI}
\end{table}

\begin{table}[t!]
\def\temptablewidth{1\textwidth}
\vspace{-12pt}
\caption{Spatial and temporal error analysis of
the TSFP method for the Dirac equation (\ref{SDEdd}) in 1D.}
{\rule{\temptablewidth}{1pt}}
\begin{tabular*}{\temptablewidth}{@{\extracolsep{\fill}}cccccc}
Spatial Errors & $h_0=2$ & $h_0/2$ &$h_0/2^2$  &  $h_0/2^3$ & $h_0/2^4$ \\
\hline
  $\varepsilon_0=1$  &  1.10  &  2.43E-1  &  2.99E-3  &  2.79E-6 & 9.45E-9\\
   {\rm order} &	--	&	2.13	&	9.01	&	32.74	&	17.18			\\
 \hline
  $\varepsilon_0/2$  &  1.06  &  1.46E-1  &  1.34E-3 &  9.61E-7 & 5.57E-9\\
   {\rm order} &	--	&	2.69	&	10.44	&	37.34	&	13.14			\\
 \hline
  $\varepsilon_0/2^2$  &  1.11  &  1.43E-1  &  9.40E-4  &  5.10E-7 & 6.50E-9\\
   {\rm order} &	--	&	2.79	&	12.33	&	42.93	&	8.86			\\
 \hline
  $\varepsilon_0/2^3$  &  1.15  &  1.44E-1  &  7.89E-4  &  3.62E-7 & 6.84E-9\\
   {\rm order} &	--	&	2.83	&	13.51	&	46.69	&	7.27			\\
 \hline
  $\varepsilon_0/2^4$  &  1.18  &  1.45E-1  &  7.62E-4  &  2.88E-7 & 7.49E-9\\
   {\rm order} &	--	&	2.85	&	13.79	&	51.44	&	6.20			\\
 \hline
  $\varepsilon_0/2^5$  &  1.19  &  1.46E-1  &  7.53E-4 &  2.59E-7  & 7.96E-9\\
   {\rm order} &	--	&	2.85	&	13.92	&	53.92	&	5.70			\\
 \hline
  $\varepsilon_0/2^6$  &  1.20  &  1.47E-1  &  7.49E-4  &  2.63E-7 & 6.90E-9\\
   {\rm order} &	--	&	2.86	&	14.01	&	53.37	&	6.17			\\
 \end{tabular*}
{\rule{\temptablewidth}{1pt}}
\begin{tabular*}{\temptablewidth}{@{\extracolsep{\fill}}cccccccc}
Temporal Errors & $\tau_0$=0.4  & $\tau_0/4$ & $\tau_0/4^2$ & $\tau_0/4^3$
& $\tau_0/4^4$ &$\tau_0/4^5$ & $\tau_0/4^6$ \\
\hline
 $\varepsilon_0=1$  & \underline{2.17E-1} & 1.32E-2 & 8.22E-4 & 5.13E-5 & 3.21E-6 & 2.01E-7 & 1.26E-8 \\
 {\rm order} &	--	& 2.02 	&	2.00	&	2.00	&	2.00	& 2.00	    & 2.00 		\\
 \hline
 $\varepsilon_0/2$  & 1.32 & \underline{6.60E-2} & 4.07E-3 & 2.54E-4 & 1.59E-5 & 9.92E-7 & 6.20E-8 \\
 {\rm order} &	--	&	2.16	&	2.00	&	2.00	&	2.00	& 2.00	    & 2.00 		\\
 \hline
 $\varepsilon_0/2^2$  & 2.50 & 3.33E-1 & \underline{1.68E-2} & 1.04E-3 & 6.49E-5 & 4.06E-6 & 2.54E-7  \\
 {\rm order} &	--	&	1.45	&	2.15	&	2.00	&	2.00	& 2.00	    & 2.00        	\\
 \hline
 $\varepsilon_0/2^3$  & 1.79 & 1.97 & 8.15E-2 & \underline{4.15E-3} & 2.57E-4 & 1.60E-5 & 1.00E-6 \\
 {\rm order} &	--	&	-0.07	&	2.30	&	2.14	&	2.01	& 2.00	    & 2.00			\\
 \hline
 $\varepsilon_0/2^4$  & 1.35 & 8.27E-1 & 8.85E-1 & 2.01E-2 & \underline{1.03E-3} & 6.35E-5 & 3.97E-6 \\
 {\rm order} &	--	&	0.35	&	-0.05	&	2.73	&	2.14	&	2.01	&		2.00		\\
 \hline
  $\varepsilon_0/2^5$  &  8.73E-1 & 2.25E-1 & 2.33E-1 & 2.49E-1 & 4.98E-3 & \underline{2.55E-4} & 1.58E-5 \\
   {\rm order} &	--	&	0.98 	&	-0.03 	&	-0.05 	&	 2.82	&	2.14 	&	2.01	 	\\
\end{tabular*}
{\rule{\temptablewidth}{1pt}}
\label{table_linear_TSSP}
\end{table}

\begin{table}[t!]
\def\temptablewidth{1\textwidth}
\vspace{-12pt}
\caption{Comparison of temporal errors of
different methods for the Dirac equation (\ref{SDEdd}) with $\vep = 1$.}
{\rule{\temptablewidth}{1pt}}
\begin{tabular*}{\temptablewidth}{@{\extracolsep{\fill}}ccccccc}
$\vep=1$ & $\tau_0$=0.1  & $\tau_0/4$ & $\tau_0/4^2$ & $\tau_0/4^3$
& $\tau_0/4^4$  & $\tau_0/4^5$\\
\hline
 LFFD  & 1.38E-1 & 8.00E-3 & 4.98E-4 & 3.11E-5 &  1.94E-6 & 1.21E-7 \\
 {\rm order} & 	--   & 	2.05	& 	2.00	& 	2.00	& 	2.00	&  2.00\\
 \hline
 SIFD1  & 1.44E-1 & 8.85E-3 & 5.53E-4 & 3.27E-5 & 2.16E-6 & 1.35E-7 \\
 {\rm order} & 	-- 	& 	2.01	& 2.00	& 2.04	& 1.96  &2.00\\
 \hline
 SIFD2  & 1.72E-1 & 1.17E-2 &  7.30E-4 & 4.05E-5 & 2.85E-6 & 1.78E-7 \\
 {\rm order} & 	-- 	& 	1.94	& 	 2.00	&  2.09		& 1.91	&2.00\\
 \hline
 CNFD  & 5.48E-2 & 3.49E-3 & 2.18E-4 & 1.34E-5  & 8.38E-7 & 5.23E-8  \\
 {\rm order} & 	-- 	& 	1.99	& 	2.00	& 	2.01	& 2.00	&2.00\\
 \hline
 EWI-FP & 1.40E-1 & 8.51E-3 & 5.33E-4 & 3.34E-5 & 2.09E-6 &  1.30E-7  \\
 {\rm order} & 	-- 	& 	2.02	& 2.00	& 	 2.00	& 	 2.00	& 2.00\\
 \hline
 TSFP  & 1.32E-2 & 8.22E-4 & 5.13E-5 & 3.21E-6 & 2.01E-7 & 1.26E-8  \\
 {\rm order} & 	-- 	& 	2.00	& 	2.00	& 	2.00	& 	2.00	& 2.00\\
\end{tabular*}
{\rule{\temptablewidth}{1pt}}
\label{table_linear_comp}
\end{table}

\begin{table}[t!]
\def\temptablewidth{1\textwidth}
\vspace{-12pt}
\caption{Comparison of temporal errors of different numerical methods
for the Dirac equation (\ref{SDEdd}) under proper $\varepsilon$-scalability. }
{\rule{\temptablewidth}{1pt}}
\begin{tabular*}{\temptablewidth}{@{\extracolsep{\fill}}cccccc}
$\ba{c}
\tau=O(\varepsilon^3)\\
\tau=O(h)\\
\ea$ &$\ba{c}
\varepsilon_0=1\\
h_0=1/8\\
\tau_0=0.1\\
\ea$ &$\ba{c}
\varepsilon_0/2\\
h_0/2\\
\tau_0/8\\
\ea$ &$\ba{c}
\varepsilon_0/2^2\\
h_0/2^2\\
\tau_0/8^2\\
\ea$ &$\ba{c}
\varepsilon_0/2^3\\
h_0/2^3\\
\tau_0/8^3\\
\ea$ &$\ba{c}
\varepsilon_0/2^4\\
h_0/2^4\\
\tau_0/8^4\\
\ea$  \\
\hline
LFFD &  1.38E-1  &  1.14E-2  & 7.01E-3  &  6.42E-3 & 6.00E-3 \\
SIFD1 &  1.44E-1  &  2.99E-2  &  1.54E-2  &  1.31E-2 &1.26E-2  \\
\hline
\hline
$\ba{c}
\tau=O(\varepsilon^3)\\
\ea$ &$\ba{c}
\varepsilon_0=1\\
\tau_0=0.1\\
\ea$ &$\ba{c}
\varepsilon_0/2\\
\tau_0/8\\
\ea$ &$\ba{c}
\varepsilon_0/2^2\\
\tau_0/8^2\\
\ea$ &$\ba{c}
\varepsilon_0/2^3\\
\tau_0/8^3\\
\ea$  &$\ba{c}
\varepsilon_0/2^4\\
\tau_0/8^4\\
\ea$   \\
\hline
SIFD2 &  1.72E-1  &  3.57E-2  & 1.63E-2  &  1.37E-2 & 1.28E-2  \\
CNFD &  5.48E-2  &  6.63E-3  & 3.55E-3  &  3.22E-3 & 3.19E-3 \\
\hline
\hline
 $\ba{c}
\tau=O(\varepsilon^2)\\
\ea$ &$\ba{c}
\varepsilon_0=1\\
\tau_0=0.1\\
\ea$ &$\ba{c}
\varepsilon_0/2\\
\tau_0/4\\
\ea$ &$\ba{c}
\varepsilon_0/2^2\\
\tau_0/4^2\\
\ea$ &$\ba{c}
\varepsilon_0/2^3\\
\tau_0/4^3\\
\ea$   &$\ba{c}
\varepsilon_0/2^4\\
\tau_0/4^4\\
\ea$     \\
\hline
EWI-FP  &  1.40E-1  &  2.37E-2  &  1.18E-2  &  1.01E-2 & 9.69E-3  \\
TSFP  &  1.32E-2  &  4.07E-3  &  1.04E-3  &  2.57E-4 & 6.35E-5 \\
\end{tabular*}
{\rule{\temptablewidth}{1pt}}
\label{table_linear_scal}
\end{table}

From Tables \ref{table_linear_LFFD}-\ref{table_linear_scal},
and additional numerical results  not shown here for brevity, we can
draw the following conclusions for the Dirac equation by using different numerical methods:

\bigskip

(i). For the discretization error in space,
for any fixed $\eps=\eps_0>0$, the FDTD methods
are second-order accurate, and resp.,
the EWI-FP and TSFP methods are spectrally accurate (cf. each row in the upper
parts of Tables \ref{table_linear_LFFD}-\ref{table_linear_TSSP} and Table \ref{table_linear_comp}).
For $0<\eps\le 1$, the errors are independent of $\eps$  for the  EWI-FP and TSFP methods
(cf. each column in the upper
parts of Tables \ref{table_linear_EWI}-\ref{table_linear_TSSP}),
 and resp., are almost independent of $\eps$ for the FDTD methods (cf. each column in the upper
parts of Tables \ref{table_linear_LFFD}-\ref{table_linear_FDTDg}).
In general, for any fixed $0<\eps\le 1$ and $h>0$,
the EWI-FP and TSFP methods perform much better than the FDTD methods
in spatial discretization.

(ii). For the discretization error in time, in the $O(1)$ speed-of-light regime,
i.e. $\eps=O(1)$, all the numerical methods including
FDTD, EWI-FP and TSFP are second-order accurate (cf. the first row in the lower
parts of Tables \ref{table_linear_LFFD}-\ref{table_linear_TSSP}).
In general, the EWI-FP and TSFP methods perform much better than the FDTD methods  in
temporal discretizations for a fixed time step.
In the non-relativistic limit  regime, i.e. $0<\eps\ll1$, for the FDTD methods, the `correct'
$\eps$-scalability is $\tau=O(\eps^3)$  which verifies our theoretical results; for
the EWI-FP and TSFP methods, the `correct'
$\eps$-scalability is $\tau=O(\eps^2)$ which again confirms our theoretical results.
In fact, for $0<\eps\le1$, one can observe clearly second-order convergence in time for
the FDTD methods only when $\tau\lesssim \eps^3$ (cf. upper triangles in the lower
parts of Tables \ref{table_linear_LFFD}-\ref{table_linear_FDTDg}), and resp., for the
EWI-FP and TSFP methods when $\tau\lesssim \eps^2$ (cf. upper triangles in the lower
parts of Tables \ref{table_linear_EWI}-\ref{table_linear_TSSP}).
In general, for any fixed $0<\eps\le 1$ and $\tau>0$,
the TSFP method performs the best, and the EWI-FP method
performs much better than the FDTD methods
in temporal discretization (cf. Table \ref{table_linear_scal}).

(iii). From Table \ref{table_linear_TSSP}, our numerical results suggest the following error bound
for the TSFP method when  $\tau\lesssim \eps^2$,
\be\label{notssp1}
 \|\Phi(t_n,\cdot)-I_M(\Phi^n)\|_{L^2}\lesssim h^{m_0}+\frac{\tau^2}{\eps^2},\qquad
 0\le n\le \frac{T}{\tau},
 \ee
which is much better than (\ref{tssp987}) for the TSFP method in the nonrelativistic limit regime.
Rigorous mathematical justification for (\ref{notssp1}) is on-going.

\bigskip

\begin{table}[t!]
\def\temptablewidth{1\textwidth}
\vspace{-12pt}
\caption{Spatial error analysis of the CNFD method for the free Dirac equation
with different  $h$.}
{\rule{\temptablewidth}{1pt}}
\begin{tabular*}{\temptablewidth}{@{\extracolsep{\fill}}cccccc}
$\varepsilon$ & $\varepsilon_0=1$ & $\varepsilon_0/2$ &$\varepsilon_0/2^2$  &  $\varepsilon_0/2^3$ & $\varepsilon_0/2^4$ \\
\hline
  $h_0=1/256$  &  1.61E-1  &  3.21E-1  &  6.35E-1  &  1.21 & 2.07\\
  $h_0/2$  &  4.03E-2  &  8.05E-2  &  1.59E-1  &  3.07E-1 & 5.43E-1\\
  $h_0/2^2$  &  1.01E-2  &  2.01E-2  &  3.99E-2 &  7.69E-2 & 1.36E-1\\
  $h_0/2^3$  &  2.52E-3  &  5.03E-3  &  9.97E-3  &  1.92E-2 & 3.41E-2\\
  $h_0/2^4$  &  6.30E-4  &  1.26E-3  &  2.47E-3  &  4.95E-3 & 8.64E-3\\
 \end{tabular*}
{\rule{\temptablewidth}{1pt}}
\label{table_linear_TDFD_spartial}
\end{table}

 From Tables \ref{table_linear_LFFD}-\ref{table_linear_FDTDg}, in the numerical example, we could not
observe numerically the $\eps$-dependence in the spatial discretization error for the FDTD methods, i.e.
$\frac{1}{\eps}$ in front of $h^2$, which was proven in Theorems \ref{thm_cnfd}-\ref{thm_sifd2}.
In order to investigate the spatial $\varepsilon$-resolution of the FDTD methods,
we consider the Dirac equation (\ref{dir}) on $\Omega=(-1, 1)$ with no electromagnetic potential
-- the free Dirac equation, i.e.
\begin{eqnarray}
A_1(t,x) \equiv 0,\qquad V(t,x) \equiv 0, \qquad x\in(-1,1), \quad t\ge0.
\end{eqnarray}
The initial data in (\ref{dir2}) is taken as
\begin{equation}
\phi_1(0,x) = e^{9\pi i(x+1)}, \qquad \phi_2(0,x) = e^{9\pi i(x+1)}, \qquad -1\le x\le 1.
\end{equation}
Table \ref{table_linear_TDFD_spartial} shows the spatial errors $e_{h,\tau_e}(t=2)$ of the CNFD method
with different $h$. The results for the LFFD, SIFD1 and SIFD2 methods are similar and
they are omitted here for brevity. From Table \ref{table_linear_TDFD_spartial}, we can conclude that the error bounds
 in the Theorems \ref{thm_cnfd}-\ref{thm_sifd2} are sharp.


\bigskip

Based on the above comparison, in view of both temporal and spatial accuracies
 and resolution capacity,  we conclude that
the EWI-FP  and TSFP methods perform much better than the FDTD methods
for the discretization of the Dirac equation, especially in the
nonrelativistic limit regime. For the reader's convenience,
we summarize the properties of different numerical methods in
Table \ref{table_linear_properties}.

\begin{table}[t!]
\def\temptablewidth{1\textwidth}
\vspace{-12pt}
\caption{Comparison of properties of different numerical methods for solving the Dirac equation
with $M$ being the number of grid points in space. }
{\rule{\temptablewidth}{1pt}}
{\small\begin{tabular*}{\temptablewidth}{@{\extracolsep{\fill}}ccccccc}
Method & LFFD  & SIFD1 & SIFD2 & CNFD & EWI-FP & TSFP \\
\hline
  Time symmetric &  Yes  &  Yes  & Yes  &  Yes  &  No & Yes \\
  Mass conservation & No  &  No  & No  &  Yes  &  No & Yes \\
  Energy conservation & No &  No  &  No &  Yes  & No & No \\
  Dispersion Relation & No &  No  & No  &  No  &  No & Yes \\ \hline
  Unconditionally stable & No &  No  & No  &  Yes  &  No & Yes \\
  Explicit scheme & Yes &  No  & No  &  No  &  Yes  & Yes \\
  Temporal accuracy & 2nd & 2nd & 2nd & 2nd & 2nd & 2nd \\
  Spatial accuracy & 2nd & 2nd & 2nd & 2nd & Spectral & Spectral \\
  Memory cost & $O(M)$ & $O(M)$ & $O(M)$ & $O(M)$ & $O(M)$ & $O(M)$ \\
  Computational cost & $O(M)$ & $O(M)$ & $O(M\ln M)$ & $\gg O(M)$ & $O(M\ln M)$ & $O(M\ln M)$ \\ \hline
  $\ba{c}
  \hbox{Resolution} \\
  \hbox{when}\, 0<\varepsilon\ll1\\
  \ea$ &$\ba{c}
   h=O(\sqrt{\varepsilon})\\
   \tau=O(\varepsilon^3)\\
   \ea$ &$\ba{c}
   h=O(\sqrt{\varepsilon})\\
   \tau=O(\varepsilon^3)\\
   \ea$  &$\ba{c}
   h=O(\sqrt{\varepsilon})\\
   \tau=O(\varepsilon^3)\\
   \ea$  &$\ba{c}
   h=O(\sqrt{\varepsilon})\\
   \tau=O(\varepsilon^3)\\
   \ea$  &$\ba{c}
   h=O(1)\\
   \tau=O(\varepsilon^2)\\
   \ea$  &$\ba{c}
   h=O(1)\\
   \tau=O(\varepsilon^2)\\
   \ea$   \\
   \end{tabular*}}
   {\rule{\temptablewidth}{1pt}}
\label{table_linear_properties}
\end{table}

\bigskip

As observed in \cite{BJP1,BJP2}, the time-splitting spectral (TSSP) method
for the Schr\"{o}dinger equation performs much better for the
physical observable, e.g. density and current, than for the wave function,
in the semiclassical limit regime with respect to the scaled Planck constnat $0<\vep\ll1$.
In order to see whether this is still valid for the TSFP method for
the Dirac equation in the nonrelativistic limit regime,
let $\rho^n=|\Phi^n_{h,\tau}|^2$, ${\bf J}^n=\frac{1}{\vep}(\Phi^n_{h,\tau})^*\sigma_1\Phi^n_{h,\tau}$
with $\Phi^n_{h,\tau}$  the numerical solution
obtained by the TSFP method with  mesh size $h$ and time step $\tau$, and define the errors
\begin{equation*}
e_{\rho}^{h,\tau}(t_n) := \|\rho^n-\rho(t_n,\cdot)\|_{l^1}=h\sum_{j=0}^{N-1}|\rho_j^n-\rho(t_n,x_j)|,
\quad  e_{{\bf J}}^{h,\tau}(t_n) := \|{\bf J}^n-{\bf J}(t_n,\cdot)\|_{l^1}=
h\sum_{j=0}^{N-1}|{\bf J}_j^n-{\bf J}(t_n,x_j)|.
\end{equation*}
Table \ref{table_linear_TSSPd} lists temporal errors $e_{\rho}^{h,\tau}(t=2)$ and
$e_{\bf J}^{h,\tau}(t=2)$ with different $\tau$
 for the TSFP method (\ref{eq:tsfp}). From this Table, we can see that
the approximations of the density and current are at the same order as for the wave function
by using the TSFP method. The reason that we can speculate is that $\rho=O(1)$ and ${\bf J} =O(\vep^{-1})$
(see details in (\ref{obser11}) or (\ref{obser12})) in
the Dirac equation, where in the Schr\"{o}dinger equation both density and current
are all at $O(1)$, when $0<\vep\ll1$.

\begin{table}[t!]
\def\temptablewidth{1\textwidth}
\vspace{-12pt}
\caption{Temporal errors for density and current of the TSFP for the Dirac equation (\ref{SDEdd}) in 1D.}
{\rule{\temptablewidth}{1pt}}
\begin{tabular*}{\temptablewidth}{@{\extracolsep{\fill}}cccccccc}
$e_{\rho}^{h,\tau}(t=2)$  & $\tau_0$=0.4  & $\tau_0/4$ & $\tau_0/4^2$ & $\tau_0/4^3$
& $\tau_0/4^4$ &$\tau_0/4^5$ & $\tau_0/4^6$ \\
\hline
 $\varepsilon_0=1$  & \underline{2.50E-1} & 1.54E-2 & 9.61E-4 & 6.01E-5 & 3.75E-6 & 2.34E-7 & 1.43E-8 \\
{\rm order}  &  	--&  	 2.01  & 2.00	&2.00	&2.00	&2.00	 	&2.01\\
\hline
 $\varepsilon_0/2$  & 1.22 & \underline{5.27E-2} & 3.21E-3 & 2.01E-4 & 1.25E-5 & 7.84E-7 & 4.92E-8 \\
{\rm order}  & 	 	--&   2.27  & 2.02	&2.00	&2.00	&2.00	 	&2.01\\
\hline
 $\varepsilon_0/2^2$  & 1.75 & 1.86E-1 & \underline{1.00E-2} & 6.20E-4 & 3.87E-5 & 2.42E-6 & 1.52E-7  \\
{\rm order}  & 	--&  1.62  & 2.11	&2.01	&2.00	&2.00	&2.00\\
\hline
 $\varepsilon_0/2^3$  & 1.11 & 1.39 & 2.95E-2 & \underline{1.53E-3} & 9.47E-5 & 5.92E-6 & 3.72E-7 \\
{\rm order}  & 	--&  -0.16  & 2.78	&2.13	&2.01	&2.00	&2.00\\
\hline
 $\varepsilon_0/2^4$  & 1.58 & 7.58E-1 & 7.81E-1 & 5.46E-3 & \underline{3.01E-4} & 1.87E-5 & 1.17E-6 \\
{\rm order}  & 	--&  0.53  & -0.02	&3.58	&2.09	&2.00	&2.00\\
\hline
 $\varepsilon_0/2^5$  & 9.59E-1 & 1.96E-1 & 2.29E-1 & 2.33E-1 & 1.20E-3 & \underline{6.76E-5} & 4.21E-6 \\
 {\rm order}  & 	--&  1.15  & -0.11	&-0.01	&3.8	&2.07	&2.00\\
\end{tabular*}
{\rule{\temptablewidth}{1pt}}
\begin{tabular*}{\temptablewidth}{@{\extracolsep{\fill}}cccccccc}
$e_{\bf J}^{h,\tau}(t=2)$ & $\tau_0$=0.4  & $\tau_0/4$ & $\tau_0/4^2$ & $\tau_0/4^3$
& $\tau_0/4^4$ &$\tau_0/4^5$ & $\tau_0/4^6$ \\
\hline
 $\varepsilon_0=1$  & \underline{1.70E-1} & 1.09E-2 & 6.83E-4 & 4.27E-5 & 2.67E-6 & 1.67E-7 & 1.02E-8 \\
 {\rm order}  & 	--&  1.98  & 2.00	&2.00	&2.00 	&2.00	&2.01  \\
\hline
 $\varepsilon_0/2$  & 9.15E-1 & \underline{6.39E-2} & 4.00E-3 & 2.50E-4 & 1.56E-5 & 9.76E-7 & 6.08E-8 \\
 {\rm order}  & 	--&  1.92  & 2.00	&2.00	&2.00 	&2.00	&2.00  \\
\hline
 $\varepsilon_0/2^2$  & 1.58 & 3.45E-1 & \underline{1.69E-2} & 1.04E-3 & 6.50E-5 & 4.06E-6 & 2.54E-7  \\
 {\rm order}  & 	--&  1.10  & 2.18	&2.01	&2.00 	&2.00	&2.00  \\
\hline
 $\varepsilon_0/2^3$  & 1.06 & 1.26 & 5.83E-2 & \underline{2.87E-3} & 1.76E-4 & 1.11E-5 & 6.94E-7 \\
 {\rm order}  & 	--&  -0.12 & 2.22	&2.17	&2.01 	&2.00	&2.00  \\
\hline
 $\varepsilon_0/2^4$  & 1.11 & 9.78E-1 & 1.05 & 2.28E-2 & \underline{1.18E-3} & 7.33E-5 & 4.58E-6 \\
 {\rm order}  & 	--&  0.09 &  -0.05	& 2.76	&2.13 	&2.00	&2.00  \\
\hline
 $\varepsilon_0/2^5$  & 4.98E-1 & 1.55E-1 & 2.22E-1 & 2.39E-1 & 4.04E-3 & \underline{2.09E-4} & 1.29E-5 \\
 {\rm order}  & 	--&  0.84 &  -0.30	& -0.05	&2.94	&2.13	&2.01 \\
\hline
\end{tabular*}
{\rule{\temptablewidth}{1pt}}
\label{table_linear_TSSPd}
\end{table}


\subsection{Dynamics of the Dirac equation in 2D}

Here we study numerically the dynamics of the Dirac equation
(\ref{SDEdd}) in 2D with a honeycomb lattice potential, i.e. we take
$d=2$ and
\be
A_1(t,\bx)=A_2(t,\bx)\equiv 0,\quad V(t,\bx)=\cos\left(\frac{4\pi}{\sqrt{3}}{\bf e}_1\cdot \bx\right)+\cos\left(\frac{4\pi}{\sqrt{3}}{\bf e}_2\cdot \bx\right)+\cos\left(\frac{4\pi}{\sqrt{3}}{\bf e}_3\cdot \bx\right),
\ee
with
\be
{\bf e}_1=(-1,0)^{T},\qquad {\bf e}_2=(1/2,\sqrt{3}/2)^T,\qquad {\bf e}_3=(1/2,-\sqrt{3}/2)^T.
\ee
The initial data in (\ref{SDEini}) is taken as
\begin{equation}
\phi_1(0,\bx)=e^{-\frac{x^2+y^2}{2}},\quad \phi_2(0,\bx)=e^{-\frac{(x-1)^2+y^2}{2}},
\quad \bx=(x,y)^T\in{\mathbb R}^2.
\end{equation}

\begin{figure}
\centerline{\psfig{figure=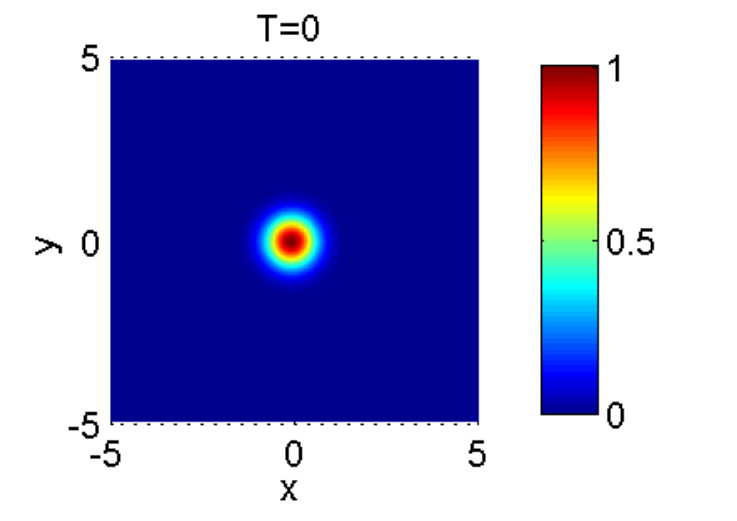,height=5cm,width=7cm,angle=0}
\psfig{figure=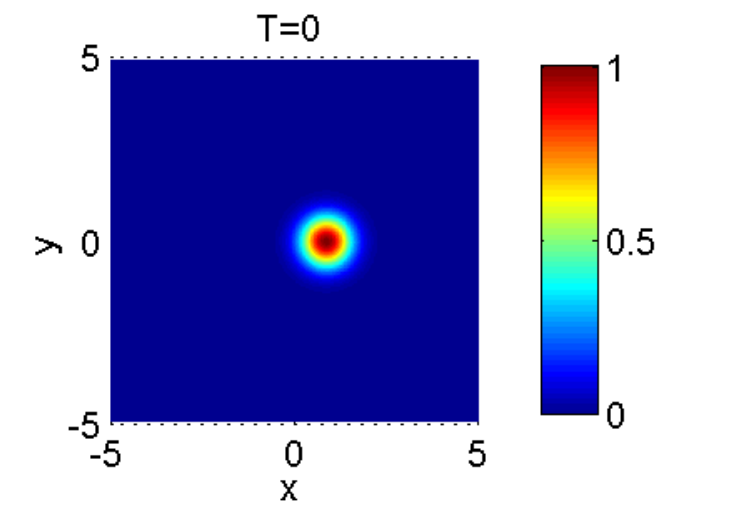,height=5cm,width=7cm,angle=0}}
\centerline{\psfig{figure=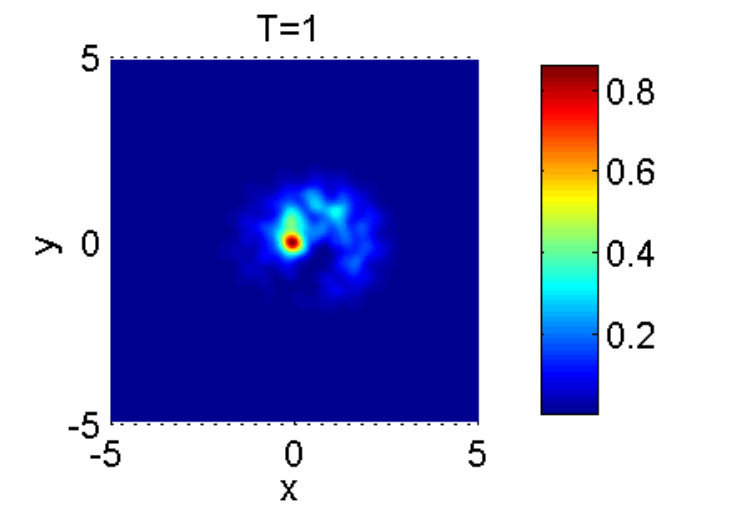,height=5cm,width=7cm,angle=0}
\psfig{figure=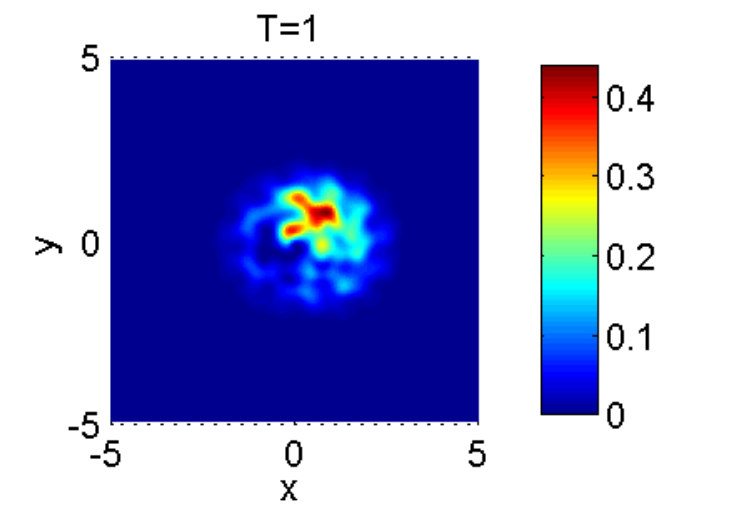,height=5cm,width=7cm,angle=0}}
\centerline{\psfig{figure=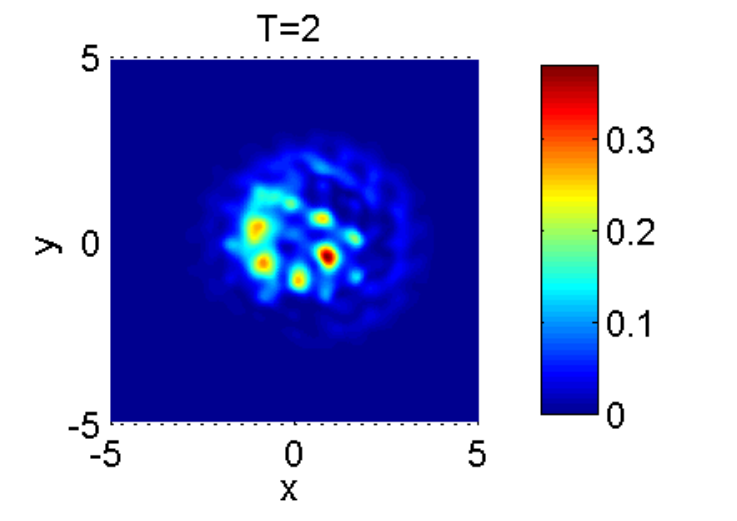,height=5cm,width=7cm,angle=0}
\psfig{figure=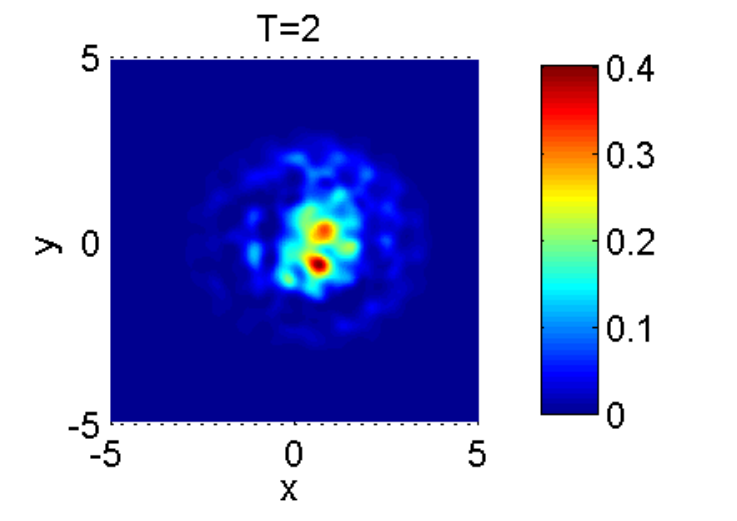,height=5cm,width=7cm,angle=0}}
\centerline{\psfig{figure=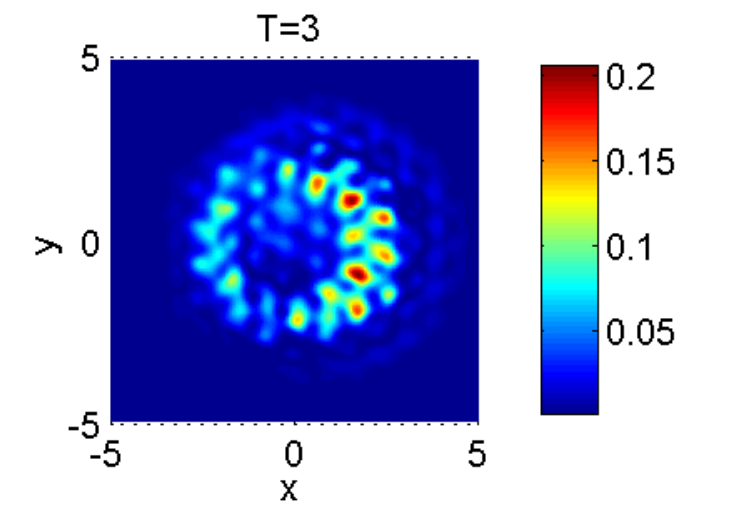,height=5cm,width=7cm,angle=0}
\psfig{figure=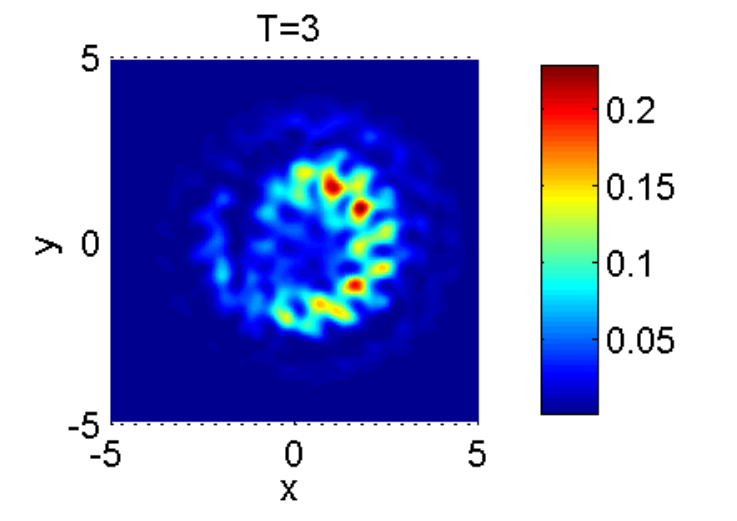,height=5cm,width=7cm,angle=0}}
\caption{Dynamics of the densities $\rho_1(t,\bx)=|\phi_1(t,\bx)|^2$(left) and $\rho_2(t,\bx)=|\phi_2(t,\bx)|^2$(right) of the Dirac equation in 2D with a honeycomb lattice potential when $\varepsilon=1$.}
\label{Dyn_HCP2_e1}
\end{figure}

\begin{figure}
\centerline{\psfig{figure=./Graphics/1stT0e1.pdf,height=5cm,width=7cm,angle=0}
\psfig{figure=./Graphics/2ndT0e1.pdf,height=5cm,width=7cm,angle=0}}
\centerline{\psfig{figure=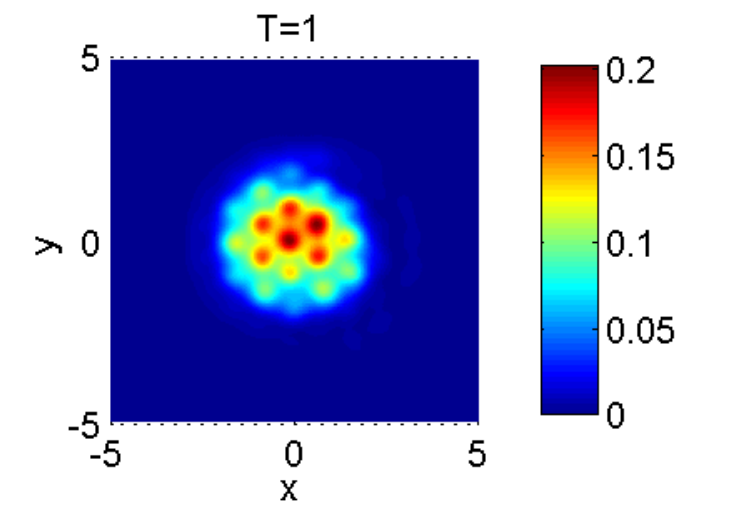,height=5cm,width=7cm,angle=0}
\psfig{figure=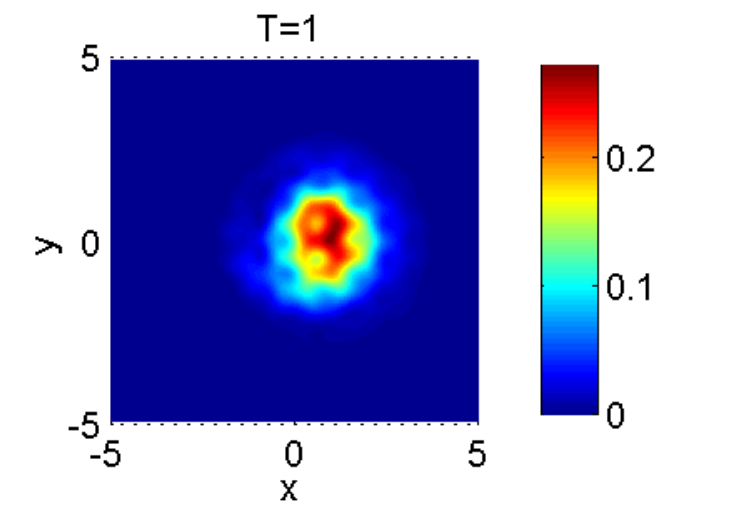,height=5cm,width=7cm,angle=0}}
\centerline{\psfig{figure=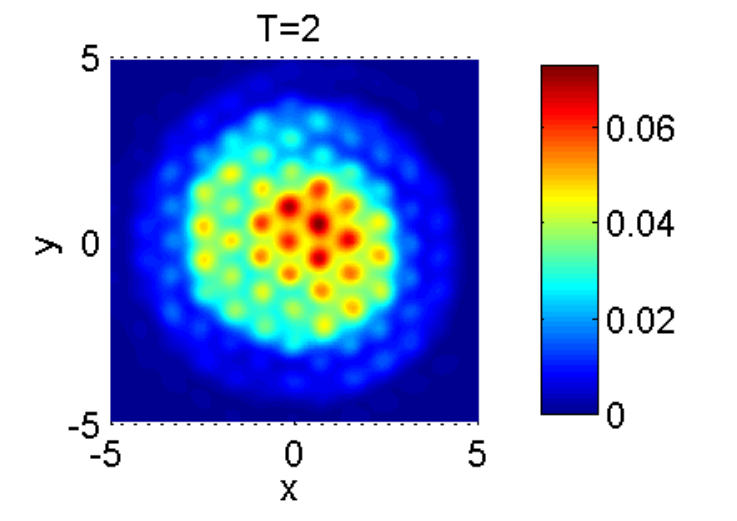,height=5cm,width=7cm,angle=0}
\psfig{figure=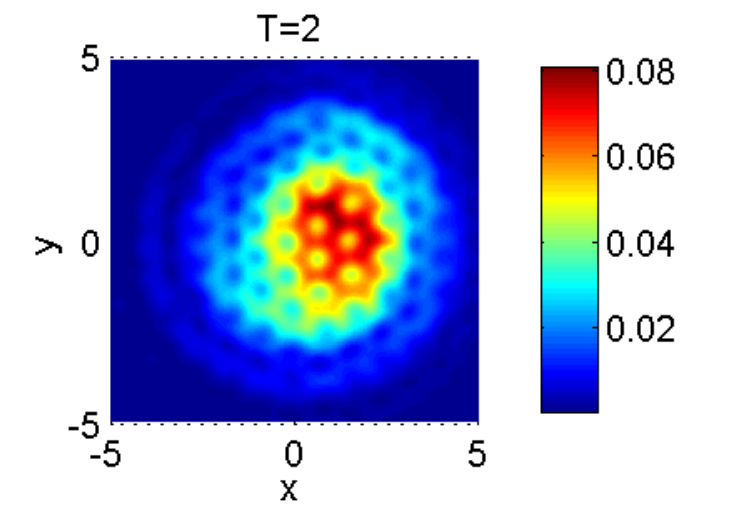,height=5cm,width=7cm,angle=0}}
\centerline{\psfig{figure=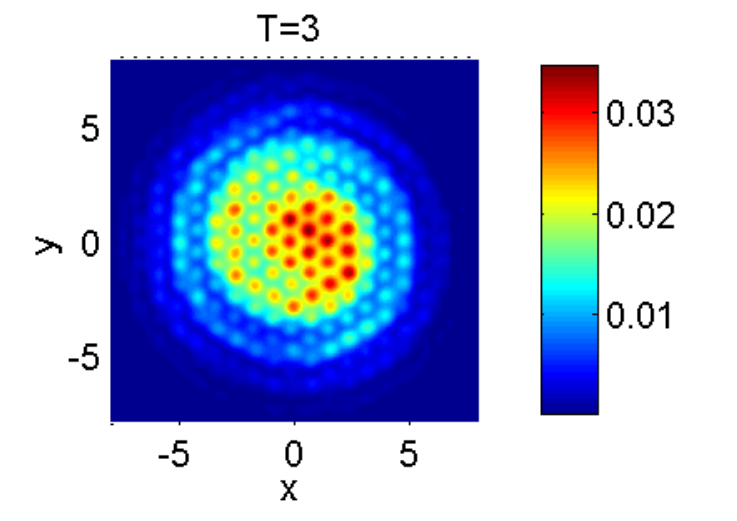,height=5cm,width=7cm,angle=0}
\psfig{figure=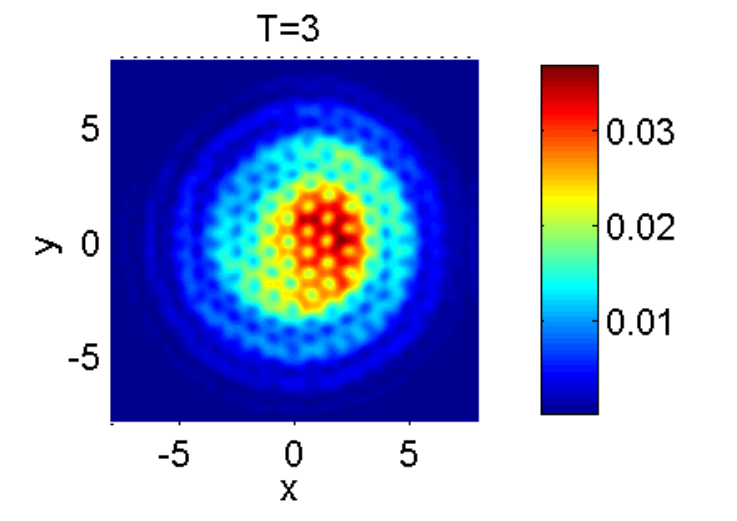,height=5cm,width=7cm,angle=0}}
\caption{Dynamics of the densities $\rho_1(t,\bx)=|\phi_1(t,\bx)|^2$(left) and $\rho_2(t,\bx)=|\phi_2(t,\bx)|^2$(right) of the Dirac equation in 2D with a honeycomb potential when $\varepsilon=0.2$.}
\label{Dyn_HCP2_e2}
\end{figure}

The problem is solved numerically on $\Omega=[-10,10]^2$ by the TSFP method with
mesh size $h=1/16$ and time step $\tau=0.01$. Figures \ref{Dyn_HCP2_e1} and \ref{Dyn_HCP2_e2}
depict the densities $\rho_j(t,\bx)=|\phi_j(t,\bx)|^2$ ($j=1,2$)
for $\varepsilon=1$ and $\varepsilon=0.2$, respectively.

From Figures \ref{Dyn_HCP2_e1}-\ref{Dyn_HCP2_e2},
we find that the dynamics of the Dirac equation depends significantly on
$\varepsilon$. In addition, the TSFP method can capture the dynamics
very accurately and efficiently.

\section{Conclusion}\label{sec6}
Three types of numerical methods based on different time integrations
were analyzed rigorously and compared numerically  for simulating the
Dirac equation in the nonrelativistic limit regime, i.e.  $0<\varepsilon\ll1$ or the speed of light goes to infinity.
The first class  consists of the second order standard FDTD methods, including energy conservative/
nonconservative and implicit/semi-implicit/explicit ones.
In the nonrelativistic limit regime, the error estimates of the FDTD methods were
rigorously analyzed, which suggest that the $\varepsilon$-scalability of the FDTD methods is
$\tau=O(\varepsilon^3)$ and  $h=O(\sqrt{\varepsilon})$. The second class applies the
Fourier spectral discretization in space and  Gautschi-type integration in time,
resulting in an EWI-FP method. Rigorous error bounds for the EWI-FP method were derived,
which show that  the $\varepsilon$-scalability of the EWI-FP method is  $\tau=O(\varepsilon^2)$
 and $h=O(1)$. The last class combines the Fourier spectral discretization in space and   splitting technique in time, which leads to a TSFP method. Based on the rigorous error analysis, the $\varepsilon$-scalability
 of the TSFP method is  $\tau=O(\varepsilon^2)$ and  $h=O(1)$, which is similar to the EWI-FP method.
 From the error analysis and numerical results, the EWI-FP and TSFP methods perform
 much better than the FDTD methods, especially in the nonrelativistic limit regime.
Extensive numerical results indicate that the TSFP method is superior than the EWI-FP in
terms of accuracy and efficiency, and thus the TSFP method is favorable for
solving the Dirac equation directly, especially in the nonrelativistic limit regime.
Finally, we studied the  dynamics of the Dirac equation in 2D with a honeycomb lattice potential
and observed some interesting dynamics for different $\eps$.

\setcounter{equation}{0}  

\begin{center}
{\bf Appendix A}. Proof of Theorem \ref{thm_cnfd} for the CNFD method
\end{center}
\setcounter{equation}{0}
\renewcommand{\theequation}{A.\arabic{equation}}

Define the local truncation error $\xi^{n}=(\xi_0^n,\xi_1^n,\ldots,\xi_M^n)^T\in X_{M}$
of the  CNFD (\ref{cnfd1}) with (\ref{bdc}) as
\begin{align}\label{Apdex1}
\xi_{j}^{n}:=&i\delta_{t}^{+}\Phi(t_n,x_{j})+
\frac{i}{\varepsilon}\sigma_1\frac{\delta_{x}\Phi(t_{n+1},x_{j})+\delta_{x}\Phi(t_n,x_{j})}{2}
-\frac{1}{\varepsilon^2}\sigma_3\frac{\Phi(t_{n+1},x_j)+\Phi(t_n,x_j)}{2}\nonumber\\
&+\left[A_{1}(t_{n+1/2},x_j)\sigma_1-V(t_{n+1/2},x_j)I_2\right]\frac{\Phi(t_{n+1},x_j)+\Phi(t_n,x_j)}{2},
\quad 0\le j\le M-1, \ n\ge0.\qquad
\end{align}
Applying the Taylor expansion in (\ref{Apdex1}), noticing (\ref{dir}) and the assumptions  (A) and (B),
 and using the triangle inequality, for $0<\eps\le 1$,  we obtain
\begin{align}
|\xi_{j}^n|\le&\frac{\tau^2}{24}\|\partial_{ttt}\Phi\|_{L^{\infty}(\overline{\Omega}_T)}+
\frac{h^2}{6\varepsilon}\|\partial_{xxx}\Phi\|_{L^{\infty}(\overline{\Omega}_T)}+
\frac{\tau^2}{8\varepsilon}\|\partial_{xtt}\Phi\|_{L^{\infty}(\overline{\Omega}_T)}
+\frac{\tau^2}{8}\left(\frac{1}{\eps^2}+V_{\rm max}+A_{1,\rm max}\right)
\|\partial_{tt}\Phi\|_{L^{\infty}(\overline{\Omega}_T)}\nonumber\\
\lesssim&\frac{\tau^2}{\varepsilon^6}+\frac{h^2}{\varepsilon}+\frac{\tau^2}{\varepsilon^5}
+\frac{\tau^2}{\varepsilon^4}\lesssim\frac{\tau^2}{\varepsilon^6}+\frac{h^2}{\varepsilon},
\qquad j=0,1,\ldots, M-1, \quad n\ge0,
\end{align}
which immediately implies
\be\label{locbd1}
\|\xi^n\|_{l^\infty}=\max_{0\le j\le M-1}|\xi_j^n|
\lesssim \frac{\tau^2}{\varepsilon^6}+\frac{h^2}{\varepsilon}, \qquad
\|\xi^n\|_{l^2}\lesssim \|\xi^n\|_{l^\infty} \lesssim
\frac{\tau^2}{\varepsilon^6}+\frac{h^2}{\varepsilon},
\qquad n\ge0, \qquad 0<\eps\le 1.
\ee
Subtracting (\ref{cnfd1}) from (\ref{Apdex1}), noticing (\ref{err}), we get for $n\ge0$
\begin{eqnarray}
\label{err1}
i\delta_t^+\bee_{j}^{n}=-\frac{i}{\varepsilon}\sigma_1\delta_x \bee_{j}^{n+1/2}+\frac{1}{\varepsilon^2}\sigma_3
\bee_{j}^{n+1/2}+\left(V_j^{n+1/2}I_2-A_{1,j}^{n+1/2}\sigma_1
\right)\bee_{j}^{n+1/2}+\xi_j^n, \qquad 0\le j\le M-1,
\end{eqnarray}
with $\bee_{j}^{n+1/2}=\frac{\bee_{j}^{n+1}+\bee_{j}^{n}}{2}$ for $j=0,1,\ldots,M$,
and the boundary and initial conditions are given as
\be\label{ibdcf}
\bee_0^{n}=\bee_{M}^{n}, \qquad \bee_{-1}^{n}=\bee_{M-1}^n, \qquad n\ge0,\qquad \qquad
\bee_j^0={\bf 0}, \qquad j=0,1,\ldots,M.
\ee
Similarly to the proof for Lemma \ref{Mass_STA_FDTD}, multiplying (\ref{err1}) from
the left by $h\tau\, \left(\bee_j^{n+1/2}\right)^*$, taking the imaginary part,
then summing for $j=0,1,\ldots,M-1$, using the triangle inequality and
Young's inequality, noticing (\ref{Paulim}), (\ref{locbd1}) and (\ref{ibdcf}), we get
\bea
\|\bee^{n+1}\|_{l^2}^2-\|\bee^{n}\|_{l^2}^2&\lesssim&\tau h\sum_{j=0}^{M-1}
|\xi_{j}^n|\,  \left(|\bee_{j}^{n+1}|+|\bee_{j}^n|\right)\lesssim\tau\left(\|\xi^n\|_{l^2}^2+
\|\bee^{n+1}\|_{l^2}^2+\|\bee^n\|_{l^2}^2\right)\nonumber\\
&\lesssim&\tau(\|\bee^{n+1}\|_{l^2}^2+\|\bee^n\|_{l^2}^2)
+\tau\left(\frac{h^2}{\varepsilon}+\frac{\tau^2}{\varepsilon^6}\right)^2, \qquad n\ge0.
\eea
Summing the above inequality for $n=0, 1, \ldots, m-1$,  we get
\begin{equation}
\|\bee^{m}\|_{l^2}^2-\|\bee^0\|_{l^2}^2\lesssim \tau\sum_{k=0}^{m}\|\bee^k\|_{l^2}^2+\tau m\left(\frac{h^2}{\varepsilon}+\frac{\tau^2}{\varepsilon^6}\right)^2, \qquad 0\le m\le \frac{T}{\tau}.
\end{equation}
Taking $\tau_0$ sufficiently small, when $0<\tau\le \tau_0$, we have
\begin{equation}
\|\bee^m\|_{l^2}^2\lesssim \tau\sum^{m-1}_{k=0}\|\bee^k\|_{l^2}^2+\tau m\left(\frac{h^2}{\varepsilon}+\frac{\tau^2}{\varepsilon^6}\right)^2\le
\tau\sum^{m-1}_{k=0}\|\bee^k\|_{l^2}^2+T\left(\frac{h^2}{\varepsilon}+
\frac{\tau^2}{\varepsilon^6}\right)^2,\qquad 0\leq m\leq \frac{T}{\tau}.
\end{equation}
Using the discrete Gronwall's inequality and noticing $\|\bee^0\|_{l^2}=0$, we obtain
\begin{equation}
\|\bee^m\|_{l^2}^2\lesssim T \left(\frac{h^2}{\varepsilon}+\frac{\tau^2}{\varepsilon^6}\right)^2
\lesssim \left(\frac{h^2}{\varepsilon}+\frac{\tau^2}{\varepsilon^6}\right)^2,
\qquad 0\leq m\leq \frac{T}{\tau},
\end{equation}
which immediately implies the error bound (\ref{ebcnfd97}).
\hfill $\Box$


\begin{center}
{\bf Appendix B}. Proof of Theorem \ref{thm_efd} for the LFFD method
\end{center}
\setcounter{equation}{0}
\renewcommand{\theequation}{B.\arabic{equation}}

Define the local truncation error $\tilde{\xi}^{n}=(\tilde{\xi}_0^n,\tilde{\xi}_1^n,\ldots,\tilde{\xi}_M^n)^T\in X_{M}$
of the  LFFD (\ref{efd1}) with (\ref{bdc}) and (\ref{bdct1}) as follows, for $0\leq j\leq M-1$,
\begin{align}\label{ApdexB1}
\tilde{\xi}_{j}^{n}:&=i\delta_{t}\Phi(t_n,x_{j})+
\frac{i}{\varepsilon}\sigma_1\delta_{x}\Phi(t_{n},x_{j})
-\frac{1}{\varepsilon^2}\sigma_3\Phi(t_{n},x_j)+\left[A_{1,j}^n\sigma_1-V_j^nI_2\right]\Phi(t_{n},x_j),\quad n\ge1,\\
\label{ApdexB_trun0}
\tilde{\xi}_{j}^{0}:&=i\delta_{t}^+\Phi(0,x_{j})+\frac{i}{\varepsilon}\sigma_1\delta_x\Phi_0(x_j)-\left(\frac{1}{\eps^2}\sigma_3
+V_j^0I_2-A_{1,j}^0\sigma_1\right)\Phi_0(x_j).
\end{align}
Applying the Taylor expansion in (\ref{ApdexB1}) and (\ref{ApdexB_trun0}), noticing (\ref{dir}) and the assumptions (A) and (B),
similarly to the proof of Theorem \ref{thm_cnfd},  we obtain
\begin{eqnarray}
|\tilde{\xi}_{j}^0|\lesssim \frac{\tau}{\varepsilon^4}+\frac{h^2}{\varepsilon}, \qquad |\tilde{\xi}_{j}^n|\lesssim\frac{\tau^2}{\varepsilon^6}+\frac{h^2}{\varepsilon},
\qquad j=0,1,\ldots, M-1, \quad n\ge1,
\end{eqnarray}
which immediately implies
\be\label{locbdB1}
\|\tilde{\xi}^n\|_{l^\infty}=\max_{0\le j\le M-1}|\tilde{\xi}_j^n|
\lesssim \frac{\tau^2}{\varepsilon^6}+\frac{h^2}{\varepsilon}, \qquad
\|\tilde{\xi}^n\|_{l^2}\lesssim \|\tilde{\xi}^n\|_{l^\infty} \lesssim
\frac{\tau^2}{\varepsilon^6}+\frac{h^2}{\varepsilon},
\qquad n\ge1, \qquad 0<\eps\le 1.
\ee
Subtracting (\ref{efd1}) from (\ref{ApdexB1}), noticing (\ref{err}), we get
\begin{eqnarray}
\label{trun_efdB}
i\delta_t\bee_{j}^{n}=-\frac{i}{\varepsilon}\sigma_1\delta_x\bee_{j}^{n}+\frac{1}{\varepsilon^2}\sigma_3
\bee_{j}^{n}+\left(V_j^{n}I_2-A_{1,j}^{n}\sigma_1
\right)\bee_{j}^{n}+\tilde{\xi}_j^n,\quad 0\leq j\leq M-1,\quad n\ge1,
\end{eqnarray}
where the boundary and initial conditions are given as
\be\label{ibdcfB}
\bee_0^{n}=\bee_{M}^{n}, \qquad \bee_{-1}^{n}=\bee_{M-1}^n, \qquad n\ge0,\qquad
\bee_j^0={\bf 0}, \qquad j=0,1,\ldots,M.
\ee
For the first step, we have
\be\label{eq:lftd:1st}
\|\bee^1\|_{l^2}=\tau\|\tilde{\xi}^0\|_{l^2}\lesssim \frac{\tau^2}{\varepsilon^4}+\frac{\tau h^2}{\varepsilon}\lesssim\frac{h^2}{\varepsilon}+\frac{\tau^2}{\varepsilon^6}.
\ee
Denote $\mathcal{E}^{n+1}$ for $n=0, 1, \ldots$ as
\begin{equation}
\mathcal{E}^{n+1}=\|\bee^{n+1}\|^2_{l^2}+\|\bee^{n}\|_{l^2}^2+2\,\text{Re}\left(\tau h\sum\limits_{j=0}^{M-1}(\bee_j^{n+1})^*\sigma_1\delta_x\bee^n_j\right)
-2\,\text{Im}\left(\frac{\tau h}{\eps^2}\sum\limits_{j=0}^{M-1}(\bee_j^{n+1})^*\sigma_3\bee_j^n\right);
\end{equation}
and under the stability condition (\ref{stclffd8}),  e.g., $\tau\leq \frac{\eps^2\tau_1 h}{\eps^2h V_{\rm max}+\sqrt{h^2+\eps^2(1+\eps h A_{1,{\rm max}})^2}}$ with $\tau_1=\frac14$,
 which implies $\frac{\tau}{h}\leq \frac14$ and $\frac{\tau}{\eps^2}\leq \frac14$, using Cauchy inequality, we can get that
\be\label{eq:equivs}
\frac12\left(\|\bee^{n+1}\|^2_{l^2}+\|\bee^{n}\|_{l^2}^2\right)\leq\mathcal{E}^{n+1}\leq \frac32\left(\|\bee^{n+1}\|^2_{l^2}+\|\bee^{n}\|_{l^2}^2\right), \qquad n\ge0.
\ee
It follows from (\ref{eq:lftd:1st}) that
\be\label{eq:lftdfst}
\mathcal{E}^1\lesssim \left(\frac{h^2}{\varepsilon}+\frac{\tau^2}{\varepsilon^6}\right)^2.
\ee
Similarly to the proof of Theorem \ref{thm_cnfd}, multiplying (\ref{trun_efdB}) from
the left by $2h\tau\, \left(\bee_j^{n+1}+\bee_j^{n-1}\right)^*$,  taking the imaginary part,
then summing for $j=0,1,\ldots,M-1$, using Cauchy inequality, noting
(\ref{locbdB1}) and (\ref{eq:equivs}),  we get for $n\ge1$,
\begin{align*}
\mathcal{E}^{n+1}-\mathcal{E}^{n}\lesssim &h\tau\sum_{j=0}^{M-1}\left((A_{1,\rm max}+V_{\rm max})|\bee_j^n|+|\tilde{\xi}_j^n| \right)(|\bee_j^{n+1}|+|\bee_j^{n-1}|)\\
\lesssim & \tau (\mathcal{E}^{n+1}+\mathcal{E}^{n})+\tau \left(\frac{h^2}{\varepsilon}+\frac{\tau^2}{\varepsilon^6}\right)^2, \qquad n\ge0.
\end{align*}
Summing the above inequality for $n=1, 2, \ldots, m-1$, we get
\be
\mathcal{E}^m-\mathcal{E}^1\lesssim \tau \sum_{k=1}^m\mathcal{E}^k+m\tau \left(\frac{h^2}{\varepsilon}+\frac{\tau^2}{\varepsilon^6}\right)^2, \qquad 1\le m\le \frac{T}{\tau}.
\ee
Taking $\tau_0$ sufficiently small, using the discrete Gronwall's inequality and noticing (\ref{eq:lftdfst}), we obtain from the above equation that
\begin{equation}
\mathcal{E}^m
\lesssim \left(\frac{h^2}{\varepsilon}+\frac{\tau^2}{\varepsilon^6}\right)^2,
\qquad 1\leq m\leq \frac{T}{\tau},
\end{equation}
which immediately implies the error bound (\ref{eblffd}) in view of (\ref{eq:equivs}).
\hfill $\Box$

\begin{center}
{\bf Appendix C}. Proof of Theorem \ref{thm_EWI} for the EWI-FP method
\end{center}
\setcounter{equation}{0}
\renewcommand{\theequation}{C.\arabic{equation}}

Define the error function $\bee^n(x)$ for $n=0,1,\ldots$ as
\begin{equation}
\label{pf_EWI_1}
\bee^n(x)=
        \begin{pmatrix}
         e_1^n(x)\\
         e_2^n(x) \\
        \end{pmatrix}
:=P_M\Phi(t_n,x)-\Phi_M^n(x)=\sum_{l=-M/2}^{M/2-1}\widehat{\bee}_l^n e^{i\mu_l(x-a)},\quad a\leq x\leq b.
\end{equation}
Using the triangular inequality and standard interpolation result, we get
\begin{equation}
\label{Error_L2}
\|\Phi( t_n,x)-\Phi_M^n(x)\|_{L^2}\leq\|\Phi( t_n,x)-P_M\Phi(t_n,x)\|_{L^2}+\|\bee^n(x)\|_{L^2}\leq
h^{m_0}+\|\bee^n(x)\|_{L^2}
\quad 0\leq n\leq\frac{T}{\tau},
\end{equation}
which means that we only need estimate $\|\bee^n(x)\|_{L^2}$.

Define the local truncation error $\xi^n(x)=\sum_{l=-M/2}^{M/2-1}\widehat{\xi}_l^ne^{i\mu_l(x-a)}\in Y_M$ of the
EWI-FP (\ref{Coe2}) for $n\ge0$ as
\be\label{eq:localerr}
\widehat{\xi}_l^n=\begin{cases}\widehat{(\Phi(\tau))}_l
-e^{-i\tau\Gamma_l/\eps^2}\widehat{(\Phi(0))}_l
+i\eps^2\Gamma_l^{-1}\left[I_2-e^{-\frac{i\tau}{\eps^2}\Gamma_l}\right]
\widehat{(G(0)\Phi(0))}_l, &n=0,\\
\widehat{(\Phi(t_{n+1}))}_l-
e^{-i\tau\Gamma_l/\eps^2}\widehat{(\Phi(t_n))}_l+iQ_l^{(1)}(\tau)\widehat{(G(t_n)\Phi(t_n))}_l
+iQ_l^{(2)}(\tau)\delta_t^-\widehat{\left(G(t_n)\Phi(t_n)\right)}_l, &n\ge1,
\end{cases}
\ee
where we write $\Phi(t)$ and $G(t)$ in short for $\Phi(t,x)$ and $G(t,x)$, respectively.

Firstly, we estimate the local truncation error $\xi^n(x)$.
Multiplying both sides of the  Dirac equation (\ref{dir})
by $e^{i\mu_l(x-a)}$ and integrating over the interval $(a,b)$,
we easily recover the equations for $(\widehat{\Phi(t)})_l$,
which are exactly the same as (\ref{ODE765}) with $\Phi_M$ being
replaced by $\Phi(t,x)$.  Replacing $\Phi_M$ with $\Phi(t,x)$, we
use the same notations $\widehat{F}_l^n(s)$  as in (\ref{eq:fdef}) and
 the time derivatives of  $\widehat{F}_l^n(s)$
enjoy the same properties  of time derivatives of $\Phi(t,x)$.
Thus, the same representation (\ref{EWIn})
holds for $(\widehat{\Phi(t_n)})_l$ with $n\ge1$.
From the derivation  of the EWI method,
it is clear that the error $\xi^n(x)$ comes from the approximations
for the integrals in (\ref{eq:inteapp0}) and (\ref{eq:inteappn}), and we have
\begin{align}
\label{EWI_Tr_1:1}
\widehat{\xi}^0_l=&-i\int_0^{\tau}e^{\frac{i(s-\tau)}{\varepsilon^2}\Gamma_l}(\widehat{F}^0_l(s)
-\widehat{F}^0_l(0))ds=-i\int_0^{\tau}\int_0^se^{\frac{i(s-\tau)}{\varepsilon^2}\Gamma_l}
\partial_{s_1}\widehat{F}^0_l(s_1)\,ds_1ds,
\end{align}
and for $n\ge1$
\begin{align}
\widehat{\xi}^n_l=&-i\int_0^{\tau}e^{\frac{i(s-\tau)}{\varepsilon^2}\Gamma_l}
\left(\int_0^s\int_0^{s_1}\partial_{s_2s_2}\widehat{F}^{n}_l(s_2)\,ds_2ds_1
+s\int_0^1\int_{\theta\tau}^\tau\partial_{\theta_1\theta_1}\widehat{F}_l^{n-1}(\theta_1)\,d
\theta_1d\theta\right)ds.\label{EWI_Tr_1:2}
\end{align}
For $n=0$, the above equalities imply $|\widehat{\xi}^0_l|\lesssim\int_0^{\tau}\int_0^s|\partial_{s_1}\widehat{F}^0_l(s_1)|ds_1ds$ and
 by the Bessel inequality  and assumptions (C) and (D), we find
\begin{align*}
\|\xi^0(x)\|_{L^2}^2=&(b-a)\sum\limits_{l=-M/2}^{M/2-1}|\widehat{\xi}_l^0|^2
\lesssim (b-a)\tau^2\int_0^\tau\int_0^s\sum\limits_{l=-M/2}^{M/2-1}|\partial_{s_1} \widehat{F}^0_l(s_1)|^2\,ds_1ds\\
\lesssim& \tau^2\int_0^\tau\int_0^s\|\partial_{s_1}(G(s_1)\Phi(s_1))\|_{L^2}^2\,ds_1ds
\lesssim \frac{\tau^4}{\varepsilon^4}.
\end{align*}
Similarly, for $n\ge1$,  we obtain
\begin{align*}
&\|\xi^n(x)\|_{L^2}^2=(b-a)\sum\limits_{l=-M/2}^{M/2-1}|\widehat{\xi}_l^n|^2\\
&\lesssim \tau^3\int_0^{\tau}\int_0^s\int_{0}^{s_1}
\sum\limits_{l=-\frac{M}{2}}^{\frac{M}{2}-1}|\partial_{s_2s_2}\widehat{F}_l^n(s_2)|^2\,ds_2ds_1ds
+\tau^3\int_0^\tau\int_{0}^1\int_{\theta\tau}^\tau s\sum\limits_{l=-\frac{M}{2}}^{\frac{M}{2}-1}
|\partial_{\theta_1\theta_1}\widehat{F}_l^{n-1}(\theta_1)|^2\,d\theta_1\,d\theta\,ds\\
&\lesssim\tau^6\|\partial_{tt}(G(t)\Phi(t))\|_{L^\infty([0,T]; (L^2)^2)}^2
\lesssim\frac{\tau^6}{\varepsilon^8},
\end{align*}
where we have used the assumptions (C) and (D). Hence, we derive that
\be\label{eq:localerrewi}
\|\xi^0(x)\|_{L^2}\lesssim \frac{\tau^2}{\varepsilon^2},\quad \|\xi^n(x)\|_{L^2}\lesssim \frac{\tau^3}{\varepsilon^4},\quad n\ge1.
\ee

Now, we look at the error equations. For each fixed
$l=-M/2,...,M/2-1$, subtracting (\ref{Coe2}) from (\ref{eq:localerr}),
we obtain the equation for the error vector function  as
\begin{equation}
\label{Error_fun_EWI}
\widehat{\bee}^0_l={\bf 0},\qquad\widehat{\bee}^1_l=\widehat{\xi}^0_l; \qquad \widehat{\bee}^{n+1}_l=
e^{-i\tau\Gamma_l/\eps^2}\widehat{\bee}^n_l+\widehat{R}^n_l
+\widehat{\xi}^n_l,\quad 1\leq n\leq\frac{T}{\tau}-1,
\end{equation}
where $R^n(x)=\sum\limits_{l=-M/2}^{M/2-1}\widehat{R}_l^ne^{i\mu_l(x-a)}\in Y_M$ for $n\ge1$ is given by
\be\label{eq:R1}
\widehat{R}^n_l=-iQ_l^{(1)}(\tau)\left(\widehat{(G(t_n)\Phi(t_n))}_l-\widehat{(G(t_n)\Phi^n_M)}_l\right)
-iQ_l^{(2)}(\tau)\left(\delta_t^-\widehat{\left(G(t_n)\Phi(t_n)\right)}_l-\delta_t^-\widehat{(G(t_n)\Phi^n_M)}_l\right).
\ee

Using the properties of the matrices $Q_l^{(1)}(\tau)$ and $Q_l^{(2)}(\tau)$, it is easy to verify that
\be\label{eq:Qbd}
\|Q_l^{(1)}(\tau)\|_2\leq \tau,\quad \|Q_l^{(2)}(\tau)\|_2\leq \frac{\tau^2}{2},\quad l=-\frac{M}{2},\ldots,\frac{M}{2}-1,
\ee
where $\|Q\|_2$ denotes the $l^2$ norm of the matrix $Q$. Combining (\ref{eq:Qbd}), (\ref{eq:R1})
and the assumption (D), we get
\begin{align}
\|R^n(x)\|_{L^2}^2=&(b-a)\sum\limits_{l=-M/2}^{M/2-1}|\widehat{R}_l^n|^2
\lesssim (b-a)\tau^2\sum\limits_{l=-M/2}^{M/2-1}\left(\left|\widehat{(\Phi(t_n))}_l-(\widehat{\Phi_M^n})_l\right|^2
+\left|\widehat{(\Phi(t_{n-1}))}_l-(\widehat{\Phi_M^{n-1}})_l\right|^2\right)\nonumber\\
\lesssim&\tau^2\sum_{k=n-1}^n\|\Phi( t_k,x)-\Phi_M^k(x)\|_{L^2}^2\lesssim \tau^2 h^{2m_0}
+\tau^2\|\bee^{n}(x)\|_{L^2}^2+\tau^2\|\bee^{n-1}(x)\|_{L^2}^2.\label{eq:rnbd}
\end{align}

Multiplying both sides of (\ref{Error_fun_EWI}) by $\left(\widehat{\bee}^{n+1}_l+e^{-i\tau\Gamma_l/\eps^2}\widehat{\bee}^n_l\right)^*$ from left, taking
the real parts and using Cauchy inequality, we obtain
\be
\left|\widehat{\bee}^{n+1}_l\right|^2-\left|\widehat{\bee}^{n}_l\right|^2
\leq
\tau\left[\left|\widehat{\bee}^{n+1}_l\right|^2+\left|\widehat{\bee}^{n}_l\right|^2\right]
+\frac{|\widehat{R}_l^n|^2}{\tau}+\frac{|\widehat{\xi}_l^n|^2}{\tau}.
\ee
Multiplying the above inequality by $b-a$ and summing together for $l=-M/2, \ldots , M/2-1$, in view of the Bessel inequality,  we obtain
\begin{align}\label{eq:rec1}
\left\|\bee^{n+1}(x)\right\|_{L^2}^2-\left\|\bee^{n}(x)\right\|_{L^2}^2\lesssim & \tau(\left\|\bee^{n+1}(x)\right\|_{L^2}^2+\left\|\bee^{n}(x)\right\|_{L^2}^2)+
\frac{1}{\tau}\|R^n(x)\|^2_{L^2}+\frac{1}{\tau}\|\xi^n(x)\|^2_{L^2},\quad n\ge1.
\end{align}
Summing (\ref{eq:rec1}) for $n=1,\ldots,m-1$, using (\ref{eq:rnbd}) and (\ref{eq:localerrewi}), we derive
\be
\left\|\bee^{m}(x)\right\|_{L^2}^2-\left\|\bee^{1}(x)\right\|_{L^2}^2\lesssim \tau\sum\limits_{k=1}^{m}\left\|\bee^{k}(x)\right\|_{L^2}^2+\frac{m\tau^5}{\varepsilon^8}+m\tau h^{2m_0},\quad 1\le m\leq\frac{T}{\tau}.
\ee
Since $\|\bee^0(x)\|_{L^2}=0$ and $\|\bee^1(x)\|_{L^2}\lesssim\frac{\tau^2}{\varepsilon^2}\lesssim\frac{\tau^2}{\varepsilon^4}$,
the discrete Gronwall's inequality will imply that for sufficiently small $\tau$,
\be\label{eq:eq:l2e}
\left\|\bee^{m}(x)\right\|_{L^2}^2\lesssim  h^{2m_0}+\frac{\tau^4}{\varepsilon^8},\qquad 1\le m\leq \frac{T}{\tau}.
\ee
Combining (\ref{Error_L2}) and (\ref{eq:eq:l2e}), we draw the conclusion (\ref{thm_eq_EWI}). \hfill $ \square$

\begin{center}
{\bf Appendix D}. Extensions of the EWI-FS  (\ref{AL1})-(\ref{Coe2}) and  TSFP (\ref{eq:tsfp}) in 2D and 3D
\end{center}
\setcounter{equation}{0}
\renewcommand{\theequation}{D.\arabic{equation}}
The EWI-FS (\ref{AL1})-(\ref{Coe2}),  EWI-FP (\ref{AL4})-(\ref{AL5}) and TSFP
(\ref{eq:tsfp}) can be easily extended to 2D and 3D with tensor  grids by modifying
the matrices $\Gamma_l$ in (\ref{eq:Gamma}) and $G(t,x)$ in (\ref{eq:Gjn}) in the TSFP case.
For the reader's convenience, we present the modifications of
 $\Gamma_l$ in (\ref{eq:Gamma}) and $G(t,x)$ in (\ref{eq:Gjn}) in 2D and 3D as follows.

{\it For the Dirac equation (\ref{SDEdd}) in 2D}, i.e. we take $d=2$ in (\ref{SDEdd}). The problem
is truncated on $\Omega=(a_1, b_1)\times(a_2, b_2)$ with mesh sizes $h_1=(b_1-a_1)/M_1$ and $h_2=(b_2-a_2)/M_2$ ($M_1,M_2$ two even positive integers) in the $x$- and $y$-direction, respectively.
The wave function $\Phi$ is a two-component vector, and the matrix $\Gamma_l$ in (\ref{eq:Gamma}) will be
replaced by
\begin{equation}
\Gamma_{jk}=\begin{pmatrix}
                   1 & \varepsilon \mu_j^{(1)}-i\eps\mu_k^{(2)} \\
                   \eps\mu_j^{(1)}+i\eps\mu_k^{(2)} & -1\\
                \end{pmatrix},\qquad \mu_j^{(1)}=\frac{2j\pi }{b_1-a_1},\quad \mu_k^{(2)}=\frac{2k\pi }{b_2-a_2},
\end{equation}
where $-\frac{M_1}{2}\leq j\leq \frac{M_1}{2}-1$, $-\frac{M_2}{2}\leq k\leq \frac{M_2}{2}-1$,
and the Schur decomposition $\Gamma_{jk}=Q_{jk}D_{jk}Q_{jk}^*$ is given as
\begin{equation}
 Q_{jk}=\begin{pmatrix}
\frac{1+\delta_{jk}}{\sqrt{2\delta_{jk}(1+\delta_{jk})}} &\frac{-\eps\mu_{j}^{(1)}+i\eps\mu_k^{(2)}}{\sqrt{2\delta_{jk}(1+\delta_{jk})}}\\
\frac{\eps\mu_{j}^{(1)}+i\eps\mu_k^{(2)}}{\sqrt{2\delta_{jk}(1+\delta_{jk})}} &\frac{1+\delta_{jk}}{\sqrt{2\delta_{jk}(1+\delta_{jk})}}
\end{pmatrix}, \quad D_{jk}=\begin{pmatrix}
\delta_{jk} &0\\
0 &-\delta_{jk}\\
\end{pmatrix},\quad \delta_{jk}=\sqrt{1+\eps^2(\mu_{j}^{(1)})^2+\eps^2(\mu_k^{(2)})^2}.
\end{equation}
The matrix $\int_{t_n}^{t_{n+1}}G(t,x)dt$ in (\ref{eq:Gjn}) becomes $\int_{t_n}^{t_{n+1}}G(t,\bx)dt$ and the Schur decomposition $\int_{t_n}^{t_{n+1}}G(t,\bx)dt=P_\bx\Lambda_\bx P_{\bx}^*$ with
$V_\bx^{(1)}=\int_{t_n}^{t_{n+1}}V(t,\bx)dt$,
$A_{l,\bx}^{(1)}=\int_{t_n}^{t_{n+1}}A_l(t,\bx)dt$ for $l=1,2$,
$\lambda_\bx^{(1)}=\sqrt{|A_{1,\bx}^{(1)}|^2+|A_{2,\bx}^{(1)}|^2}$,
$\Lambda_\bx={\rm diag}(\Lambda_{\bx,-},\Lambda_{\bx,+})$,
$\Lambda_{\bx,\pm}=V_\bx^{(1)}\pm \lambda_\bx^{(1)}$, and  $P_\bx=I_2$ if $\lambda_\bx^{(1)}=0$ and otherwise
\be
 P_\bx=\begin{pmatrix}\frac{1}{\sqrt{2}}&\frac{A_{1,\bx}^{(1)}-iA_{2,\bx}^{(1)}}{\sqrt{2}\lambda_\bx^{(1)}}\\
 \frac{A_{1,\bx}^{(1)}+iA_{2,\bx}^{(1)}}{\sqrt{2}}&\frac{1}{\sqrt{2}\lambda_\bx^{(1)}}
 \end{pmatrix}.
\ee

{\it For the Dirac equation (\ref{SDEd}) in 3D}, i.e. we take $d=3$ in (\ref{SDEd}).
The problem is truncated on $\Omega=(a_1, b_1)\times(a_2, b_2)\times(a_3, b_3)$ with mesh sizes $h_1=(b_1-a_1)/M_1$, $h_2=(b_2-a_2)/M_2$ and $h_3=(b_3-a_3)/M_3$ ($M_1,M_2,M_3$ three even positive integers) in
 $x$-, $y$- and $z$-direction, respectively.
The wave function $\Psi$ is a four-component vector, and the matrix $\Gamma_l$ in (\ref{eq:Gamma})
will be replaced by  $\Gamma_{jkl}$  as:
\begin{equation}
\Gamma_{jkl}=\begin{pmatrix}
                   1 & 0 & \eps\mu_l^{(3)} & \eps\mu_j^{(1)}-i\eps\mu_k^{(2)} \\
                   0 & 1 & \eps\mu_j^{(1)}+i\eps\mu_k^{(2)} & -\eps\mu_l^{(3)} \\
                   \eps\mu_l^{(3)} & \eps\mu_j^{(1)}-i\eps\mu_k^{(2)} & -1 & 0 \\
                   \eps\mu_j^{(1)}+i\eps\mu_k^{(2)} & -\eps\mu_l^{(3)} & 0 & -1 \\
                \end{pmatrix},
\end{equation}
where $-\frac{M_1}{2}\leq j\leq\frac{M_1}{2}-1,-\frac{M_2}{2}\leq k\leq\frac{M_2}{2}-1,-\frac{M_3}{2}\leq l\leq\frac{M_3}{2}-1$  and
\begin{equation}
\mu_j^{(1)}=\frac{2j\pi }{b_1-a_1},\quad \mu_k^{(2)}=\frac{2k\pi }{b_2-a_2},\quad \mu_l^{(3)}=\frac{2l\pi }{b_3-a_3}.
\end{equation}
The eigenvalues of $\Gamma_{jkl}$ are
$$
\delta_{jkl}, \delta_{jkl}, -\delta_{jkl}, -\delta_{jkl},\quad \text{with}\quad \delta_{jkl}=\sqrt{1+\varepsilon^2\left|\mu_j^{(1)}\right|^2+\eps^2\left|\mu_{k}^{(2)}\right|^2+
\eps^2\left|\mu_l^{(3)}\right|^2}.
$$
 The corresponding eigenvectors are
 \begin{equation*}
 \bold{v}^{(1)}_{jkl}=\begin{pmatrix}
1+\delta_{jkl}\\
 0\\
 \eps\mu_l^{(3)}\\
 \eps\mu_j^{(1)}+i\eps\mu_{k}^{(2)}\end{pmatrix},\,
 \bold{v}^{(2)}_{jkl}=\begin{pmatrix}
 0\\
 1+\delta_{jkl}\\
 \eps\mu_j^{(1)}-i\eps\mu_{k}^{(2)}\\
 -\eps\mu_l^{(3)}\end{pmatrix},\,
 \bold{v}^{(3)}_{jkl}=\begin{pmatrix}
 -\eps\mu_l^{(3)}\\
 -\eps\mu_j^{(1)}-i\eps\mu_{k}^{(2)}\\
 1+\delta_{jkl}\\
 0\end{pmatrix},\,
 \bold{v}^{(4)}_{jkl}=\begin{pmatrix}
 -\eps\mu_j^{(1)}+i\eps\mu_{k}^{(2)}\\
 \eps\mu_l^{(3)}\\
 0\\
1+\delta_{jkl}\end{pmatrix}.
 \end{equation*}
Then the Schur decomposition $\Gamma_{jkl}=Q_{jkl}D_{jkl}Q^*_{jkl}$ is given as
 \begin{equation*}
 D_{jkl}=\mathrm{diag}(\delta_{jkl},\delta_{jkl},-\delta_{jkl},-\delta_{jkl}),\quad
 Q_{jkl}=\frac{1}{\sqrt{2\delta_{jkl}(1+\delta_{jkl})}}\left(\bold{v}^{(1)}_{jkl},
 \bold{v}^{(2)}_{jkl},\bold{v}^{(3)}_{jkl},\bold{v}^{(4)}_{jkl}\right).
 \end{equation*}
The matrix $\int_{t_n}^{t_{n+1}}G(t,x)dt$ in (\ref{eq:Gjn}) becomes $\int_{t_n}^{t_{n+1}}G(t,\bx)dt$ and the Schur decomposition $\int_{t_n}^{t_{n+1}}G(t,\bx)dt=P_\bx\Lambda_\bx P_{\bx}^*$ with
$V_\bx^{(1)}=\int_{t_n}^{t_{n+1}}V(t,\bx)dt$,
$A_{l,\bx}^{(1)}=\int_{t_n}^{t_{n+1}}A_l(t,\bx)dt$ for $l=1,2,3$,
$\lambda_\bx^{(1)}=\sqrt{|A_{1,\bx}^{(1)}|^2+|A_{2,\bx}^{(1)}|^2+|A_{3,\bx}^{(1)}|^2}$,
$\Lambda_\bx={\rm diag}(\Lambda_{\bx,-},\Lambda_{\bx,-},\Lambda_{\bx,+},\Lambda_{\bx,+})$,
$\Lambda_{\bx,\pm}=V_\bx^{(1)}\pm \lambda_\bx^{(1)}$, and  $P_\bx=I_4$ if $\lambda_\bx^{(1)}=0$ and otherwise
\begin{equation*}
{\bf u}_\bx^{(1)}=\begin{pmatrix}\frac{1}{\sqrt{2}}\\0\\ \frac{A_{3,\bx}^{(1)}}{\sqrt{2}\lambda_\bx^{(1)}}\\ \frac{A_{1,\bx}^{(1)}+iA_{2,\bx}^{(1)}}{\sqrt{2}\lambda_\bx^{(1)}}\end{pmatrix},
\quad {\bf u}^{(2)}=\begin{pmatrix}0\\ \frac{1}{\sqrt{2}}\\ \frac{A_{1,\bx}^{(1)}-iA_{2,\bx}^{(1)}}{\sqrt{2}\lambda_\bx^{(1)}}\\ \frac{-A_{3,\bx}^{(1)}}{\sqrt{2}\lambda_\bx^{(1)}}\end{pmatrix},\quad
{\bf u}^{(3)}=\begin{pmatrix}\frac{-A_{3,\bx}^{(1)}}{\sqrt{2}\lambda_\bx^{(1)}}\\  \frac{-A_{1,\bx}^{(1)}-iA_{2,\bx}^{(1)}}{\sqrt{2}\lambda_\bx^{(1)}}\\
\frac{1}{\sqrt{2}}\\ 0\end{pmatrix},\quad
{\bf u}^{(4)}=\begin{pmatrix}\frac{-A_{1,\bx}^{(1)}+iA_{2,\bx}^{(1)}}{\sqrt{2}\lambda_\bx^{(1)}}\\  \frac{A_{3,\bx}^{(1)}}{\sqrt{2}\lambda_\bx^{(1)}}\\0\\ \frac{1}{\sqrt{2}}\end{pmatrix}.
\end{equation*}

For the Dirac equation (\ref{SDEd}) in 2D, we simply let $\mu_{l}^{(3)}=0$, $A_3(t,\bx)\equiv0$
in the above  3D case; and for the Dirac equation (\ref{SDEd}) in 1D, we
let $\mu_{k}^{(2)}=\mu_{l}^{(3)}=0$, $A_2(t,\bx)=A_3(t,\bx)\equiv 0$ in the above  3D case. Then the EWI-FP (\ref{AL4})-(\ref{AL5}) and  TSFP (\ref{eq:tsfp}) can be designed
accordingly for the  Dirac equation (\ref{SDEd}) in 2D and 1D.

\section*{Acknowledgements}
W. Bao and X. Jiao acknowledge support from  the Ministry of Education of Singapore grant R-146-000-196-112.
Y. Cai was partially  supported by NSF grants DMS-1217066 and DMS-1419053.  Q. Tang acknowledge
the support from the  ANR project BECASIM ANR-12-MONU-0007-02.
Part of this work was done when the authors were visiting
the Institute for Mathematical Sciences at the National University of Singapore in 2015.



\end{document}